\documentclass[10pt]{amsart}
\usepackage[utf8]{inputenc}
\usepackage{enumerate}
 \pdfoutput=1

\usepackage{enumerate}
\usepackage{amsmath, amsthm}
\usepackage{amsfonts}
\usepackage{amssymb}
\usepackage{graphicx}
\usepackage{url}

\usepackage[left=3cm, right=3cm]{geometry}
\usepackage{array}
\makeatletter
\newcolumntype{e}[1]{
   >{\minipage[t]{\linewidth}%
     \NoHyper
     \let\\\tabularnewline
     \enumerate
        \addtolength{\rightskip}{0pt plus 50pt}
        \setlength{\itemsep}{-\parsep}}%
   p{#1}%
   <{\@finalstrut\@arstrutbox\endenumerate
     \endNoHyper
     \endminipage}}

\newcolumntype{i}[1]{
   >{\minipage[t]{\linewidth}%
        \let\\\tabularnewline
        \itemize
           \addtolength{\rightskip}{0pt plus 50pt}%
           \setlength{\itemsep}{-\parsep}}%
   p{#1}%
   <{\@finalstrut\@arstrutbox\enditemize\endminipage}}

\AtBeginDocument{%
    \@ifpackageloaded{hyperref}{}%
        {\let\NoHyper\relax\let\endNoHyper\relax}}
\makeatother

\AtEndDocument{\bigskip{\footnotesize%
  \textsc{Westfälische Wilhelms-Universität Münster, Institut für Mathematische Logik und Grundlagenforschung
  	Einsteinstraße 62, 48149 M\"{u}nster, Deutschland
} \par  
  \textit{E-mail address}, G\"{o}nen\c{c} Onay  \texttt{onay@uni-muenster.de} }}
\renewcommand{\epsilon}{\varepsilon}
\renewcommand{\le}{\leqslant}
\renewcommand{\ge}{\geqslant}
\renewcommand{\leq}{\leqslant}
\renewcommand{\geq}{\geqslant}

\newcommand{\jump}[1]{\operatorname{Jump}(#1)}
\newcommand{\lcm}{\operatorname{lcm}}
\newcommand{\hens}{\operatorname{hens}}
\newcommand{\h}{\operatorname{h}}
\newcommand{\rv}[1]{\operatorname{rv}(#1)}
\newcommand{\nqrt}[2]{\sqrt[\varphi^{#1}]{#2}}
\newcommand{\ssi}{\Leftrightarrow}
\newcommand{\implique}{\Rightarrow}
\newcommand{\vers}{\rightarrow}
\newcommand{\tq}{\;\: | \:\:}
\newcommand{\mesp}{\;\:}
\newcommand {\vfi}{\varphi}
\newcommand{\mr}[2]{\text{M}_{R}(#1,#2)}
\newcommand{\iffl}{\longleftrightarrow}
\newcommand{\iffs}{\leftrightarrow}

\newcommand{\Nn}{\mathbb{N}}
\newcommand{\Zz}{\mathbb{Z}}
\newcommand{\Qq}{\mathbb{Q}}
\newcommand{\Rr}{\mathbb{R}}
\newcommand{\Cc}{\mathbb{C}}
\newcommand{\Kk}{\mathbb{K}}
\newcommand{\Hh}{\mathbb{H}}
\newcommand{\F}{\mathbb{F}}
\newcommand{\K}{\mathbb{K}}

\newtheorem{theorem}{Theorem}[section]
\newtheorem{lemma}[theorem]{Lemma}
\newtheorem{prop}[theorem]{Proposition}
\newtheorem{proposition}[theorem]{Proposition}
\newtheorem{corr}[theorem]{Corollary}

\newtheorem{defn}[theorem]{Definition}
\newtheorem{fait}[theorem]{Fact}

\newtheorem{remark}[theorem]{Remark}
\theoremstyle{definition}

\newtheorem{obs}[theorem]{Observation}
\newtheorem{nota}[theorem]{Notation}

\author{G\"{o}nen\c{c} Onay}   
\thanks{Research partially founded by DFG through SFB 878.}

\makeatletter
\def\blfootnote{\gdef\@thefnmark{}\@footnotetext}
\makeatother
\title{$\mathbb{F}_p((X))$ is decidable as a module over the ring of additive polynomials
}

\begin{document}

\sloppy

\begin{abstract}
Let $p$ be a prime number, $K$ be the henselization of the rational functions over the finite
field $\mathbb{F}_p$ and $R$ be the ring of additive polynomials over K. We show that the field of Laurent series
over $\mathbb{F}_p$ is decidable seen as an R-module. Moreover, we provide a recursively enumerable axiom system
(satisfied by $K$) in the language of $R$-modules together with a unary predicate for the valuation ring,
modulo which every positive primitive formula is equivalent to a universal formula. Consequently the $R$-module
theory of the field of Laurent series is model-complete in this language and admits $K$ as its prime
model.	
\end{abstract}
\maketitle
\blfootnote{\textup{2010} \textit{Mathematics Subject Classification}: 11U05 (primary), 12L05, 12J10	03C68 (secondary)}

\section{Introduction}

The decidability and the axiomatization of the field of Laurent
series over the finite field $\mathbb{F}_p$ is a longstanding 
problem (cf. for instance \cite{ku-elementary} and \cite{kuhlmanntame}). On the other hand, its characteristic $0$ analog $\mathbb{Q}_p$ is  axiomatized and is decidable by the Ax-Kochen and Ershov Theorem, proved independently by Ax-Kochen in USA and by Ershov in Russia (see \cite{ax-kochen1} and \cite{ersov}) . As a consequence, $\mathbb{Q}_p$ and $\mathbb{F}_p((X))$ have elementary  ultraproducts (over all primes $p$).  
The first and maybe the most famous application of Ax-Kochen and Ershov theorem (called also Ax-Kochen theorem after  \cite{ax-kochen1}) is to give a  corrected form of a conjecture by Artin. This corrected form states that 
given any integer $d>0$ there is some integer 
$m$, such that for any prime $p>m$, every homogeneous polynomial of degree $d$ over $\mathbb{Q}_p$ with $>d^2$ variables has a non-trivial zero in $\mathbb{Q}_p$.
This application uses simply that fact that for every prime $p$ the field $\mathbb{F}_p((X))$  is $C_2$ 
(i.e. every homogeneous polynomial with strictly more
        variables than the square of its degree has a non-trivial zero)\footnote{The original conjecture of Artin was claiming that $\mathbb{Q}_p$ was $C_2$ for all prime $p$, this is refuted by Ternejian.}.

Despite the (asymptotic -) analogy between $\mathbb{Q}_p$ and $\mathbb{F}_p((X))$, for a fixed prime $p$ 
the field theory of $\mathbb{F}_p((X))$ remains  unknown. Moreover,  Kuhlmann proved that the naive translation into positive characteristic of the 
complete theory of $\mathbb{Q}_p$ is incomplete (see \cite{ku-elementary}, we will give more details on this in the following lines).
Denef and Schoutens showed that the existential theory of $\mathbb{F}_p((X))$, in the language of rings with a constant symbol for $X$, is decidable assuming  resolution of singularities in positive characteristic (see \cite{denefsch}).  Recently Anscombe and Fehm showed unconditionally that the existential theory of $\mathbb{F}_p((X))$, in the language of rings (without the parameter $X$), is decidable \cite{sylvie-fehm}.

After the works of Kuhlmann and van den Dries, incompleteness of this naive theory  can be expressed 
using only  properties of   additive polynomials (as we explain in the following lines), that is, in  the language of $S$-modules  where $S$ 
is the ring of  additive polynomials over $\mathbb{F}_p(X)$.
Rohwer, in his thesis (see \cite{rohwer}), shows that the complete theory of $\mathbb{F}_p((X))$ as an $S$-module is model-complete 
in the language of $S$-modules together with a predicate for the valuation ring\footnote{However the meaning of Conjecture A.4. in \cite{rohwer} remains unclear for us.}. 
However, he does not provide an axiom system for this theory and does not show a decidability result.

In this article we are considering the same problem (up to passing to the ring of additive polynomials
over the henselization of 
$F_p(X))$ and we think that our philosophy is similar to Rohwer’s one.  On the other hand, we have had been rather inspired by  the articles \cite{ddp1}, \cite{vdd-kuhlmann} and \cite{sylvie-fehm}.

Notice that  B\'elair and Point also studied the module theory of
some (valued) fields satisfying  strong divisibility conditions (see
\cite{belair-point10} and \cite{belair-point15}).

\subsection*{Main results}
Let $K$ be the henselization of the field of rational functions over $\mathbb{F}_p$ and  $$\vfi:K\to K,\mesp x\mapsto x^{p^k}$$ 
    be some (fixed-) power of the Frobenius map. We set $R$ to be the ring of $\vfi$-polynomials, that is, additive polynomials 
    whose monomials are of the form $aT^{{p^k}^i}$ ($a \in K$), equipped with the
    composition and the usual addition. Let $M:=\mathbb{F}_p((X))$ be the field of Laurent series over $\mathbb{F}_p$, let $L$ be the language of (right-) $R$-modules
    and finally let $L_{\mathcal{O}}$ be the language $L$ together with the unary predicate $\mathcal{O}$ (for the valuation ring $\mathbb{F}[[X]]$).

    We prove in particular the following results in the present article.

\begin{enumerate}
\item There is a recursively enumerable $L_{\mathcal{O}}$-theory $T_1$ such that $M\models T_1$, and any completion of $T_1$ is model complete. Moreover,  $K$ is the prime model of the complete $L_{\mathcal{O}}$-theory of $M$  (see Theorem \ref{th:pp2universal} and following corollaries). 

\item Both the  complete $L$- and $L_{\mathcal{O}}$-theories of $M$ are decidable (see Theorem \ref{thm:finaldecidable}).
\end{enumerate}

\subsection*{Presentation of the problem and our strategy}
Following Kuhlmann (\cite{ku-elementary}), let us first explain why the naive adaptation of the theory
of $p$-adics into positive characteristic is incomplete. We denote this theory by $T_{naive}$, written in 
language of rings together with a unary predicate for the valuation ring 
and a constant for the uniformizer $X$,

We recall that $(U,v) \models T_{naive}$ if and only if
\begin{enumerate}[-]\label{naivfpt}
\item $(U,v)$ is a  henselian, non trivially valued  field of characteristic $p>0$,
\item the residue field $U/v$ is $\mathbb{F}_p$,
\item the value group $\Gamma$ is a $\mathbb{Z}$-group, 
\item $(U,v)$ is defectless,
\item $v(X)=\min \Gamma_{>0}$.
\end{enumerate}

Set $M=\mathbb{F}_p((X))$. We know that:
\begin{equation}\label{as-formula}
M=\bigoplus_{i=1}^{p-1}M^pX^i + \wp(M) + \mathcal{O}_M, \tag{*} 
\end{equation}
where $\wp$ is the Artin-Scherier map $x\mapsto x^p -x$  and $\mathcal{O}_M$ is the valuation ring of $M$ (see \cite{ku-elementary}, Lemma 1.2).
But the equality (\ref{as-formula}) fails to hold  for some extension  $N\supset M$, which is a  model of  $T_{naive}$. Hence $T_{naive}$ is incomplete.

Consider the polynomial $F(z_0,\ldots,z_{p-1}):=z_0^p-z_0 +\sum_{i=1}^{p-1} z_i^pX^i$. From the equality above  and using Hensel's Lemma, one can deduce that for any $a$, the set 
$\{v(F(x)-a) \tq x \in M\}$ has a maximum in $\mathbb{Z}\cup \{\infty\}$: we say that the
image of $F$ has the {\it optimal approximation property} (see \cite{vdd-kuhlmann}). From this observation Kuhlmann suggests a candidate for a complete axiomatization
of the theory of $M$ in the language of rings,  which is essentially $T_{naive}$ together with sentences
saying that the image of every multi-variable additive polynomial has the optimal approximation property (see \cite{ku-elementary} and \cite{kuhlmanntame}).

Now we present our approach. We work in the language 
$L_{\mathcal{O}}$. 

Let us illustrate in an example the main ideas of the proof of Theorem \ref{th:pp2universal}, which states that modulo the theory $T_1$
every {\it positive primitive} formula is equivalent to a universal one. This theorem immediately shows that every completion of
$T_1$ is model-complete (see Section \ref{sec:baur-monk} and Corollary \ref{corr:baur-monk}).

Instead of considering the multivariate polynomial $F$ above, we use the equation (\ref{as-formula}) to study the image of the Artin-Scherier map $\wp$ and see the sum 
$$C(M):=\sum_{i=1}^{p-1}M^pX^i$$ as {\bf the pseudo-complement}  to the image 
$\wp(M)$. Note  that $$\wp(M)\cap C(M) \subseteq \mathcal{O}_M,$$ hence this intersection is {\it small} 
with respect to the valuation metric; that is why we call $C$ the pseudo-complement. In addition, using Hensel's lemma we have that 
$$\mathcal{O}_M \cap \wp(M)=\mathfrak{m}$$ where 
$\mathfrak{m}$ is the maximal ideal $X\mathbb{F}_p[[X]]$.

At this point, we are able to describe the set $\wp(M)$, definable a priori by an existential formula, by the following formula:
\begin{equation}\label{elim-as}
 \phi(x): \forall y  \mesp ((\exists z \mesp \wp(z)=x-y\wedge y\in C + \mathcal{O}) \rightarrow y \in \mathfrak{m}).   \tag{**}
\end{equation}
 It is easy to see that this formula is equivalent to a universal 
 formula in the language $L_\mathcal{O}$ (hence also in the language of rings with the parameter $X$) since
 $C$ is existentially definable.

Motivated by the above example, we consider the general case. We introduce the notion of a ball (in this section) as an $L_\mathcal{O}$-formula of the form 
$$B(x): x.X^\gamma \in \mathcal{O}$$ for some fixed $\gamma \in \mathbb{Z}$. Let $Q$ be a  finite set of  polynomials (scalars from $R$).  Abusively, we also denote by $Q$ the formula which defines the sum of the images of $q\in Q$.  Our strategy is to assign to $Q$ a positive existential  formula $D$ and balls $B_1,B_2$ in a computable way,  such that in every model $N\models T_1$,
    \begin{enumerate}[(i)]
    \item \label{task1} $D(N)  + Q(N)=N$ and $Q(N)\cap D(N) \subseteq B_2(N)$
    \item \label{task2} $B_1(N) \subseteq  B_2(N)$, such that  $B_1(N)\cap Q(N)$ is definable by a universal formula which depends only on $T_1$. 
    \end{enumerate}
 For instance, the above example suggests that for $Q=\{\wp\}$, 
 $B_1=\mathcal{O}$, $B_2=\mathcal{O}.X$ and 
 $D=C + \mathcal{O}$ suit. Then  prove that, modulo $T_1$, every positive primitive formula is equivalent to a universal one (similar to (\ref{elim-as}) above).
 
In the general case, for the task (\ref{task2}),  we cannot always ensure that $B_1(N)\cap Q(N)$  is a ball, or more generally quantifier free definable in $L_\mathcal{O}$: to see this it is enough to consider the image of $x\mapsto x^p$. To handle  this fact, we introduce 
{\bf $p$-th roots of coordinate functions}, called $\lambda$-functions. For instance, for $x\in \mathbb{F}_p((X))$ we write (using the fact that $1,X,\ldots, X^{p-1}$ is basis of 
$\mathbb{F}_p((X))$ over the subfield  consisting of $p$-th powers)
$$ x=x_0^p + x_1^pX+ \dots + x_{p-1}^pX^{p-1}$$
    and define $$\lambda_i(x):=x_i \quad  (i=0,\ldots, p-1).$$

Notice that $\lambda$-functions are both universally and existentially definable in $L$. We prove then $B_1(N)\cap Q(N)$  is positive quantifier-free definable in the language $L$ together with $\lambda$-functions (see Theorem \ref{boule_up}).

This article is organized as follows:
    generalities and definitions about the non-commutative ring $R$, its matrix ring and $R$-modules are given 
     Section 2. In the third section we prove an  equivalent version of Hensel's lemma (see Theorem \ref{hensel-onay}) -we do not know if it is really new-
    and use it to axiomatize what we call {\it henselian filtered $R$-modules}. This part  consists of finding a ball (like $\mathfrak{m}$ above) where the trace of the definable set considered can be defined by a positive quantifier-free formula in the language $L$ together with the $\lambda$-functions. In the fourth section, we introduce the notion of a valued module, and study   pseudo-complements in order to satisfy the task (\ref{task1}) above. In the final section we apply the previous results to obtain  the main theorems. Here we use some general model theory of modules (e.g. Baur-Monk elimination) and some recent facts from \cite{koeningsman-sylvie} and \cite{sylvyarno17}

   In the course of this article we introduce several languages and theories.  For the sections 2-4 (i.e. except the last section), it is worth noting that $K$ and $\mathbb{F}_p((X))$, considered as structures with respect to  these languages are  models of all of these theories. When we say that something is {\it computable} we mean that there is a Turing machine which can compute it (see for instance \cite{koenigsmann2014undecidability} for the related definitions about computability).

\subsection*{Acknowledgements}
This is a long-standing work, which began with my Ph.D thesis in early 2007. I would like to thank  my advisor Fran\c{c}oise Point for suggesting me the topic and  for giving the idea to introduce $\lambda$-functions. I thank  Fran\c{c}oise Delon, also my advisor; without her help, this work could not be finished.
I want to thank  Franz-Viktor Kuhlmann for having invited me, as early as  2008 and for sharing his ideas about the topic. Lou van den dries and Matthias Aschenbrenner have kindly discussed with me  some aspects of the present article in last times. Finally, I want to thank  Sylvy Anscombe and Arno Fehm for their encouragement to write this manuscript.

\section{Preliminaries}

 Let $p$ be a prime number and $k$ a positive integer, we set $d:=p^k$. Let $K$ be the henselization of  the field of rational functions over $\mathbb{F}_d$, with respect to the $X$-adic valuation; denoted as $$K=\mathbb{F}_d(X)^{h}.$$ Let
$$\vfi:x\mapsto x^{d} $$ be the $k$-th power of the Frobenius  endomorphism of $K$. Recall that $K$ is a finite extension of $K^{\vfi}$ of dimension $d$. \footnote{Note that the contents of this section can be generalized to any field  and  to any endomorphism
satisfying  similar hypotheses. In particular one can 
consider $\mathbb{F}_q(X), \mathbb{F}_q((X)), \ldots etc$. As it is our main interest and for readability reasons,  we preferred to stick to one case.} We fix the notation $v_K$ for the $X$-adic valuation on $K$.

  Set $\alpha(0):= \{1\}$ and let $$\alpha:=\alpha(1)=(\alpha_0,\alpha_1,\ldots, \alpha_{d-1})$$ be  a basis of  the $K^{\vfi}$-vector space of $K$. One can think of the basis $(1,X, \ldots, X^{d-1})$.    It is easy to see that for all $n\geq 1$, $\alpha$  induces canonically
the (ordered -) basis $\alpha(n)$ of the $K^{\varphi^{n}}$-vector space $K$. We identify $d^n$ with the set of functions $\{0,\dots,n-1\} \to \{0,\dots,d-1\}$ ordered lexicographically. We denote by $\alpha_k$ the $k$-th element of the basis $\alpha(n)$.


We define the ring \begin{equation}
R:=K\langle t  \mid ta^\vfi=at \rangle .  
\end{equation}
the  (right-) $K$-algebra with the indeterminate $t$, subject to the commutation rule
$ta^\vfi=at$ for all $a \in K$. When 
$k=1$, $R$ is isomorphic to the ring of additive polynomials over $K$, equipped with addition and composition. Every non zero  $q\in R$ can be  written as
\begin{equation}\label{R-polynomial}
q=t^na_n + t^{n-1}a_{n-1} + \dots + a_0 
\end{equation}
where the $a_i$ are from $K$ and $a_n\neq 0$. Since the image of $q$ in $K[T]$ is
$$Q(T)=a_nT^{d^n}+\dots+a_0T$$ under the aforementioned isomorphism, ``the constant term $a_0$ of $q$" gives rise to the linear term $a_0T$ of $Q(T)$.
\begin{defn}
	An element  $q$ as in (\ref{R-polynomial}) is said to be separable if $a_0\neq 0$.
\end{defn}
\begin{remark}
	Note that $q$ is separable if and only if its image $Q(T) \in K[T]$ is separable.
\end{remark}
The integer $n$ in the expression (\ref{R-polynomial}) is called the degree of $q$ and we set the degree of $0$ as $-1$. 
With  the degree function $R$ is right euclidean: for every non zero $r,q \in R$ there exists
$q'$ and a unique $r'$ of degree $<\deg(q)$ such that $r=qq' +r'$. As a consequence, 
least common multiple $\lcm(q,r)$ and greatest common divisor $\gcd(q,r)$ are well defined.
In particular $R$ is right Ore. Note that if $\vfi$ is onto then $R$ is also left euclidean and 
if it is not (which is the case that we consider here), $R$ is not even left Ore. The reader can see \cite{cohn} or \cite{valmod1} Section 2, for more details.

For our quantifier elimination Theorem \ref{boule_up} we will import some definitions and results from \cite{ddp1} on the ring $R$ and 
on $R$-modules. We give the proofs of  some results not to be self-contained but rather, to initiate the reader with the style of the computations that will permit us to refer   easily  to \cite{ddp1}.
\begin{fait}
	The field $K$ is recursively enumerable since $\mathbb{F}_d(X)$ is, and $K$ is formed by adding to $\mathbb{F}_d(X)$ the unique root  of the  each polynomial satifying the hypothesis of the Hensel's Lemma . Consequently the ring $R$ and its matrix rings are recursively enumerable.
\end{fait}

\begin{lemma}\label{lem:decompq}
Let $q\in R$ and $n$ be a positive integer. Let $\alpha$ be a basis of the
$K^\vfi$-vector space $K$. Then 
\begin{enumerate}[1.]
\item  \label{qi} $q$ can be uniquely written as:
  \begin{equation}
\sum_{i \in d^n} q_i\alpha_i \quad (q_i \in R),
\end{equation}

\item for all $n>0$ there is an endomorphism $\nqrt{n}{\cdot}:(R,+) \to (R,+)$ such that
    if  $q=\sum_{i \in d^n} q_i\alpha_i$  then 
        \begin{equation}
t^nq = \sum_{i \in d^n} \nqrt{n}{q_i}t^n\alpha_i.
    \end{equation}
In addition, if $q$ is separable there exists $i \in d^n$ such that
    $\nqrt{n}{q_i}$ is separable.
\end{enumerate}\medskip
\end{lemma} 
\begin{proof}
\item[1.]  Let $0 \leq k\leq s=\deg{q}$. For each monomial $t^ka_k$, by expressing $a_k$ with respect to the basis $\alpha(n)$,  $q$ can be written as
$$t^k(\sum_{i\in d^n} a_{k,i}^{\varphi^n}\alpha_i) .$$ Set 
$q_{k,i}=t^ka_{k,i}^{\varphi^n}$ and  $q_i = \sum_{0\leq k \leq s} q_{k,i}.$ 
So we get 
$$q=\sum_{i \in d^n}q_i\alpha_i.$$

\item[2.] Let $a \in K$. Write 
$$a=\sum_{i \in d^n} a_i^{\vfi^n}\alpha_i .$$ We define
$$\nqrt{n}{a}:=\sum_{i \in d^n} a_i\alpha_i$$ and extend it to $r \in R$ by applying it to its coefficients. By definition $\nqrt{n}{\cdot}$ preserves the addition.

To see that \begin{equation}\label{nqrt} 
t^nq = \sum_{i \in d^n} \nqrt{n}{q_i}t^n\alpha_i
\end{equation} it is enough to show that
$$t^nq_i=\nqrt{n}{q_i}t^n .$$ 
Note that
 $\nqrt{n}{q_{k,i}}=t^ka_{k,i}$ and hence  
 $$\nqrt{n}{q_{k,i}}t^n=t^{n+k}a_{k,i}^{\vfi^{n}}=t^{n}q_{k,i} .$$ 
Using additivity,
we get $\nqrt{n}{q_i}t^n=t^nq_i$.

Suppose now $q$ is separable. If none of the  $\nqrt{n}{q_i}$ is separable, then  for all  $i$, we have  $\nqrt{n}{q_i}=tq_i'$. But then 
$t^nq \in t^{n+1}R$ by the equality (\ref{nqrt}). Hence
$q \in tR$,  which is a contradiction.
\end{proof}


\begin{defn}\label{matriss} An $m\times n$  matrix  $A=(q_{i,j})$ over $R$ is said to be 
\begin{enumerate}[1.]
\item  lower triangular if, 
    $j>i$ implies $q_{i,j}=0$,
\item   lower triangular diagonally  separable 
if it is lower triangular, $n\leq m$ and the  $q_{ii}$ ($i\leq n$) are separable,
\item lower triangular separable  if $A=(A_1,0)$ where $A_1$ is an $m\times k$ lower triangular diagonally separable matrix, 
and $0$ is the $m\times l$ null matrix with $k+l=n$.

\end{enumerate}
\end{defn}

\begin{prop}\label{mattri}
For any matrix $A$ there exists an invertible matrix  $P$  with coefficients 
in  $\{0,1\}$ and an invertible matrix
$Q$ such that $PAQ$ is lower triangular.
\begin{proof}
See the Proposition $6.1$ in  \cite{ddp1}.
\end{proof}
\end{prop}

\subsection{$R$-modules}
For the rest of the article we will always understand the expression $R$-module as right $R$-module. In an $R$-module $M$, scalar multiplication will be denoted as $x.r$, for $x\in M$ and $r\in R$.
\begin{defn}
Let $M$ be an  $R$-module and $\beta$ a basis of the $K^\vfi$-vector-space $K$. 
\begin{enumerate}[1.]
\item $M$ is said to be $t\beta$-decomposable if $x\mapsto x.t$ is injective and we have
\begin{equation}
M=\bigoplus_{i \in d}M.t\beta_i.
\end{equation}
(where $\oplus$ indicates the direct sum as abelian subgroups). 
\item $M$ is said to be $t$-decomposable if it is $t\alpha$-decomposable with
$\alpha=(1,X,\dots, X^{d-1})$.  Furthermore, for $i\in d$ we then define 
$\lambda_i(x)=x_i$ where 
$$x=x_0.t\alpha_0 + \dots + x_i.t\alpha_i + \dots x_{d-1}.t\alpha_{d-1}.$$
\end{enumerate}\medskip
\end{defn}

\begin{remark}\label{lambdascompose}
The direct sum 
$$M=\bigoplus_{i \in d}M.t\beta_i$$ induces the direct sum below for every positive $s$
$$M=\bigoplus_{i \in d^s} M.t^s\beta_i.$$
\end{remark}
\begin{nota}
For the rest of the article we set  
$\alpha:=\{1,X,\dots,X^{d-1}\}$.
\end{nota}

\begin{remark}\label{extunideflambdas}
The functions $\lambda_i$ are both existentially and universally definable in the language of right $R$-modules:\begin{equation}
y=\lambda_i(x) \iffl \forall x_1, \dots,\forall x_i, \dots, \forall x_d \quad x=\sum_{j=1}^d x_i.t.\alpha_j \longrightarrow x_i=y
\end{equation}
and
\begin{multline}
y=\lambda_i(x) \iffl \exists x_1, \dots, \exists x_{i-1}, \exists x_{i+1}, \dots, \exists x_d \quad \\ 
x=\sum_{j=1}^{i-1} x_j.t\alpha_j + y.t\alpha_i + \sum_{j=i+1}^d x_j.t\alpha_j.
\end{multline}

For a positive $s$, by Remark \ref{lambdascompose} above, we get canonically the $\lambda$ functions of level $s$,
defined for all $i\in d^s$, in an obvious way.
\end{remark}

{\bf The language $L(\lambda)$.} For the rest of the article we let $L(\lambda)$ be the language of $R$-modules together with  the functions (-symbols) $\lambda_i$. 	

\begin{defn} We denote by $T_\lambda$ the $L(\lambda)$-theory of $t$-decomposable $R$-modules, that is, the theory of $R$-modules together with the axioms
	
{\bf $\lambda$-decomposition :}\label{ax:lambdas} 
$$\begin{matrix}  \forall x &  x = \sum \lambda_i (x) . t \alpha_i \\
\forall  x \forall (x_i )_{(i\in d)} & (x = \sum_i x_i. t \alpha_i \to  \bigwedge_i x_i = \lambda_i (x)). \end{matrix}$$
\end{defn}

\begin{lemma}\label{Lambdaterms} In any $t$-decomposable $R$-module $M$, any $L(\lambda)$-term can be evaluated on the tuple $(x_i)_i$ from $M$,   as 
$$\sum_i \sum_j \lambda_j(x_i).r_{ij}$$ where 
$r_{ij} \in R$.
\end{lemma}
\begin{proof}
This is  Corollary 3.3 in \cite{ddp1}.
\end{proof}

\begin{lemma}\label{lem:passagelambdasterms}
 Let  $m > 0$, and $q_ j , q_ j' \in R$ such that  $q_j = t^m .q_j'$ . Then the equation
$\sum y_j.q_j = u$ is equivalent to
$$ \bigwedge_{i \in d^{m}}\sum_j y_i.\nqrt{n}{q'_{j_i}}=\lambda_i(u).$$ 
 in any $t$-decomposable $R$-module.
\end{lemma}
\begin{proof}
This is  Lemma 3.4 in \cite{ddp1}
\end{proof}

\begin{obs}
Let $q=(q_0, \dots, q_{n-1})$ be a non zero tuple from $R$. We set
$$e:=\min\{k \tq q_i \in t^kR, \mesp \text{for all}\, i\}.$$
Notice that  $e=0$ means that at least one of the $q_i$ is separable. 

Suppose $e>0$, one can write 
$$q_i=t^eq'_i=\sum_{k \in d^e} \nqrt{e}{q'_i}_{k}t^e\alpha_k .$$
Since at least one of the $q'_i$ is separable, for some $(i,k)$, 
the polynomial  $\nqrt{e}{{q'_{i}}_{k}}$ is separable by Lemma \ref{lem:decompq} (2.).
Set $q_{i,k}= \nqrt{e}{q'_i}_k $  and  let $q^e_\lambda$ be the  $n\times d^e$-matrix $(q_{i,k})$ whose $i$-th line consists of 
the sequence $(q_{i,k})_{k \in d^e}$.
By  this process we have replaced the tuple $q$ by a matrix which  has at least one separable coefficient. Iterating this process and using  Lemma \ref{lem:passagelambdasterms} above we get the following result. 
\end{obs}

\begin{lemma}\label{ddp6.4}
Let  $A$ be a non zero $n\times k$ matrix over $R$. Then, 
the system  $y.A = u$ is equivalent to 
$$y.PQ = w(u)$$ modulo $T_{\lambda}$,
where $P$ is a permutation matrix (i.e. invertible with coefficient in $\{0,1\}$),  $Q$ is lower triangular separable  and $w$ is a tuple 
consisting of $L({\lambda})$-terms. 
\end{lemma}
\begin{proof}
This is a reformulation of  Lemma  6.4.4 in \cite{ddp1}.
\end{proof}

\subsection{Baur-Monk Elimination}\label{sec:baur-monk}

The following is a reminder of Theorem A.1.1, Corollary A.1.2 and the discussion which follows in \cite{hodges}, p. 653-656.

Let $\mathcal{L}$ be any language which contains the language $\{+,-,0\}$ of abelian groups.  A {\it positive primitive formula} ($p.p.$) $\phi$ of 
$\mathcal{L}$, is the one of the form
$$\exists \bar{y} \mesp \left(\bigwedge_i \psi_i(\bar{x},\bar{y})\right)$$
where the $\psi_i$ are atomic.

A {\it group-like $\mathcal{L}$-structure} $A$, is an $\mathcal{L}$-structure whose base set is a group
with respect to $\{+,-,0\}$. A {\it basic formula} for an $\mathcal{L}$-theory $T$, is a $p.p.$ formula which defines  a subgroup of the corresponding cartesian power of any $\mathcal{L}$-structure $A\models T$.

Note that if $S$ is any ring and $\mathcal{L}$ is the language of $S$-modules then in any $S$-module $N$ any $p.p.$ formula defines a subgroup
of the corresponding cartesian power of $N$, hence any $p.p.$ formula is basic for the theory of $S$-modules.

Let $T$ be an $\mathcal{L}$-theory such that every model of $T$ is group-like, 
and every $p.p.$ formula is a basic formula for $T$. An invariant 
$\mathcal{L}$-sentence  of $N\models T$ is  an $\mathcal{L}$-sentence 
$\Theta$ satisfied by $N$, 
such that for some   $p.p.$ formulas $G(x)$ and $H(x)$ of one variable $x$,   for some $m\in \mathbb{N}$,
$$T\models \Theta \Longleftrightarrow  \vert G/G\wedge H \vert =m$$
where the right-hand side of the equivalence is an abbreviation of the formula:
\begin{multline}
\exists  (x_i)_{(i=1..m)}  (\bigwedge_i G(x_i) \wedge \\
\bigwedge_{i\neq j} \neg (G\wedge H)(x_i-x_j)\wedge \forall z \mesp (G(z) \rightarrow \bigvee_i (G\wedge H)(z-x_i))).
\end{multline}
An invariant sentence is an invariant sentence of some $N$.
\begin{theorem}[Baur-Monk]\label{thm:baur-monk}
Every $\mathcal{L}$-formula is equivalent modulo $T$ to a boolean combination of
$p.p.$ formulas and invariant sentences. Hence for all models $N,M \models T$, 
$N \equiv M$ if and only if $N$ and $M$ have same invariant sentences. 
\end{theorem}

\begin{corr}\label{corr:baur-monk}
A completion of $T$ is model-complete if and only if every $p.p.$ formula is equivalent to a universal formula modulo $T$.
\end{corr}

\section{The tropical action of $R$ on $\mathbb{Z}$ and Filtration}

We recall some elementary facts about henselian valued fields. Let $(F,v)$ be a valued field with value group $\mathfrak{G}$ and valuation ring
$\mathcal{O}$. We set $\Gamma=\mathfrak{G}\cup\{\infty\}$,    extend the usual addition of $\mathfrak{G}$ to $\Gamma$ by letting 
$$
\infty+\infty=\infty + a=a+\infty=\infty
$$ for all $a \in \mathfrak{G}$.
 We recall 
the {\it tropicalisation} of a one variable polynomial $Q(T)$ over $F$:
Write $$Q(T)=\sum a_iT^i.$$ Then tropicalisation of $Q$ is the map
$$Q_v: \Gamma \to \Gamma$$
$$\gamma \mapsto  \min_i\{i\gamma + v(a_i)\}.$$

A {\it jump value} (or a tropical zero) of $Q$ is some $\gamma \in \Gamma$ such that
$$ \vert \{i \tq  i\gamma + v(a_i)= Q_v(\gamma)\} \vert \geq 2 .$$
We denote by $\jump{Q}$ the set of jump values of $Q$. Note that this set is finite and has at most $n-1$ element if $Q$ is of degree $n$.

\begin{fait}[Newton's Lemma] Let $(F,v)$ be a valued field and $f$ be a polynomial with coefficients from the valuation ring  
$\mathcal{O}$. Consider the following property $h(f)$ of $f$,
\begin{multline}
 v(f(0))>2v(f'(0)) \Rightarrow (\exists b \mesp 
 f(b)=0 \quad \text{and} \quad v(b)=vf(0)-vf'(0)) \tag{{\it h(f)}}
\end{multline}
Then $(F,v)$ is henselian if and only if $h(f)$ holds for every $f$ over the valuation ring $\mathcal{O}$ of $F$.
\end{fait}
\begin{proof}
See \cite{prestel} Theorem 4.1.3.
\end{proof}

Consider a polynomial $G$,  such that 
$G(0)=0$ and $G'(0)\neq 0$. Then $G$ is of the form 
$$G(T)=aT + \text{sum of monomials of higher degree}.$$
Consider the set $A_1:=\{ \gamma \tq  G_v(\gamma) <   (G-aT)_v(\gamma) \}$. This is a non-empty final segment of $\Gamma$. Also let 
$$A_2:= \{\gamma + v(a) \tq \gamma \in A_1\}.$$ 
We set 
$B_1(G):=v^{-1}(A_1)$ and $B_2(G):=v^{-1}(A_2)$. 

Note that $B_i$  $(i=1,2)$ are convex subsets of $F$, that is, inverse images by $v$ of convex subsets of $vF$.

\begin{theorem}\label{hensel-onay} $(F,v)$ is henselian if and only if  each polynomial $G$ with coefficients in the valuation ring  $\mathcal{O}$, such that
$G(0)=0$ and $G'(0)\neq 0$,  induces a bijection $B_1(G) \to B_2(G)$.
\end{theorem}
\begin{proof}
Suppose $(F,v)$ is henselian. By the definitions of $B_1$ and $B_2$ we have 
$G(B_1)\subseteq B_2$. We will show the converse inclusion. Let $z \in B_2$ and set $f=G-z$. We have 
$vf(0)=v(z)=\gamma + v(a)$ where $\gamma \in A_1$. Hence 
$ G_v(\gamma) =\gamma + v(a) <  (G-aT)_v(\gamma)$. Since a monomial of $G-aT$ is at least of degree $2$ and $G$ is over $\mathcal{O}$, we have 
$$\gamma + v(a)<2\gamma.$$
It follows that $\gamma > v(a)$ and $v(z)>2v(a)$. By the above fact, there is a root $b$ of 
$f$, of valuation $\gamma$. In other words, $G(b)=z$ and hence $G(B_1) \supseteq B_2$. 

Now we show that $G \restriction B_1$ is 1-1. Let $z$ and $f$ be as above and $G(x_1)=G(x_2)=z$ with $x_1,x_2 \in B_1$. Then   $x_1,x_2$ have the same  valuation $\gamma$, but it is well-known that $f$ has at most one root of valuation $\gamma$ (or it can be proven by using Taylor expansions that $v(G(x_1-x_2))>G_v(v(x_1-x_2)$ but this is not possible since $v(x_1-x_2)$ is bigger than then all the jump values of $G$).

The converse is similar and easy (and we won't need).   
\end{proof}

\begin{corr}
	Let $G$ be any polynomial over $(F,v)$ such that $G(0)=0$ and $G'(0)\neq 0$ then there exist convex subsets $B_1$, $B_2$ such that
$G\restriction B_1:B_1 \to B_2$ is a bijection.
\end{corr}
\begin{proof}
Divide $G$ by the  coefficient which has the minimal valuation among the coefficients of $G$. Let $H$ be the obtained polynomial (which has coefficients over $\mathcal{O}$). Then $B_1:=B_1(H)$ and $B_2:=aB_2(H)$ suits. 
\end{proof}

\begin{nota}
	Let $G$ be as above. If the value group of $F$ is discrete of rank 1, then $A_1,A_2$ are closed intervals; in this case, we set $$\h(G):=\min A_1 \quad \text{and} \quad \hens(G):=\min{A_2}.$$
\end{nota}
{\bf The tropical action of $R$.} From now on we set 
$$\Gamma:=\mathbb{Z}\cup \{\infty\}$$ and we equip $\Gamma$ with a right action of $R$ using tropicalisations:

Let 
$q=\sum_i t^ia_i \in R$.  We define the tropicalisation 
of $q$, as the tropicalisation of the $Q(T)=\sum_i a_iT^{d^{i}}$, that  is,  as  the map 
\begin{equation}
\cdot q:\Gamma \to \Gamma; \quad \gamma \mapsto \gamma\cdot q= \min_i\{d^i\gamma + v(a_i)\}.
\end{equation}
In particular 
$\cdot q$ is strictly increasing  if  $q\neq 0$.
Note that $\gamma \cdot a =  (aT)_v(\gamma) =\gamma + v(a)$ for $a \in K$.

We set  $\jump{q}:=\jump{Q}$, $\h(q):=\h(Q)$ and $\hens(q):=\hens(Q)$. For instance,
if $q=t-1$ then $A_1=A_2:=\{\gamma \in \Gamma \tq \gamma>0\}$. Hence
$\h(q)=\hens(q)=1$.

For our interests, we also introduce the tropicalisations of the $\lambda$-functions. Remark that we have for all $x\in K$,
$x=\sum_{i \in d} \lambda_i(x)^dX^i$, 
       $$v_K(x)=\min_i\{v_K(\lambda_i(x)^d)+i\},$$
and the minimum is attained for a unique $i\in d$.
We define 
$$ \lambda_i(\gamma):=\begin{cases} \frac{\gamma - i}{d}\quad \text{if} \mesp \gamma \in d\mathbb{Z} + i \\
        0 \quad \text{else}
        \end{cases}$$
and we set 
$$\lambda(\gamma):=\sum_i \lambda_i(\gamma).$$

\begin{remark}\label{rem:lambdaregular} For $x\in K$ 
$$v_K(x)\geq \gamma \Leftrightarrow \quad \text{for all}\mesp i\in d \quad v_K(\lambda_i(x)) \geq \lambda(\gamma).$$
\end{remark}
\begin{proof} Note that $\Leftarrow$ is clear and if $\gamma=v_K(x)$ then the assertion is trivial. Suppose $v_K(x)>\gamma$. Let $j$ be such that 
$v_K(\lambda_j(x).t^{d}\alpha_j)=dv_K(\lambda_j(x))+j=v_K(x)$. Let $i$ be such that 
$\gamma=d\lambda(\gamma)+i$.  It follows that
$$d(v_K(\lambda_j(x))-\lambda(\gamma))>j-i.$$
Hence $v_K(\lambda_j(x))-\lambda(\gamma)>-1$ and
$v_K(\lambda_j(x))-\lambda(\gamma)\geq 0$.
\end{proof}
\begin{corr}\label{corr:lambas(s)}For all non zero $x\in K$ and $\gamma \in \Gamma$, 
$$v_K(x)\geq \gamma \Leftrightarrow \quad \text{for all}\mesp i\in d^s \quad v_K(\lambda_i(x)) \geq \lambda^s(\gamma)$$
where $\lambda^s=\underbrace{\lambda\circ\ldots\circ\lambda}_{s-times}.$
\end{corr}

\begin{remark}Let $\gamma\in \Gamma$ and $r,q \in R$. Then
\begin{enumerate}[1.]
\item $\gamma\cdot rq = (\gamma \cdot r)\cdot q$,
\item $\gamma\cdot (r + q)\geq\min\{\gamma\cdot r, \gamma \cdot q\}$,
\item $\cdot r$ is strictly increasing for all non zero $r$, 
\item $\infty \cdot r =  \gamma \cdot 0 = \infty$ for all $r$ and $\gamma$.
\end{enumerate}
\end{remark}
\begin{proof}
This follows by direct computations from the definition. Note that (1) follows more generally from  \cite{gonenc}, Corollary 4.1.14.
\end{proof}

\begin{fait}
	The theory of $\Gamma$ together with the tropical action of $R$ and
	$\lambda$, is decidable since this structure  is definable in the ordered abelian group structure of $\mathbb{Z}$, together with $\infty$ and constants for the elements of $\mathbb{Z}$.
\end{fait}

\begin{defn}	
	We call this structure the tropical structure of $\Gamma$.	
\end{defn}
\subsection{Quantifier Elimination Near 0}

\subsection*{\bf The languages $L_{\mathcal{O}}$ and $L_{\mathcal{O}}(\lambda)$} The language  
$L_{\mathcal{O}}$ is obtained by adding to $L$ a unary predicate $\mathcal{O}$, and the language $L_\mathcal{O}(\lambda)$ is the language $L_\mathcal{O}$ together with the functions (-symbols) $\lambda_i$ added to $L_\mathcal{O}$. We want to study the divisibility conditions for $t$-decomposable $R$-modules which can be seen
as  properties reflecting a kind of henselianity, analog to Theorem \ref{hensel-onay} just above.

\begin{nota}
For the rest of this article, for $\gamma \in \Gamma$,
we write $P_\gamma$ 
 for the predicate  $\mathcal{O}.X^{\gamma}$, i.e. in any $L_\mathcal{O}$-structure
 $M$, $x\in P_\gamma$ if and only if $x.X^{-\gamma} \in \mathcal{O}$.
\end{nota}
\begin{defn}
An henselian filtered module, is a $t$-decomposable $R$-module which is an $L_\mathcal{O}(\lambda)$-structure satisfying the following axioms:

\item[0. {\bf Balls:}] $P_{\infty}=\{0\}$  and the $P_{\gamma}$ form a chain of subgroups decreasing with $ \gamma$  such that  the  inclusions are proper.
\item[1. {\bf Ultrametric:}]  $\forall x  \forall y  \mesp x \in P_{\gamma}\wedge y \in P_{\delta}  \rightarrow (x.r + y) \in P_{\min\{\gamma\cdot r,\delta\}}$   for all $\gamma,\delta \in \Gamma$ and $r\in R$.
\item[2. {\bf Regularity:}]   $\forall x \mesp  x \in P_{\gamma} \leftrightarrow x.r \in P_{\gamma \cdot r}$, for all $\gamma$ and $r\neq 0$ such that
$\gamma \notin \jump{r}$.
\item [3. $\lambda$-{\bf regularity:}]  $\forall x \, x\in P_{\gamma} \leftrightarrow \bigwedge_i \lambda_i(x) \in P_{\lambda(\gamma)}$, for all $\gamma$. 
\item [4. {\bf Henselianity:}] $\forall x  \in P_{\hens(s)}\setminus \{0\} \exists ! \, y \in P_{\h(s)} \mesp y.s=x$ for all separable $s$.

\end{defn}
Note that axiom 3 implies that
$$\forall x \, x\in P_{\gamma} \leftrightarrow \bigwedge_{i\in d^s} \lambda_i(x) \in P_{\lambda^s(\gamma)},$$ a consequence analog to the one expressed in Corollary
\ref{corr:lambas(s)}.   
\begin{nota}
	We denote by $T_{hens}$ the theory of henselian filtered modules.

We also isolate some theories of the $R$-modules that are already considered in \cite{ddp1}:
\begin{itemize}
	\item we denote by $T_{free}$ the theory $T_\lambda$ together with the following (scheme of-) axioms:
	\label{ax:separablyclosed} $$\forall x \exists y \mesp x=y.s,$$
	 for all separable $s$,
	\item  we denote by $T^{0}_{free}$, the $L_{\lambda}$-theory of torsion-free non zero models of $T_{free}$.
\end{itemize}
\end{nota}
\begin{remark}\label{P_M} In a henselian filtered module $M$, we denote by 
$P_M$ the subgroup defined by  the intersection 
$$P_M:=\bigcap_{\gamma\neq \infty} P_\gamma(M).$$
Note that by  ultrametric and regularity axioms, $P_M$ is a 
$L(\lambda)$-substructure of $M$ and it is torsion-free as an $R$-module. Moreover, it is straightforward to check that by Hensel's axioms 
$$P_M \models T^0_{free}.$$
\end{remark}
\begin{proposition}\label{prop:t0freedecidable}  Given $L$-$p.p.$ formulas	$\mathfrak{a}(x)$ and $\mathfrak{b}(x)$  of one variable $x$, 
the quotient $\mathfrak{a}/(\mathfrak{a}\wedge \mathfrak{b})$ is either trivial or infinite in every model of $T^0_{free}$. In particular $T^0_{free}$ is complete and decidable. 
	\begin{proof}
This follows from Lemma 6.8 in \cite{ddp1}.
	\end{proof}
	
\end{proposition}

\begin{defn}
	A ball of $L_\mathcal{O}$, is an atomic formula $W(x_1,\ldots,x_k)$ of the form
	$$W(x_1,\ldots, x_k): \bigwedge_{i=1}^k x_i \in P_{\gamma_i}.$$ We will write it rather as a product 
	$$W=P_{\gamma_1}\times \dots \times P_{\gamma_k}$$ of predicates. 
	$W$ is said to be proper if none of the ${\gamma_i}$ is equal to $\infty$. 
\end{defn}

\begin{remark}
	A  positive primitive  formula $\phi(\bar{x})$ of $L_{\mathcal{O}}$ is equivalent to one in the form
	$$\exists \bar{y} \mesp\bar{x}.B-\bar{y}.A\in \mathcal{O}^k\times\{0\}^n$$ which it self is equivalent to a formula 
	$$\exists \bar{y} \mesp\bar{x}.B'-\bar{y}.A'\in W$$ 
	where $A,B,A',B'$ are matrices over $R$ and $W$ is a ball. Hence any $p.p.$ formula is a basic formula. For our purposes will rather use the latter
	equivalence.
\end{remark}

Note that modulo $T_{hens}$, the set of  definable sets  by an $L_{\mathcal{O}}$-$p.p.$ formula contains the set
of definable sets by an $L$-$p.p.$ formula 
since $W$ can be
chosen equal to $P_{\infty}^k$ for some $k$.  

In this section, we will prove the following:
\begin{theorem}\label{boule_up}
	Let $\phi(x_1,\dots,x_m)$ be a $p.p.$ formula of
	$L_\mathcal{O}$. Then there is some computable $\delta\neq \infty$ and some positive quantifier free
	$L(\lambda)$-formula $\psi$ which depends only on the theory $T_\lambda$, such that, with
	$V=P_\delta^m$,  we have
	$$T_{hens} \models \phi \wedge V \leftrightarrow \psi \wedge V .$$ 
	Moreover,
	If $\phi_1$ is another $p.p.$ formula with the same arity, and
	if $\phi$ and $\phi_1$ are equivalent modulo $T^0_{free}$, then for some computable $\gamma$ and with $W:=P_\gamma^m$ 
	$$ T_{hens} \models \phi\wedge W \leftrightarrow \phi_1 \wedge W .$$
\end{theorem}

\begin{lemma}\label{croissant} Given $\delta\neq \infty$ and an $L(\lambda)$-term $u(\bar{x})$, there exists $\gamma\neq \infty$ such that 
\begin{equation}
T_{hens} \models \forall \bar{x} \mesp \bar{x} \in P_\gamma^{\vert \bar{x} \vert} \rightarrow u(\bar{x}) \in P_{\delta}.
\end{equation}
\end{lemma}
\begin{proof} Put $u(\bar{x})$ in the form $\sum_{i,j} \lambda_j(x_i).r_{ij}$ using Lemma \ref{Lambdaterms}. 
    By regularity, $\lambda$-regularity and ultrametric axioms, for any $\rho$,
$$\bar{x} \in P_\rho^{|\bar{x}|} \rightarrow u(\bar{x}) \in P_{\min\{\lambda(\rho)\cdot r_{ij}\}}$$ holds. Since all the tropical functions $\lambda$ and $\cdot r_{ij}$ are unbounded and increasing,
    one can choose $\gamma$ such that $\gamma \cdot r_{ij} \geq \delta$ for all $i$.
\end{proof}
\begin{remark}\label{rem:gamma_computable}
	The value $\gamma$  is computable from $\delta$ and $u(\bar{x})$ in the tropical structure of $\Gamma$ .
\end{remark}

As a consequence of the henselianity axioms, we  observe the following.
\begin{lemma}
Given  $\gamma\neq \infty$ and a separable $s$ there exists  $\delta\neq \infty$ such that

$$ T_{hens} \models \forall y \mesp  [(y \in P_\delta) \rightarrow (\exists x  \in P_\gamma \wedge x.s=y)].$$

\end{lemma}
\begin{proof}
If $ \gamma \leq \h(s)$ then set $\delta=\hens(s)$. Otherwise $\gamma \notin \jump{s}$ and by the 
regularity axioms $x \in P_{\gamma} $ if and only if $x.s \in P_{\gamma \cdot s}$. 
Since $P_\gamma \subseteq P_{\h(s)}$, $P_{\gamma \cdot s} \subseteq  P_{\hens(s)}$. 
Hence if $y \in P_{\gamma \cdot s}$ then the unique solution $x$ such that $x.s=y$ lies in 
$P_\gamma$.  So $\delta=\gamma \cdot s$ fits for our requirements.
\end{proof}

\begin{remark}\label{rem:delta_computable}
	As above, $\delta$ is computable in the tropical structure $\Gamma$. 
\end{remark}

\begin{corr}\label{bouledivisible}
Let $A=(a_{ij})$  be an $m\times k$ lower triangular diagonally separable matrix (in particular  $k \leq m$). Then,  
for every proper ball $W$, there exists a proper computable ball $W_1 $ such that  
$$T_{hens} \models \forall \bar{x} \mesp (\bar{x}\in W_1 \rightarrow \bar{x}\in  W.A) .$$ 
\end{corr} 
\begin{proof}
Write $W = \prod_{i=1}^{m} P_{\delta_i}$. By Remark \ref{croissant},  for fixed $i$, chose  $\gamma_i  \in \Gamma\setminus\{\infty\}$ such that,
$$\sum_{j \neq i} P_{\gamma_i}.a_{ji} \subset P_{\delta_i}.$$
 Since the $a_{ii}$ are separable, by the above lemma there exist proper balls $U_i$ such that whenever   $z_i \in U_i$,  
 there exists $y_i \in P_{\gamma_i}$, such that $z_i = y_i.a_{ii}$; hence 
$$z_i - \sum_{j=1}^{n} y_i.a_{ji}= \sum_{j\neq i} y_i.a_{ji} \in P_{\delta_i}$$ for all $1\leq j \leq k$. Take $W_1=
\prod_{i} U_i.$ 
\end{proof}

{\bf Proof Theorem \ref{boule_up}.}
Let $\bar{y}:=(y_1,\dots,y_k)$, $\bar{x}:=(x_1,\dots,x_m)$,  $\phi(\bar{x})$,
    $W:=\prod P_{\gamma_i}$ and
$$\phi(\bar{x}): \exists \bar{y} \mesp \bar{x}.B -\bar{y}.A \in W.$$

Let $I:=\{i \tq \gamma_i=\infty \}$ and $J:=\{j \tq \gamma_j\neq \infty \}$, 
 and $A_I$  be the matrix formed by the columns $C_{i \in I}$ (resp. 
     $A_J$ be the matrix formed by the columns $C_{j\in J}$ of $A$). We set
     $u(\bar{x}):=\bar{x}.B$ and  by $u_I(\bar{x})$ (resp.  $u_J(\bar{x})$) we denote the
     tuple formed by $I$-coordinates (resp. $J$-coordinates) of $u(\bar{x})$. We may assume that $\phi$ is of the form
     \begin{equation}\label{firstred}
      \exists \bar{y} \mesp \left( {u_I(\bar{x})} = \bar{y}.A_I \wedge u_J(x) -
\bar{y}.A_J \in W_J \right)
     \end{equation}
where $W_J$ is the obvious projection of $W$ to its non zero coordinates.
By Lemma \ref{ddp6.4}, there exists a lower triangular separable matrix  $\widetilde{A}_I=(S,0)$ 
 such that the formula 
 $$u(\bar{x})=\bar{y}.A_I$$
is equivalent, modulo $T_\lambda$, to 
$$
(t_1(u_I(\bar{x})),\dots,t_{n-l}(u_I(\bar{x})))=\bar{y}.PS \wedge (t_{n-l+1}(u_I(\bar{x})), \dots ,t_{n}(u_I(\bar{x})))=0,
$$ 
 where $P$ is a permutation matrix and the $t_1(u_I(x)),\ldots t_n(u_I(\bar{x}))$ are some $L(\lambda)$-terms.
	   
     By remark \ref{croissant}, chose a proper ball $U'$ such that $U'.A_J \subseteq W_J$
     and by Corollary \ref{bouledivisible} chose a proper ball $V'$ such that $V' \subseteq U'.PS$.
     
     By remark \ref{croissant} again, we choose $V$ such that for all 
     $\bar{a} \in V$,
     $$(t_1(u_I(\bar{a})),\dots,t_{n-l}(u_I(\bar{a}))) \in V'$$ and  $u_J(\bar{a}) \in W_J$.

     Hence:
     \begin{equation} T_{hens} \models \forall \bar{x} \mesp  \left(\phi(\bar{x})\wedge \bar{x}\in V\right) \longleftrightarrow \left((t_{n-l+1}(u_I(\bar{x})), \dots ,t_{n}(u_I(\bar{x})))=0\wedge \bar{x}\in V\right).
     \end{equation}
     
 Denote by $\psi$ the formula 
 $$(t_{n-l+1}(u_I(\bar{x})), \dots ,t_{n}(u_I(\bar{x})))=0.$$
 
 The first statement is now proved.

 Now if $\phi_1$ is another $p.p.$ formula, modulo $T^0_{free}$, $\phi_1$ is equivalent  to some positive quantifier free formula $\psi_1$. Since $T^{0}_{free}$ is decidable, we have an algorithm which checks if $T^0_{free} \models \psi \leftrightarrow \psi_1$. 
 Modifying this algorithm we can remember the finitely many
 non zero $r \in R$, such that algorithm uses the axiom
 $$\forall x \mesp  x\neq 0 \to  x.r\neq 0,$$ and the finitely many separable $s$, such that, the algorithm uses the axiom 
 $$\forall x \exists y \mesp y.s=x.$$ Choose $\gamma$ bigger than
 all the $\max \jump{r}$ and $\max\{\h(s),\hens(s)\}$ for all the $r$ and $s$ as above. Set $W:=P_\gamma^m$. Hence for any non zero $x \in P_\gamma$,
 $x.r\neq 0$ and there exists   $y \in P_\gamma$, with $y.s=x$. It follows that the same algorithm computes a proof of
$ W\wedge \psi(\bar{x}) \leftrightarrow W\wedge \psi_1(\bar{x})$ from $T_{hens}$.  In particular,
 we have $\phi_1 \wedge W \leftrightarrow W\wedge \phi(W)$ modulo $T_{hens}$.
 \begin{remark}
 	By passing to an $\omega_1$-saturated model
 	$M$, since $P_M\models T^0_{free} $, it is trivial that if $T^0_{free} \models \phi \leftrightarrow \phi_1$, then for some proper $W$, $$\phi(W) = \phi_1(M).$$ What we show above is that the decidability of $T^0_{free}$ yields the computability of $W$.
 \end{remark}
 
 \begin{corr}\label{corr:boule_up}
 	Let $\phi$ and $\phi_1$ be $p.p.$ formulas of $L_{\mathcal{O}}$ and $\psi$ and $\psi_1$ are given as in the proof of the above theorem. Then, there is a computable ball $W$ such that 
 	$$\vert (\phi\wedge W)/(\phi_1 \wedge W)\vert =1,$$ or 
 	for all proper ball $V \subseteq W$, 
 	 	$$\vert (\phi\wedge V)/(\phi_1 \wedge V)\vert = \infty.$$
 \end{corr}
 \begin{proof}
 Follows by Proposition \ref{prop:t0freedecidable} and by the above theorem.
 \end{proof}
 

       \section{pseudo-complements}
       We will introduce the notion of a {\it valued module}  to study henselian filtered modules. These are  $R$-modules  $M$, equipped with a function $v:M\to \Gamma$ inducing the ultrametric topology.  After investigating the elementary properties of valued modules, we will get the consequences that can be expressible in the language $L_{\mathcal{O}}$.
 
\begin{defn}\label{def:valmod} A valued module  is a $t$-decomposable $R$-module together with a surjective map $v:M \to \Gamma$ such that for all $x,y\in M$, 

\begin{enumerate}[1.]
\item  $v(x\pm y)\geq \min\{v(x),v(y)\}$
\item  $v(x)=\infty \leftrightarrow x=0$
\item  $v(x)\notin \jump{r} \to v(x.r)=v(x)\cdot r$,  for all $r\in R$.
\end{enumerate}\medskip
\end{defn}

\begin{remark}\label{rvR} Let $r=t^na + \dots + t^ka_k \in R$ where monomials are written following decreasing degrees,  then
\begin{enumerate}[1.]
\item $v(x.t^ia_i)=v(x)\cdot t^ia_i$, for all $x \in M$,
\item  $v(x.r)=v(x)\cdot t^ka_k <v(x.(r-t^ka_k))$ whenever $v(x)>\max \jump{r}$,
\item  $v(x.r)=v(x)\cdot t^na_n<v(x.(r-t^na_k))$ whenever $v(x)<\min \jump{r}$.
\end{enumerate}
\end{remark}
\begin{proof}
1. Follows from  Definition \ref{def:valmod} (3.) since a monomial has no jump value.

2. By Definition \ref{def:valmod} (3.), if $v(x)>\max \jump{r}$ then $v(x.r)=v(x)\cdot r$. 
Let $\gamma > \max \jump{r}$. Then for some $i$, $\gamma\cdot r=\gamma \cdot t^ia_i<\gamma \cdot t^ja_j$ for all $j\neq i$.  In other words the line  $\{(\delta,d^i\delta +v(a_i)\}_\delta$ does not intersect any other line $\{(\delta,d^j\delta +v(a_j)\}_\delta$ in the area $(\max \jump{r}, \infty]\times \Gamma$. This can happen only if $i<j$ for all $j\neq i$. Hence $i=k$.

3. The proof is very similar to (2.)
\end{proof}

 Let $(M,v)$ be a  valued module.
We define the equivalence relation
$\mathbf{RV}$ on $M$ by 
$$x \mathbf{RV} y \iffl v(x)=v(y)<v(x-y).$$ We denote the $\mathbf{RV}$-class of an element $x$ by $\rv{x}$ whereas $\rv{A}$ stands for the set $\{\rv{x} \tq x\in A\}$ for $A\subseteq M$.
We also set the notation $P_\gamma$  for the closed ball of radius $\gamma$ centered at $0$.
   
\begin{defn} For subgroups $A$ and $B$ of $M$, we say that $A$ and $B$
\begin{itemize}
\item   are $m$-immediate (m stands for \emph{mutually}) if $\rv{A\setminus P_\gamma}=\rv{B\setminus P_{\gamma}}$ for some $\gamma \in \Gamma$, and we write
$A \approx B$,
\item   are pseudo-orthogonal if $vA\cap vB \subseteq [\gamma, \infty]$ for some $\gamma \in \Gamma$, and we write $A | | B$.
\end{itemize}\medskip
\end{defn}
\begin{remark}
   $A\cap C \subseteq P_{\gamma}$ for some $\gamma$, whenever $A || C$.
\end{remark}
\begin{remark}\label{rem:modulogamma} If $A\approx B$ and $B\subseteq A$ then for some $\gamma$,
$A+P_\gamma=B+P_\gamma$.
\end{remark}
\begin{proof}
Let $\gamma$ be such that  $\rv{A\setminus P_\gamma}=\rv{B\setminus P_{\gamma}}$.
It is enough to show $\{x\in A \tq v(x)\leq \gamma\} \subseteq B+P_\gamma$.
We proceed by induction w.r.t. the dual order on the initial segment 
$(-\infty ,\gamma]$:
Let $a\in A$. If $v(a)=\gamma$ there is nothing to do. 
Suppose $v(a)<\gamma$ and for all $a' \in A$ of valuation $>v(a')$ there is some
$b' \in B$ such that $a-b'\in P_\gamma$. Since $A\approx B$, there is $b\in B$ such that
$v(a-b)>v(a)$. Since $b\in B$ and $B\subseteq A$, $a-b \in A$. Now by applying 
the induction hypothesis to $a-b$ we have $a-b-b'\in P_\gamma$ for some $b'\in B$. That is, $a \in B+P_\gamma$.
\end{proof}
\begin{defn} A pseudo-complement of a subgroup $A$ is a subgroup $C$ such that,
    for some $\gamma$ 
    \begin{equation}
    M=A + C + P_\gamma \quad \text{and} \quad A|| C.
    \end{equation}
\end{defn}

\begin{remark}\label{rem:invimage}
	If $A,C$ and $P_\gamma$ are as above and $f:M\to M$ is an additive map such that $f^{-1}(P_\gamma) \subseteq P_\delta$ for some $\delta$, then $f^{-1}(C)$ is a pseudo-complement to $f^{-1}(A)$. In particular, this is the case when 
	$f$ is given by a scalar multiplication.
\end{remark}

\begin{remark}\label{rem:+ballpseudocomp}
If $C$ is a pseudo-complement of $A$ then $C$ is a pseudo-complement
	of $A+P_{\gamma}$ for any $\gamma$.
\end{remark}

\begin{remark}\label{rem:samecomplement} It is straightforward to see that if $A$ and $B$ have the same pseudo-complement then $A \approx B$. The following lemma establishes the converse using that the value set is $\mathbb{Z}\cup \{\infty\}$.
\end{remark}
\begin{lemma}\label{changecomp} If $A\approx B$ then $C$ is a pseudo-complement of $A$ if and only if it is a pseudo-complement of 
$B$.
\end{lemma}
\begin{proof}
Suppose $C$ is a pseudo-complement of $A$  satisfying
$$   M=A + C + P_\gamma$$
and $A\approx B$. Let $\delta$ be such that  $\rv{A\setminus P_\delta}=\rv{B\setminus P_\delta}$. We may assume that  $\delta \leq \gamma$. We claim that 
$$M=B+C+P_\delta .$$

	We proceed by induction on $\{\beta \in \Gamma \tq \beta\leq \delta\}$ as in the proof of Remark \ref{rem:modulogamma}. If $v(z)=\delta$ then trivially $z \in B + C + P_\delta$. Let $z \in M$ and suppose that for all $x$ with $v(x)>v(z)$ there exist 
$(b,c,x_\delta) \in B\times C\times P_{\delta} $ such that 
$$x=b+c+x_{\delta}.$$
Write $z=a+c+z_\delta$ for some $(a,c,z_\delta) \in A\times C \times P_\delta$ (in fact $z_\delta$ can be chosen in $P_\gamma \subseteq P_\delta$). If $v(a)\geq \delta$ there is nothing to do. If $\delta>v(a)>v(z)$ then by induction hypothesis
$a=b+c'+z'_\delta$ with $(b,c',z_\delta) \in B\times C \times P_\delta$, and hence 
$$z=b+c+c'+z_\delta + z'_\delta. $$ The only possibility which remains to be considered is $v(z)=v(a)<\delta$ since $A || C$ implies that
$v(a+c+z_\delta)=\min\{v(a),v(c)\}$ whenever $v(a)<\delta$. Pick  $b\in B$ such that
 $\rv{b}=\rv{a}$. Then, since $v(a-b)>v(a)$, by induction hypothesis the equality
$$a-b=b'+c'+z'_\delta$$ holds for some $ (b,c',z_\delta') \in B\times C \times P_\delta$. Hence $z \in B+C+P_\delta$. 

Now to see that $B || C $ we claim that $vB \cap vC \in [\delta, \infty]$. Suppose for a contradiction that
for some $b \in B $ and $c\in C$, $v(b)=v(c)<\delta$. Then we can choose $a \in A$ such that
$v(a)=v(b)$ hence $v(a)=v(c) \in vA \cap vC$. But then $v(a) \in [\gamma, \infty]$. 
This is a contradiction since $\gamma \geq \delta$.  
\end{proof}
\begin{lemma} Let $A,A',B,B'$ be  such that 
$A \approx A'$,  $B\approx B'$ and $A || B$. Then the following holds: 
\begin{enumerate}[1.]
\item  $A' || B'$, 

\item  every pseudo-complement of $A+B$ is a  pseudo-complement of $A'+B'$.
\end{enumerate}\medskip
\end{lemma}
\begin{proof}
1. Let $\gamma=v(a')=v(b')$ with $a'\in A$ and $b'\in B$. Let $\delta$ be such that $vA\cap vB \subseteq [\delta,\infty]$ and $\gamma_1$, $\gamma_2$ are respectively the values yielding $A\approx A'$ and $B\approx B'$. We claim that  
$\gamma\geq \min\{\delta,\gamma_1,\gamma_2\}$. Let $a \in A$ and $b\in B$ be such that $\rv{a}=\rv{a'}$ and
$\rv{b}=\rv{b'}$. Then $\gamma=v(a)=v(b)$. Hence $\gamma$ must be $\geq \delta$.

2. Let $C$ be a pseudo-complement of $A+B$. Let $\delta,\gamma_1,\gamma_2$ be as above. Let
	$z=a+b$ with $a\in A$ and $b \in B$ such that $v(z)<\min\{\delta,\gamma_1,\gamma_2\}$. Then $v(a)\neq v(b)$ and $v(z)=\min\{v(a),v(b)\}$ necessarily. Suppose $v(a)=v(z)$.
Then $\rv{a}=\rv{z}$ hence for some $a'\in A'$, $\rv{z}=\rv{a'}$.  If $v(z)=v(b)$ we can choose  in the same way $b'$ such that $\rv{b'}=\rv{z}$. Hence $A+B \approx A' + B'$.  We apply now Lemma \ref{changecomp}.
\end{proof}

\begin{defn}
	A valuation independent basis $\beta$ of the $K^{\vfi^{n}}$-vector space $K$,  is a basis such that the members of $\beta$ have all different valuations in the finite set $\{0,\ldots,d^{n}-1\}$. 
\end{defn}
For example, the basis  $\alpha$ of the $K^{\vfi}$-vector space $K$ is valuation independent.
\begin{remark}
If $\beta$ is a valuation independent basis of the $K^\vfi$-vector space $K$ then for all $s>0$, $\beta(n)$ is a valuation independent basis of $K^{\vfi^n}$-vector space $K$. Moreover,  any valued module is 
$t\beta$-decomposable for all valuation independent $\beta$.
\end{remark}


\begin{lemma}
Let $g_1,\dots,g_m$ be all of the same degree $s$  such that the leading coefficients $b_1, \dots, b_m$ have distinct valuations  in
	$\{0,\dots,d^s-1\}$. Then $M.g_i || M.g_j$ whenever $i\neq j$. 
\end{lemma}
\begin{proof}
	By \ref{rvR}, we have $v(x.g_i)=v(x)\cdot g_i=v(x)\cdot t^sb_i=p^{d^s}v(x)+v_K(b_i)$ whenever $v(x)<\jump{g_i}$ for $i=1\dots m$. In particular, 
if $v(x)<\min\{\jump{g_i},\jump{g_j}\}$ for a $j\neq i$ then $v(x.g_i)\neq v(x.g_j)$.
\end{proof}
\begin{lemma}\label{lem:addtoMta} Let $\beta$ be a valuation independent basis of the $K^{\vfi}$-vector space $K$. Let $q=t^sa+\dots \in R$ be of degree $s$. 
Then there is a unique $j=j(q)$ such that $\beta_j \in \beta(s)$ and $M.q \approx M.t^s\beta_j$. As a consequence $C:=\sum_{j' \neq j(q)} M.t^s\beta_{j'}$ is a pseudo-complement for $M.q$.
\end{lemma}
\begin{proof}
Write
$$a=\sum_{j \in d^s} a_j^{\varphi^s}\beta_j. $$

Let $j$ be such that $v_K(a)=v_K(a_j^{\vfi^s}\beta_j)<v_K(a_{j'}^{\vfi^s}\beta_{j'})$ for all $j' \neq j$ (such a $j$ exists since $\beta$ is valuation-independent). 
Then $$v(x.t^sa)=v(x)\cdot t^s(a_j^{\vfi^s}\beta_j)<v(x)\cdot t^s(a_{j'}^{\vfi^s}\beta_{j'})$$ for all $x\in M$ and $j'\neq j$. 

 Let $\gamma < \min
\jump{q}$. Then by Remark \ref{rvR}, for all $x$ of such that $v(x)\leq \gamma$, we have $\rv{x.q}=\rv{x.t^sa}$, hence 
$$\rv{x.q}=\rv{x.t^sa_j^{\vfi^s}\beta_j}=\rv{(x.a_j).t^s\beta_j}.$$ Since $x\mapsto x.a_j$ is a bijection it follows that $M.q \approx M.t^s\beta_j$. The consequence follows by Lemma \ref{changecomp}.
\end{proof}

\begin{lemma}\label{vddku}
Given a non zero matrix $Q$ with coefficient over $R$, say with $k$ rows, there exists a matrix $Q'$ such that
the first column of $Q'$ consists of polynomials which have all the same degree $s$,
	such that the leading coefficients in this column have distinct valuations  in $\{0,\ldots, d^s-1\}$  and
$$M^k.Q =  M^k.Q' \mesp .$$
\end{lemma}
\begin{proof} This proof is essentially a slight generalization  of the proofs of 
 Lemma 3 and Lemma 4 in \cite{vdd-kuhlmann}.
 
 Let $Q=(q_{ij})$ be a matrix over $R$ with $k$-many rows.
	We will proceed by induction on $f=\sum_{\{(i,j)\tq q_{ij}\neq 0\}}(\deg(q_{ij})+1)$.

	Since $Q$ is non zero, $f>0$. Suppose $f=1$. We may  assume 
	$q_{11}=c \in K^{\times}$ and all other entries of $Q$ are zero. Then 
	$\bar{x}=\bar{y}.Q$ for some $\bar{y}$ if and only if, all the coordinates of $\bar{x}$ except possibly the first one, are $0$. Hence we can take $Q'$ the matrix which has $1$ at the position
	$(1,1)$ and has all other entries equal to zero.
 
	Now we suppose $f>1$. Let $e:=\max\{\deg{q_{ij}}\}\geq 0$. We may suppose that $q_{11}$ has degree $e$.
Set $$e_{ij}:=\deg{q_{ij}},  \mesp e_i:=e_{i1}$$ and 
$$c_{ij}:=\text{the leading coefficient of}\mesp q_{ij},  \mesp c_i=c_{i1}.$$
 
 \noindent
 {\bf Claim 1}: We may assume that for all $a_1,\ldots a_k \in K$, not all are $0$,
 $\sum_i a_i^{\vfi^{e_i}}c_i \neq 0$.
 
 Suppose $\sum_i a^{\vfi^{e_i}}c_i = 0$. We may also suppose that $a_1=1$.
 We define for all $j$,  
 $$\widetilde{q_{1j}}=\sum_{i=1}^{k} t^{e-e_{i}}a_iq_{ij}.$$

Since $a_1=1$ and $e=e_1$ 
\begin{equation}\label{eq:matricepassage}
	\widetilde{q_{1j}}=q_{1j} + t^{e-e_2}a_2q_{2j} + \dots + t^{e-e_k}a_kq_{kj}.
\end{equation}  
We also define $q'_{ij}$ by the equality
$$q_{ij}=t^{e_i}c_{ij} + q'_{ij}.$$

We claim that $\widetilde{q_{11}}$ has degree  $<e$: We have
\begin{equation}
 \widetilde{q_{11}}=\sum_{i=1}^k t^{e-e_i}t^{e_i}a_i^{\vfi^{e_i}}c_i + t^{e-e_i}a_iq'_{i1} = \sum_{i=1}^k t^ea_i^{\vfi^{e_i}}c_i + t^{e-e_i}a_iq'_{i1}.
\end{equation}
Since each $q'_{i1}$ has degree $<e_i$, each  $t^{e-e_i}a_iq'_{i1}$ has degree $<e$. Hence the coefficient $t^e$ in $\widetilde{q_{11}}$ is 
$\sum_i a_i^{\vfi^{e_i}}c_i=0$ and $\deg (\widetilde{q_{11}})<e$.

Let $\tilde{Q}$ be the matrix where we have replaced $q_{1j}$ by  $\widetilde{q_{1j}}$. Now the sum of degrees of the non zero entries of $\tilde{Q}$ is less than the sum of degrees of the non zero entries of $Q$. Hence in order to apply an induction, it is enough to prove that the solvability of the system
$\bar{x}=\bar{y}.Q $ is equivalent to the solvability of $ \bar{x}=\bar{z}.\tilde{Q}$. 
But this follows by expressing the equations (\ref{eq:matricepassage}) by the equality
$$\tilde{Q}=PQ$$
where $$P=\begin{pmatrix} 
1 & t^{e-e_2}a_2 & t^{e-e_3}a_3 & \cdots & t^{e-e_k}a_k \\
0 & 1            & 0            & \cdots &       0       \\
\vdots & 0 & \ddots \\
0 & \cdots & \cdots & \cdots & 1
\end{pmatrix}					  $$ 
is clearly invertible in $R$.

{\it Claim 2}: Assume Claim 1. We may assume that the polynomials 
$(q_{i1})_i$ of the first  column have same degree and leading coefficients
of the $(q_{i1})_i$ are $K^{\vfi^e}$ linearly independent.

We will show that we can change $Q$ to some $S$, possibly having more rows, such that 
the system $\bar{x}=\bar{y}.Q$ is equivalent to $\bar{x}=\bar{z}.S$ with the first column of $S$ has required properties.

Recall that $q_{i1}$ has degree $e_i$ with $e_1=e$. For all  $1\leq i \leq k$,  using the basis $\alpha(e-e_i)$, we write the equality of terms 
	$$y_i=\sum_{u \in d^{e-e_i}}\lambda_u(y_i).t^{e-e_i}\alpha_u .$$
Now $$y_i.q_{i1}=y_i.t^{e_i}c_i + y_i.q'_{i1} = \sum_u \lambda_u(y_i).(t^e\alpha_u^{\vfi^{e_i}}c_i + r_{i1}(u))$$
where $r_{i1}(u)= t^{e-e_i}\alpha_uq'_{i1}$ which has degree $<e$. Set $z_i(u)=\lambda_u(y_i)$ and 
$$s_{i1}(u):=t^e\alpha_u^{\vfi^{e_i}}c_i + r_{i1}(u).$$ 
Since the leading coefficient of $s_{i1}(u)$ is $\alpha_u^{\vfi^{e_i}}c_i$, for any $i$ and $u$, 
$s_{i1}(u)$ has the degree  $e$.

For $j>1$ and $u\in d^{e-e_i}$, set $s_{ij}(u):=t^{e-e_i}\alpha_u.q_{ij}$.   
Note that we keep $e_i=e_{i1}$ but $j$ varies. Let $S$ be the matrix obtained from $Q$, by replacing 
$i$-th row by the matrix $(s_{ij}(u))_{u,j}$ where $u$ is the row-index and $j$ is the column index.
Then the system $(x_j=\sum_i y_i.q_{ij})_j$ is equivalent to the system
$(x_j=\sum_i \sum_u \lambda_u(y_i).s_{ij}(u))_j$, which can be written as
$$\bar{x} = \left((\lambda_u(y_1))_u,\ldots,(\lambda_u(y_k))_u \right).S .$$

Now we will show that the leading  coefficients $\alpha_u^{\vfi^{e_i}}c_i$'s are 
$K^{\vfi^{e}}$-linearly independent.

Suppose $$\sum_i \sum_u a_{u_i}^{\vfi^e}\alpha_u^{\vfi^{e_i}}c_i= 0$$ for some
tuple $(a_{u_i})_{u,i}$ from $K$. It follows   by Claim 1 that, 
$\sum_u (a_{u_i}^{\vfi^{e-e_i}}\alpha_u)=0$ for each $i$. Since the $\alpha_u$ are $K^{\vfi^e}$-linearly independent and $\vfi$ is injective
    $a_{u_i}=0$ for all $i,u$. The Claim 2 is  proved.

We assume now that the first column of $Q$ consists of 
polynomials of degree $e$ with leading coefficients being $K^{\vfi^e}$-linearly independent.
By section 3  and by the last paragraph of Lemma 4 of \cite{vdd-kuhlmann}, there exists an invertible matrix over $P$ over $K$
such that $P.Q_1$, where $Q_1$ the first column of $Q$, consists of polynomials with leading coefficients has all different valuations in 
$\{0,\ldots,d^{e}-1\}$. Hence considering $P.Q$ finishes the proof.
\end{proof}

For the following lemma, we will use the above result with $Q$ a column matrix and then in the following corollary we will use it in the whole generality.

\begin{lemma}\label{lem:pseudocompbeta}
Let $A \subseteq M$, of the form $A=\sum_i M.q_i$ then for some integer $s$, $A$ has a pseudo-complement of the form 
$\oplus_{i\in I} M.t^sb_i$ where $I\subseteq d^s$, and the $b_i$ are valuation independent.  
\end{lemma}
\begin{proof}
Use the above lemma to chose $g_j$ such that $\sum_j M.g_j =A$,  all of degree $s$ with leading coefficients $b_j$'s have different valuations in $\{0,\dots,d^s-1\}$. Then
by the lemma \ref{lem:addtoMta}, $M.g_j \approx M.t^sb_j$ for all $j$. Complete the $b_j$ to a valuation independent basis of $K^{\vfi^s}$-vector space $K$. We write the new elements of this basis as  the $b_i$.    Let $C:=\sum_{i} M.t^sb_i$. 
Since $C\oplus \sum M.t^sb_j=M$, $C$ is in particular a pseudo-complement for $\sum M.t^jb_j$. Hence by Lemma
\ref{changecomp}, $C$ is a pseudo-complement for $A$. 
\end{proof}
\begin{remark}\label{rem:pseudocomplements} 
 The pseudo-complement $C$ is $p.p.$ definable  subgroup by an $L$-formula $C(x)$, which does not depend on $M$. In other words, if $\phi(x)$ is the formula 
$$\exists y_1, \dots , y_m \mesp x=\sum_i y_i.q_j$$ then in any valued module $(M,v)$,
$C(M)$ is a pseudo-complement to $\phi(M)$.

\end{remark}

\begin{theorem}\label{pseudocomplements}
Let $\phi(\bar{x})$ be a $p.p.$ formula of $L_\mathcal{O}$ of the form 
$$\bar{x} - \bar{y}.Q \in W,$$ where $Q$ a matrix with coefficients from $R$ and $W$ a ball. Set $k:=\vert \bar{x} \vert$. Then there exists computable $\gamma$ and an existential $L$-$p.p.$ formula $D(\bar{x})$ such that 
$$M^k=\phi(M^k) + D(M^k) + P_\gamma^k$$ $\&$
$$D(M^k)\cap \phi(M^k) \subseteq P_\gamma^k$$
for all valued module $(M,v)$.
\end{theorem}
\begin{proof}

 By Lemma \ref{vddku}
we may assume that the first column $Q_1$ of $Q$ consists of polynomials having same degree $s$ with leading coefficients having distinct valuations in $\{1,\dots, d^s-1\}$. In addition, by Remark \ref{rem:+ballpseudocomp},
we may suppose that $W=\{0\}^k$.

 Take $C$ a pseudo-complement to $M.Q_1$ as in the above lemma. 
Let $x=(x_1, \dots, x_k)$ such that $x_1 \in C \cap M.Q_1$. Then $x_1 \in 
P_\delta$ for some computable $\delta$. Write $x_1$ also as
$$x_1=y_1.q_{11} + y_2.q_{21} + \dots y_mq_{m1}.$$
Since the leading coefficients of the $q_{1j}$ have different valuations
in  $\{0,\dots, d^s-1\}$,  if 
$$v(y_j)\leq \rho:=\min\{\min_{ij} \{\jump{q_{ij}}\}\}-1;$$ we have $v(y_j.q_{1j})\neq v(y_{j'}.q_{1j'})$ for $j\neq j'$. Hence the $y_j$ can not have indefinitely small valuations since $v(x_1)\geq \delta$. Now  for $i>1$, since $x_i=\sum_{j}y_j.q_{ij}$ we have 
$v(x_i)\geq \min_j \{\rho\cdot q_{ij}, \delta\}$. Setting $\gamma:=\min_{ij}\{\rho\cdot q_{ij}\}$ and $D:=C\times M^{k-1}$ yields our claim since $D$ and $\gamma$ depends only to $Q$.
\end{proof}

\begin{remark}\label{boundedsolutions}
The above proof shows that whenever a matrix $Q$ with the first column  consists of polynomials whose dominant coefficients are valuation independents and $\gamma \in \Gamma$, as above are given, there is a computable $\delta \in \Gamma$, such that   
$$\bar{y}.Q \in P_\gamma^m \Rightarrow \bar{y} \in P_\delta^k  .$$  
\end{remark}
\subsection*{\bf  The theory $T_{\Psi}$}

Let $(M,v)$ be a valued module and set 
    $$A:=\sum M.q_i$$
  where the $q_i$ are all of degree $s$ such that for a  valuation independent basis $\beta$  of $K^{\vfi}$-vector space $K$, the leading coefficients of the $q_i$ are from the basis $\beta(s)$. Let $\mathbb{A}(\beta)$ be the tree consisting of subgroups $M.t^i\beta_i(n) \, (i,n \in \omega)$, ordered by inclusion. At the level $n$ we have the subgroups 
$M.t^i\beta_i(n) \mesp (i \in d^n)$. 

\begin{defn}
	We call the pseudo-complement for $A:=\sum_j M.q_j$ 
	the unique pseudo-complement which can be written as the sum of some elements of level
	$\deg(q_i)$ of the tree $\mathbb{A}(\beta)$.
\end{defn} 
This definition only depends on  the decomposition 
of the leading coefficients of the $(q_i)_i$ in the basis $\beta(n)$. Hence, given $A$ as above, by Lemmas \ref{vddku} and \ref{lem:pseudocompbeta}, there is a recursive function 
\begin{equation}\label{metaf}
\mathfrak{f}: \bigcup_n K^n \to \bigcup_n K^n
\end{equation} 
which computes the basis $\beta$ and the pseudo-complement $C$ of  in every valued 
$t\beta$-decomposable $R$-module $M$. We  write 
$C_{\mathfrak{f}(q_1,\ldots,q_n)}$ for the pseudo-complement computed by $\mathfrak{f}$.

On the other hand, the ball $P_\gamma$ such that $M=A+C + P_\gamma$, can be chosen 
by letting $\gamma:=\min\{\min_i\{\jump{q_i}\} \} - 1$. This yields another recursive function
\begin{equation}\label{metaj}
\mathfrak{j}:  \bigcup_n K^n \to \Gamma.  
\end{equation} 

Hence we can express the statement of Corollary \ref{pseudocomplements} by introducing an axiom $\psi(q_1,\ldots, q_n)$  for any matrix $Q$ which has the first column  $(q_,\ldots,q_n)$, in the language 
$L_{\mathcal{O}}$, which says that 
 \begin{equation}
 M^k=M^k.Q + C_{\mathfrak{f}(q_1,\ldots,q_n)}\times M^{k-1} + P_{\mathfrak{j}(q_1,\ldots,q_n)}^k \mesp \quad \text{and}
\end{equation}
$$ M^k.Q \cap  C_{\mathfrak{f}(q_1,\ldots,q_n)}\times M^{k-1} \subseteq P_{\mathfrak{j}(q_1,\ldots,q_n.)}^k.$$

 Let $T_\Psi$ be the $L_\mathcal{O}$-theory of  $R$-modules together with the sentences  $\psi(q_1,\ldots, q_n)$. Hence
 $T_\Psi$ is recursively enumerable.

\begin{theorem}\label{approx}
Let $Q$ be a  $m\times k$  matrix over $R$, $W=\prod_{i=1}^k P_{\gamma_i}$  and $\phi$ 
	be the $L_{\mathcal{O}}$-formula
$$\phi(x_1,\ldots,x_k): \exists y_1\ldots y_m \mesp (x_1,\ldots x_k) - (y_1,\ldots, y_m).Q \in W.$$  Then for some  computable $\gamma$, and some positive primitive  $L_\mathcal{O}$-formula $D(x)$ we have  
$$M^k =  \phi(M^k) + D(M^k) + \prod_{i=1}^k P_\gamma   \quad \& \quad D(M) \cap \phi(M^k)\subseteq \prod_{i=1}^k P_\gamma.$$ 
for all $M\models T_{\Psi}$.
\end{theorem}

\section{Decidability and model completeness of $\mathbb{F}_d((X))$}
We will introduce a new theory $T_1$, containing $T_{\Psi}$, augmented by sentences counting the number of solutions of the $p.p.$ formula 
$$\bar{x}-\bar{y}.Q \in B$$ in $B_0$ modulo $B_1$, for proper balls $B_1 \subseteq B_0$.

For  $Q$, an $m\times n$ matrix over $R$, such that its first column consists of 
polynomials whose leading coefficients are valuation independent and $B$ a ball, set $$A(\bar{x}):= \exists  \bar{y} \mesp \bar{x}-\bar{y}.Q \in B.$$ Let $\delta$ be the value computed by 
Theorem \ref{boule_up}, such that $A\wedge P_\delta^m$ is quantifier-free definable
in language  $L(\lambda)$. 

   Let    $B_0:=P_\gamma^m$, where $\gamma$ is given by Theorem \ref{approx}.  $B_1:=P_\delta^m$,  and 
   $$k:=\vert (A(K)\cap B_0(K))/A(K)\cap B_1(K)\vert$$ if $\gamma\leq \delta.$ We set the sentence  $\theta(A,\gamma,\delta)$ expressing 
    \begin{equation}
     k=  \left| \frac{(A\wedge B_0)(K)}{(A \wedge B_1)(K)}\right| .
    \end{equation}
Let  $\Theta$ the be set of sentences  $\theta(A,\gamma,\delta)$ and we set 
 $$T_1:=T'_{hens}\cup T_{\Psi} \cup \Theta$$
 where $T'_{hens}$ is the $L_\mathcal{O}$-theory composed by the axioms of $T_{hens}$, where we have replaced any $L(\lambda)$-term by its equivalent modulo the theory of $R$-modules in language $L$ (recall Remark \ref{extunideflambdas}). 

 Note that $T_1$ implies that $P_{\gamma +1}/P_{\gamma}$ has exactly $d$ elements for all $\gamma \neq \infty$.
\begin{remark} $\Theta$ and hence $T_1$ is a recursively enumerable theory. In fact, by Remark \ref{boundedsolutions}, if $\bar{x},\bar{y},\gamma$ are such that
	$\bar{y}.Q=\bar{x}\in P_{\gamma}^n$,  then $\bar{y}$ is in some $P_{\theta}^k$ for a computable $\theta$, hence searching the solutions 
	$\bar{y}$, of $\bar{y}.Q=x$ with $\bar{x} \in P_\gamma^n$ can be bounded to searching $y$'s in some ball. Hence
	searching such solutions modulo another ball can be done in some finite $\mathbb{F}_d$-vector space by an algorithm.
\end{remark}
\begin{prop}\label{laurentsatisfyT1}
Let $M:=\mathbb{F}_d((X))$ then as an $L_\mathcal{O}$-structure $M\models \Theta$  and hence $M \models T_1$.
\end{prop}
\begin{proof}
By \cite{kuhlmanntame}, Theorem 5.14, $K$ is existentially closed as a ring  in $M$. In addition, by \cite{sylvyarno17} Corollary 6.18, there exists an  existential ring-formula
without parameters which defines the maximal ideal both in $K$ and in $M$. Since the valuation ring is the complement of the set of inverses of the elements in the maximal ideal, we have a universal ring-formula  which defines uniformly the valuation ring in $K$ and in $M$. Hence the balls 
centered at $0$ are definable  universally with the parameter  $X$ both in $K$ and $M$.

Suppose $\vert (A(K)\cap B_0(K))/(A\cap B_1(K))\vert = k$. Consider the sentence $$\sigma: \vert (A\wedge B_0)/(A\wedge B_1)\vert \geq k.$$
Since by theorem \ref{boule_up}$, (A(K)\cap B_1(K))$ is definable both universally
and existentially in $L_{\mathcal{O}}$,  and and the valuation ring is universally definable both in $K$ and $M$,  $\sigma$ is equivalent
to an existential ring-formula with parameters in $K$. Since this quotient is finite, we must have $\vert (A(M)\cap B_0(M))/A(M)\cap B_1(M)\vert = k$.
\end{proof}

To prove Theorem \ref{th:pp2universal}, we will use  a lemma from Rohwer's thesis 
(\cite{rohwer}, Lemma 8.2). This lemma is a generalization of the following fact:

\begin{obs} In an abelian group $G$ with existentially definable subgroups $A,B$ such that $A+B=G$, if $A\cap B$ is definable by a universal formula then $A$  is definable by the following universal formula $\psi(x)$:

\begin{equation}
\psi(x): \forall y \mesp (x-y \in A \wedge y \in B) \rightarrow y \in A\wedge B.
\end{equation}
By iterating this observation,  we have the following lemma.
\end{obs}
\begin{lemma}[Rohwer]
	Let $\mathcal{T}$ be a theory expanding  the theory of abelian groups. For each
	$M \models \mathcal{T}$ and for  $A, A_c, A_m, A_s, B_0, B_1$  definable subgroups of $M$ satisfying  the following configuration,
\end{lemma}
\begin{enumerate}[1.]
\item $A + A_c=M$,
\item $A\cap A_c \subseteq B_0$,
\item $A\cap B_1 = A_s\cap B_1$,
\item $A\cap B_0 \subseteq A_m \subseteq A + B_1$.
\end{enumerate}\medskip
where $A,A_c,B_1$ are definable by existential formulas, and $A_m, A_s$ 
by universal formulas (where all formulas in question do not depend on $M$), $A$ is definable by a universal formula (which does not depend on $M$).   
\begin{proof} See \cite{rohwer} Lemma 8.2. 
\end{proof}

\begin{theorem}\label{th:pp2universal}
Any $p.p.$ formula of $L_{\mathcal{O}}$ is equivalent modulo $T_1$ to a universal $L_\mathcal{O}$-formula.   
\end{theorem}\label{modcomplete}
\begin{proof}
 Let $$A(\bar{x}): \exists \bar{y} \mesp \bar{x}.S -\bar{y}.Q \in W$$ be a $p.p.$ formula. For our purposes, we may assume that 
$S$ is the identity matrix since me way replace $A(x)$ by
	$$\forall \bar{z} \mesp (\bar{z}=\bar{x}.S \to \exists \bar{y} \mesp \bar{z}-\bar{y}.Q \in W).$$  
	By Lemma \ref{vddku}, we may suppose that the first column of $Q$ consists of polynomials which have coefficients
	in a valuation independent basis. 
	By Theorem \ref{boule_up}, there is some proper ball $B_1$ determined by $T_1$ such that the formula $A_s:=A\wedge B_1$ is equivalent 
to  a quantifier free $L(\lambda)$-formula. Hence by Remark \ref{extunideflambdas}, $A_s$ is equivalent to a universal
$L$-formula. Now, let  $B_0$ be given by Theorem \ref{approx} such that
$$A \wedge A_c \to B_0$$ where $A_c$ is of the form $A_c=D+P_\gamma^l$ with $D$    as in Theorem \ref{approx}.

\emph{Claim} :  Set $A_m:=(B_0 \wedge  A) + A_s$. Then $A_m$ is equivalent to a  universal 
$L_{\mathcal{O}}$-formula modulo $T_1$.

\underline{proof of the claim}:
 Let $M \models T_1$. Let $k$ be the cardinality of $$(B_0(M) \cap A(M))/
 A_s(M)$$ and $y_1, \dots, y_k$ ($y_i$ are tuples of variables)  be  representatives of the classes.
Note that $k$ is determined by a sentence 
in $\Theta$, hence depends only on the theory $T_1$.
	Then $$\{y_1, \dots y_k \}  + A_s(M) = A(M)\cap B_0(M)+A_s(M).$$ Moreover,
	for all $z_1,\dots,z_k \in B_0$ satisfying $z_i-z_j \notin A_s(M)$ for $i\neq j$,
	$$\{ z_1, \dots z_k \}  + A_s(M) = A(M)\cap B_0(M) + A_s(M)$$
if and only if,  
	$$\{z_1, \dots z_k \}  + A_s(M) \subseteq A(M) \cap B_0(M) + A_s(M).$$ Hence 
the formula 
\begin{multline}
	x\in B_0 \wedge \forall y_1, \dots, \forall y_k \\ (
  (\bigwedge_{i=1}^k y_i \in B_0 ) \wedge ( \bigwedge_{i\neq j} y_i - y_j \notin A_s) \wedge  \\
	[\exists_{i=1}^{k}z_i (\mesp \bigwedge_{i=1}^k z_i\in A \wedge \bigwedge_{i=1}^k y_i - z_i \in A_s)]) \longrightarrow \bigvee_{i=1}^k y_i -x \in A_s.
\end{multline}
is equivalent to $A_m$, which is equivalent to a universal $L_\mathcal{O}$-formula, thus the claim is proved.

Now $A,A_c,A_m,A_s,B_0,B_1$ are in the configuration of Rohwer's Lemma above.
\end{proof}

\begin{corr}\label{decidable}
	Every completion of $T_{1}$ is model-complete in language $L_{\mathcal{O}}$.  In particular, the complete $L_\mathcal{O}$-theory of $\mathbb{F}_d((X))$ is model-complete.
\end{corr}
\begin{proof}
	It follows by Corollary \ref{corr:baur-monk}.
\end{proof}
Note that 
$K$ embeds  (via an $L_{\mathcal{O}}$-embedding) to any model $N\models T_{1}$. In fact, choose any $s \in N\setminus\{0\}$ such that $s.t=s$. Consider 
the $L$-embedding $k\mapsto s.k$. We show that it is an $L_{\mathcal{O}}$-embedding:  It is easy to see that the assertion
$$ \mathcal{O} = A \oplus \mathcal{O}.X$$ where  $A$ 
is the formula $x.t=x$, is a consequence of $T_{1}$ and it is clear that 
if $a\in \mathcal{O}_K$ then $s.a \in \mathcal{O}_N$. Suppose now that 
$s.a \in \mathcal{O}_N$ for some $a \in K$. Then
$s.a=y  +  m$ where $y.t=y$ and $m\in \mathcal{O}_N.X$. Then
$(s.a).(t-1)=m.(t-1) \in \mathcal{O}_N.X$. Note that 
$(s.a).(t-1)=s.ta^{\vfi} -sk=s.(a^{\vfi} -a) \in \mathcal{O}.X$ since $s.t=s$.
This can only happen if $(a^{\vfi} -a) \in \mathcal{O}_K$, only if 
$a\in \mathcal{O}_K$.

\begin{corr}\label{K}
The models of $T_1$ in which $K$ is existentially closed as an $L_\mathcal{O}$-structure are elementary equivalent to $K$. 
\end{corr}

\begin{proof}

	Let $N\models T_{1}$ in which $K$ is existentially closed. 
 Let $A, B$ be $p.p.$ formulas (of the language $L_{\mathcal{O}}$) with one free variable  such that  
 $$T_1\models \forall x \mesp A_1(x) \rightarrow A(x).$$
	
Set $k:=\vert A(K)/A_1(K) \vert$. It is enough to show that 
 we have 
	$$\vert A(N)/A_1(N) \vert =k.$$ 
	
Let $l>k$ and suppose
$$M \models \exists x_1,\ldots \exists x_l \mesp (\bigwedge_i x\in A)\wedge
	(\bigwedge_{i\neq j} x_i-x_j \notin A_1).$$ 
		Since modulo $T_{1}$ the formula $x \in A_1$ is equivalent to a universal formula, its negation is equivalent to an  existential $L_{\mathcal{O}}$-formula. Since
		$K$ is existentially closed in $N$, there exists at least $l$ element in $A(K)/A_1(K)$. Contradiction.

\end{proof}

In particular, we have:
\begin{corr}
$K$ is the prime model of the complete theory of $\mathbb{F}_d((X))$.
\end{corr}

\begin{proof}
$K$ is existentially closed in $\mathbb{F}_d((X))$ as an $L_{\mathcal{O}}$-structure since it is existentially closed in $\mathbb{F}_d((X))$ as a ring. Hence $$K\equiv \mathbb{F}_d((X)).$$ Since any completion of $T_1$ is model-complete $K$ is an $L_{\mathcal{O}}$ elementary substructure of $\mathbb{F}_d((X))$.
\end{proof}

\subsection*{Decidability of $\mathbb{F}_d((X))$}

The following fact is an easy exercise using Hensel's lemma:
\begin{fait}
	$\mathbb{F}_d[[X]]$ is definable by the $L$-$p.p.$ formula 
	$$ \exists y \mesp x.tX=y.(t-1)$$ inside $\mathbb{F}_d((X))$. 
\end{fait}

Hence, the decidability of $\mathbb{F}_d((X))$ as an 
$L_\mathcal{O}$-structure and as an $L$-structure are equivalent. Recall that the $L$-theory of $K$ (hence of $\mathbb{F}_d((X)))$, are given by the sentences
stating that 
$$\vert A(K)/(A(K)\wedge B(K)) \vert =k \quad (k\in \mathbb{N}\cup \{\infty\})$$
where $A, B$ are $L$-$p.p.$ formulas of  
with one free variable.

We will show that  a recursively enumerable subset of these sentences forms 
a complete axiom system which implies all of them. Hence the $L$-theory of $K$ is decidable.

For $A$ and $B$ as above, we set $D:=A\wedge B$. Let 
$A$ be given by $\exists \bar{y} \mesp x.p=\bar{y}\bar{q}$ and 
$B$ by $\exists \bar{y} \mesp x.r=\bar{y}\bar{s}$.  We may suppose that 
both $p$ and $r$ are unitary. Hence, by Lemma \ref{lem:addtoMta}, for  $k=\deg(p)$ and $s=\deg(r)$  we have
$$K.p \approx K.t^k$$ and $$K.r \approx K.t^s .$$


Now by Remark \ref{rem:samecomplement}, \ref{rem:invimage} and Lemma \ref{changecomp}, $A\approx D$ if and only if  the preimage by $.t^k$ of $C_{\mathfrak{j}(\bar{q})}$(which is a pseudo-complement of $A$) is equal to
the preimage  of $C_{\mathfrak{j}(\bar{s})}$ by $.t^s$ (which is a pseudo-complement of $D$). Hence it is decidable if $A\approx D$.

Notice that  If  $A \not \approx D$ then
$A/D$ is infinite:  Let $a\in A \setminus D$.  Then for  $a' \in A$ with $v(a')<v(a)$  either $a'$ or $a-a'$ is not in $D$. 

Now suppose $A\approx D$. It follows by Remark \ref{rem:modulogamma} that
$$ A+P_\alpha=D+P_\alpha $$  for some computable $\alpha$: $\alpha$ can be chosen less than every jump values of all polynomials appearing in the definition of $A$ and $B$. 
Hence $\vert A/D \vert = \vert A\cap P_\alpha / D\cap P_\alpha\vert $.

Now by Corollary \ref{corr:boule_up}, there is a computable $\gamma$ such that $A/D$ is finite if and only if
 $$A\cap P_\gamma = D \cap P_\gamma$$

Consider the following algorithm: given $A$ and $B$  setting $D=A\wedge B$, the algorithm check if $A\approx B$, if not it  sets
$\vert A/D \vert =\infty$. Otherwise computes $\alpha $ such that $A+P_\alpha=B+P_\alpha$.  Then it computes the value $\gamma$
 so that it can check whether $A\cap P_{\gamma}/D\cap P_{\gamma}$  is trivial or infinite. If it is infinite it sets 
 $\vert A/D \vert = \infty$. If not, by Remark \ref{boundedsolutions} the algorithm can compute the number of the elements of 
  \begin{equation}\label{t0}
  \frac{(A\cap P_\alpha)/P_{\gamma}}{(D\cap P_\alpha)/P_{\gamma}}.
  \end{equation}
  
  Hence the $L$-theory consisting of sentences of type $(\ref{t0})$ can be
  recursively enumerable and implies the $L$-theory of $K$, hence of $\mathbb{F}_d((X))$. We have proved:
  \begin{theorem}\label{thm:finaldecidable}
  	Both the $L$-  and the $L_{\mathcal{O}}$-theories of $\mathbb{F}_d((X))$ are decidable.
  \end{theorem}
  
  \bibliographystyle{amsalpha}
  \bibliography{onayg}
  
\end{document}

\documentclass[10pt]{amsart}
\usepackage[utf8]{inputenc}
\usepackage{enumerate}
\pdfoutput=1

\usepackage{enumerate}
\usepackage{amsmath, amsthm}
\usepackage{amsfonts}
\usepackage{amssymb}
\usepackage{graphicx}
\usepackage{url}

\usepackage[left=3cm, right=3cm]{geometry}
\usepackage{array}
\makeatletter
\newcolumntype{e}[1]{
	>{\minipage[t]{\linewidth}%
		\NoHyper
		\let\\\tabularnewline
		\enumerate
		\addtolength{\rightskip}{0pt plus 50pt}
		\setlength{\itemsep}{-\parsep}}%
	p{#1}%
	<{\@finalstrut\@arstrutbox\endenumerate
		\endNoHyper
		\endminipage}}

\newcolumntype{i}[1]{
	>{\minipage[t]{\linewidth}%
		\let\\\tabularnewline
		\itemize
		\addtolength{\rightskip}{0pt plus 50pt}%
		\setlength{\itemsep}{-\parsep}}%
	p{#1}%
	<{\@finalstrut\@arstrutbox\enditemize\endminipage}}

\AtBeginDocument{%
	\@ifpackageloaded{hyperref}{}%
	{\let\NoHyper\relax\let\endNoHyper\relax}}
\makeatother

\AtEndDocument{\bigskip{\footnotesize%
		\textsc{ Instit\"{u}t f\"{u}r Mathematishe Logik und
			Grundladenforchung Fachbereich und Informmatik 
			Einsteinstra{\ss}e 62, 48149 M\"{u}nster, Deutschland
		} \par  
		\textit{E-mail address}, G\"{o}nen\c{c} Onay  \texttt{onay@uni-muenster.de} }}
\renewcommand{\epsilon}{\varepsilon}
\renewcommand{\le}{\leqslant}
\renewcommand{\ge}{\geqslant}
\renewcommand{\leq}{\leqslant}
\renewcommand{\geq}{\geqslant}

\newcommand{\jump}[1]{\operatorname{Jump}(#1)}
\newcommand{\lcm}{\operatorname{lcm}}
\newcommand{\hens}{\operatorname{hens}}
\newcommand{\h}{\operatorname{h}}
\newcommand{\rv}[1]{\operatorname{rv}(#1)}
\newcommand{\nqrt}[2]{\sqrt[\varphi^{#1}]{#2}}
\newcommand{\ssi}{\Leftrightarrow}
\newcommand{\implique}{\Rightarrow}
\newcommand{\vers}{\rightarrow}
\newcommand{\tq}{\;\: | \:\:}
\newcommand{\mesp}{\;\:}
\newcommand {\vfi}{\varphi}
\newcommand{\mr}[2]{\text{M}_{R}(#1,#2)}
\newcommand{\iffl}{\longleftrightarrow}
\newcommand{\iffs}{\leftrightarrow}

\newcommand{\Nn}{\mathbb{N}}
\newcommand{\Zz}{\mathbb{Z}}
\newcommand{\Qq}{\mathbb{Q}}
\newcommand{\Rr}{\mathbb{R}}
\newcommand{\Cc}{\mathbb{C}}
\newcommand{\Kk}{\mathbb{K}}
\newcommand{\Hh}{\mathbb{H}}
\newcommand{\F}{\mathbb{F}}
\newcommand{\K}{\mathbb{K}}

\newtheorem{theorem}{Theorem}[section]
\newtheorem{lemma}[theorem]{Lemma}
\newtheorem{prop}[theorem]{Proposition}
\newtheorem{proposition}[theorem]{Proposition}
\newtheorem{corr}[theorem]{Corollary}
\newtheorem*{intth}{Theorem}
\newtheorem*{intcor}{Corollary}
\newtheorem{defn}[theorem]{Definition}
\newtheorem{fait}[theorem]{Fact}

\newtheorem{remark}[theorem]{Remark}
\theoremstyle{definition}
\newtheorem{exemple}[theorem]{Example}
\newtheorem{obs}[theorem]{Observation}
\newtheorem{nota}[theorem]{Notation}

\author{G\"{o}nen\c{c} Onay}   
\thanks{Research partially founded by DFG through SFB 878.}

\makeatletter
\def\blfootnote{\gdef\@thefnmark{}\@footnotetext}
\makeatother
\title{$\mathbb{F}_p((X))$ is decidable as a module over the additive polynomials
}

\begin{document}

	\sloppy
	
	\begin{abstract}
		Let $R$ be the (non commutative-) ring of additive polynomials over the field $K:=\mathbb{F}_p(X)^h$, the henselization of the
		field $\mathbb{F}_p(X)$. We show that the (right-) $R$-module theory of the field  $\mathbb{F}_p((X))$ is decidable. Moreover, we provide a recursively enumerable axiom
		system $T_1$ (satisfied by $\mathbb{F}_p((X))$)
		in the language $L_{\mathcal{O}}$, the language of $R$-modules together with a predicate $\mathcal{O}$ for the 
		valuation ring $\mathbb{F}_p[[X]]$, and show that every positive primitive formula is equivalent to a universal formula modulo $T_1$. As an $L_{\mathcal{O}}$-structure,  $\mathbb{F}_p((X))$ is also decidable and the   $L_\mathcal{O}$-theory of  $\mathbb{F}_p((X))$ is model-complete admitting   $K$ as its prime model.
	\end{abstract}
	\maketitle
	\blfootnote{\textup{2010} \textit{Mathematics Subject Classification}: 11U05 (primary), 12L05, 12J10	03C68 (secondary)}

	\section{Introduction}
	
	The decidability and the axiomatization of the field of Laurent
	series over the finite field $\mathbb{F}_p$ is a longstanding 
	problem. On the other hand, its characteristic $0$ analog $\mathbb{Q}_p$ is  axiomatized and is decidable by the Ax-Kochen and Ershov Theorem. As a consequence the similarities between $\mathbb{Q}_p$ and $\mathbb{F}_p((X))$ can be expressed in terms of ultralimits when $p \to \infty$.  
	The first and maybe the most famous application is to give a  corrected form of a conjecture by Artin. This corrected form states that 
	given any integer $d>0$ there is some integer 
	$m$, such that for any prime $p>m$, every homogeneous polynomial of degree $d$ over $\mathbb{Q}_p$ with $>d^2$ variables has a non-trivial zero in $\mathbb{Q}_p$.
	This application uses simply that fact that for every prime $p$ the field $\mathbb{F}_p((X))$  is $C_2$ 
	(i.e. every homogeneous polynomial with strictly more
	variables than the square of its degree has a non-trivial zero)\footnote{The original conjecture of Artin was claiming that $\mathbb{Q}_p$ was $C_2$ for all prime $p$, this is refuted by Ternejian.}.

	Despite the (asymptotic -) analogy between $\mathbb{Q}_p$ and $\mathbb{F}_p((X))$, for a fixed prime $p$ 
	the field theory of $\mathbb{F}_p((X))$ remains  unknown. Moreover,  Kuhlmann proved that the naive translation into positive characteristic of the 
	complete theory of $\mathbb{Q}_p$ is incomplete (see \cite{ku-elementary}, we will give more details on this in the following lines).
	Denef and Schoutens showed that its existential theory in the language of rings with a constant symbol for $X$ is decidable assuming  resolution of singularities in positive characteristic (see \cite{denefsch}).  Recently Anscombe and Fehm showed unconditionally that its existential theory in the language of rings is decidable \cite{sylvie-fehm}.

	After the works of Kuhlmann and van den Dries, incompleteness of this naive theory  can be expressed 
	using only  properties of   additive polynomials, hence in  the language of $S$-modules  where $S$ 
	is the ring of  additive polynomials over $\mathbb{F}_p(X)$.
	We believe that Rohwer, in his thesis (see \cite{rohwer}), shows that the complete theory of $\mathbb{F}_p((X))$ as an $S$-module is model-complete 
	in the language of $S$-modules together with a predicate for the valuation ring\footnote{However the meaning of Conjecture A.4. remains unclear for us.}. 
	However, he does not provide an axiom system for this theory and does not show its decidability. We think that our philosophy is similar to Rohwer's one but for our results we have had been rather inspired by  the articles \cite{ddp1}, \cite{vdd-kuhlmann} and \cite{sylvie-fehm}.
	Notice that  B\'elair and Point also studied the module theory of
	some (valued -) fields satisfying  strong divisibility conditions (see
	\cite{belair-point10} and \cite{belair-point15}).

	\subsection*{Main results}
	Let $K$ be the henselization of the field of rational functions over $\mathbb{F}_p$ and  $$\vfi:K\to K,\mesp x\mapsto x^{p^k}$$ 
	be some (fixed-) power of the Frobenius map. We set $R$ to be the ring of $\vfi$-polynomials, that is, additive polynomials 
	whose monomials are of the form $aT^{{p^k}^i}$, equipped with the
	composition and the usual addition. Let $M:=\mathbb{F}_p((X))$ be the field of Laurent series over $\mathbb{F}_p$, let $L$ be the language of (right-) $R$-modules
	and finally let $L_{\mathcal{O}}$ be the language $L$ together with the unary predicate $\mathcal{O}$ (for the valuation ring $\mathbb{F}[[X]]$).
	
	We prove in particular the following results in the present article.
	
	\begin{enumerate}
		\item There is a recursively enumerable $L_{\mathcal{O}}$-theory $T_1$, $M\models T_1$, such that any completion  is model complete. Moreover,  $K$ is the prime model of the complete $L_{\mathcal{O}}$-theory of $M$  (see Theorem \ref{th:pp2universal} and following corollaries). 
		
		\item Both the  complete $L$- and $L_{\mathcal{O}}$-theories of $M$ are decidable (see Theorem \ref{thm:finaldecidable}).
	\end{enumerate}

	\subsection*{Strategy and organization of the article}
	Following Kuhlmann (\cite{ku-elementary}), let us first explain why the naive adaptation of the theory
	of $p$-adics into positive characteristic is incomplete. We denote this theory by $T_{naive}$, written in 
	language of rings together with a unary predicate for the valuation ring 
	and a constant for the uniformizer $X$,

	We recall that $(U,v) \models T_{naive}$ if and only if
	\begin{enumerate}[-]\label{naivfpt}
		\item $(U,v)$ is a  henselian, non trivially valued  field of characteristic $p>0$,
		\item the residue field $U/v$ is $\mathbb{F}_p$,
		\item the value group $\Gamma$ is a $\mathbb{Z}$-group, 
		\item $(U,v)$ is defectless,
		\item $v(X)=\min \Gamma_{>0}$.
	\end{enumerate}
	
	Set $M=\mathbb{F}_p((X))$. We know that:
	\begin{equation}\label{as-formula}
		M=\bigoplus_{i=1}^{p-1}M^pX^i + \wp(M) + \mathcal{O}_M, \tag{*} 
	\end{equation}
	where $\wp$ is the Artin-Scherier map $x\mapsto x^p -x$  and $\mathcal{O}_M$ is the valuation ring of $M$.
	But the equality (\ref{as-formula}) fails to hold  for some extension  $N\supset M$, which is a  model of  $T_{naive}$. Hence $T_{naive}$ is incomplete.

	Consider the polynomial $F(z_0,\ldots,z_{p-1}):=z_0^p-z_0 +\sum_{i=1}^{p-1} z_i^pX^i$. From the equality above  and using Hensel's Lemma, one can deduce that for any $a$, the set 
	$\{v(F(x)-a) \tq x \in M\}$ has a maximum in $\mathbb{Z}\cup \{\infty\}$: we say that the
	image of $F$ has the {\it optimal approximation property} (see \cite{vdd-kuhlmann}). From this observation Kuhlmann suggests a candidate for a complete axiomatization
	of the theory of $M$ in the language of rings,  which is essentially $T_{naive}$ together with sentences
	saying that the image of every multi-variable additive polynomial has the optimal approximation property (see \cite{ku-elementary} and \cite{kuhlmanntame}).

	Now we present our approach. We work in the language 
	$L_{\mathcal{O}}$. 
	
	Let us illustrate in an example the main ideas of the proof of Theorem \ref{th:pp2universal}, which states that modulo the theory $T_1$
	every {\it positive primitive} formula is equivalent to a universal one. This theorem immediately shows that every completion of
	$T_1$ is model-complete (see Section \ref{sec:baur-monk} and Corollary \ref{corr:baur-monk}).

	Instead of considering the multivariate polynomial $F$ above, we use the equation (\ref{as-formula}) to study the image of the Artin-Scherier map $\wp$ and see the sum 
	$$C(M):=\sum_{i=1}^{p-1}M^pX^i$$ as {\bf the pseudo-complement}  to the image 
	$\wp(M)$. Note  that $$\wp(M)\cap C(M) \subseteq \mathcal{O}_M,$$ hence this intersection is {\it small} 
	with respect to the valuation metric; that is why we call $C$ the pseudo-complement. In addition, using Hensel's lemma we have that 
	$$\mathcal{O}_M \cap \wp(M)=\mathfrak{m}$$ where 
	$\mathfrak{m}$ is the maximal ideal $X\mathbb{F}_p[[X]]$.
	
	At this point, we are able to describe the set $\wp(M)$, definable a priori by an existential formula, by the following formula:
	\begin{equation}\label{elim-as}
		\phi(x): \forall y  \mesp ((\exists z \mesp \wp(z)=x-y\wedge y\in C + \mathcal{O}) \rightarrow y \in \mathfrak{m}).   \tag{**}
	\end{equation}
	It is easy to see that this formula is equivalent to a universal 
	formula in the language $L_\mathcal{O}$ (hence also in the language of rings with the parameter $X$) since
	$C$ is existentially definable.

	Motivated by the above example, we consider the general case. We introduce the notion of a ball (in this section) as an $L_\mathcal{O}$-formula of the form 
	$$B(x): x.X^\gamma \in \mathcal{O}$$ for some fixed $\gamma \in \mathbb{Z}$. Let $Q$ be a  finite set of  polynomials (scalars from $R$).  Abusively, we also denote by $Q$ the formula which defines the sum of the images of $q\in Q$.  Our strategy is to assign to $Q$ a positive existential  formula $D$ and balls $B_1,B_2$ in a computable way,  such that in every model $N\models T_1$,
	\begin{enumerate}[(i)]
		\item \label{task1} $D(N)  + Q(N)=N$ and $Q(N)\cap D(N) \subseteq B_2(N)$
		\item \label{task2} $B_1(N) \subseteq  B_2(N)$, such that  $B_1(N)\cap Q(N)$ is definable by a universal formula which depends only on $T_1$. 
	\end{enumerate}
	For instance, the above example suggests that for $Q=\{\wp\}$, 
	$B_1=\mathcal{O}$, $B_2=\mathcal{O}.X$ and 
	$D=C + \mathcal{O}$ suit. Then  prove that, modulo $T_1$, every positive primitive formula is equivalent to a universal one (similar to (\ref{elim-as}) above).
	
	In the general case, for the task (\ref{task2}),  we cannot always ensure that $B_1(N)\cap Q(N)$  is a ball, or more generally quantifier free definable in $L_\mathcal{O}$: to see this it is enough to consider the image of $x\mapsto x^p$. To handle  this fact, we introduce 
	{\bf $p$-th roots of coordinate functions}, called $\lambda$-functions. For instance, for $x\in \mathbb{F}_p((X))$ we write (using the fact that $1,X,\ldots, X^{p-1}$ is basis of 
	$\mathbb{F}_p((X))$ over the subfield  consisting of $p$-th powers)
	$$ x=x_0^p + x_1^pX+ \dots + x_{p-1}^pX^{p-1}$$
	and define $$\lambda_i(x):=x_i \quad  (i=0,\ldots, p-1).$$
	
	Notice that $\lambda$-functions are both universally and existentially definable in $L$. We prove then $B_1(N)\cap Q(N)$  is positive quantifier-free definable in the language $L$ together with $\lambda$-functions (see Theorem \ref{boule_up}).

	This article is organized as follows:
	generalities and definitions about the non-commutative ring $R$, its matrix ring and $R$-modules are given 
	Section 2. In the third section we prove an  equivalent of Hensel's lemma (see Theorem \ref{hensel-onay}) -we don't know if it is really new-
	and use it to axiomatize what we call {\it henselian filtered $R$-modules}. This part  consists of finding a ball (like $\mathfrak{m}$ above) where the trace of the definable set considered can be defined by a positive quantifier-free formula in the language $L$ together with the $\lambda$-functions. In the fourth section, we introduce the notion of a valued module, and study   pseudo-complements in order to satisfy the task (\ref{task1}) above. In the final section we apply the previous results to obtain  the main theorems. Here we use some general model theory of modules (e.g. Baur-Monk elimination) and some recent facts from \cite{koeningsman-sylvie} and \cite{sylvyarno17}
	
	In the course of this article we introduce several languages and theories.  For the sections 2-4, (i.e. except the last section), it is worth noting that $K$ and $\mathbb{F}_p((X))$, realized with the natural interpretation of these languages are  models of all of these theories. When we say that something is {\it computable} we mean that there is a Turing machine which can compute it.
	
	\subsection*{Acknowledgements}
	This is a long-standing work, which began with my Ph.D thesis in early 2007. I would like to thank  my advisor Fran\c{c}oise Point for suggesting me the topic and  for giving the idea to introduce $\lambda$-functions. I thank  Fran\c{c}oise Delon, also my advisor; without her help, this work could not be finished.
	I want to thank  Franz-Viktor Kuhlmann for having invited me, as early as  2008 and for sharing his ideas about the topic. Lou van den dries and Matthias Aschenbrenner have kindly discussed with me  some aspects of the present article in last times. Finally, I want to thank  Sylvy Anscombe and Arno Fehm for their encouragement to write this manuscript.

	\section{Preliminaries}

	Let $p$ be a prime number and $k$ a positive integer, we set $d:=p^k$. Let $K$ be the henselization of  the field of rational functions over $\mathbb{F}_d$, with respect to the $X$-adic valuation; denoted as $$K=\mathbb{F}_d(X)^{h}.$$ Let
	$$\vfi:x\mapsto x^{d} $$ be the $k$-th power of the Frobenius  endomorphism of $K$. Recall that $K$ is a finite extension of $K^{\vfi}$ of dimension $d$. \footnote{Note that the contents of this section can be generalized to any field  and  to any endomorphism
		satisfying  similar hypotheses. In particular one can 
		consider $\mathbb{F}_q(X), \mathbb{F}_q((X)), \ldots etc$. As it is our main interest and for readability reasons,  we preferred to stick to one case.} We fix the notation $v_K$ for the $X$-adic valuation on $K$.

	Set $\alpha(0):= \{1\}$ and let $$\alpha:=\alpha(1)=(\alpha_0,\alpha_1,\ldots, \alpha_{d-1})$$ be  a basis of  the $K^{\vfi}$-vector space of $K$. One can think of the basis $(1,X, \ldots, X^{d-1})$.    It is easy to see that for all $n\geq 1$, $\alpha$  induces canonically
	the (ordered -) basis $\alpha(n)$ of the $K^{\varphi^{n}}$-vector space $K$. We identify $d^n$ with the set of functions $\{0,\dots,n-1\} \to \{0,\dots,d-1\}$ ordered lexicographically. We denote by $\alpha_k$ the $k$-th element of the basis $\alpha(n)$.
	

	We define the ring \begin{equation}
		R:=K\langle t  \mid ta^\vfi=at \rangle .  
	\end{equation}
	the  (right-) $K$-algebra with the indeterminate $t$, subject to the commutation rule
	$ta^\vfi=at$ for all $a \in K$. When 
	$k=1$, $R$ is isomorphic to the ring of additive polynomials over $K$, equipped with addition and composition. Every non zero  $q\in R$ can be  written as
	\begin{equation}\label{R-polynomial}
		q=t^na_n + t^{n-1}a_{n-1} + \dots + a_0 
	\end{equation}
	where the $a_i$ are from $K$ and $a_n\neq 0$. Since the image of $q$ in $K[T]$ is
	$$Q(T)=a_nT^{d^n}+\dots+a_0T$$ under the aforementioned isomorphism, ``the constant term $a_0$ of $q$" gives rise to the linear term $a_0T$ of $Q(T)$.
	\begin{defn}
		An element  $q$ as in (\ref{R-polynomial}) is said to be separable if $a_0\neq 0$.
	\end{defn}
	\begin{remark}
		Note that $q$ is separable if and only if its image $Q(T) \in K[T]$ is separable.
	\end{remark}
	The integer $n$ in the expression (\ref{R-polynomial}) is called the degree of $q$ and we set the degree of $0$ as $-1$. 
	With  the degree function $R$ is right euclidean: for every non zero $r,q \in R$ there exists
	$q'$ and a unique $r'$ of degree $<\deg(q)$ such that $r=qq' +r'$. As a consequence, 
	least common multiple $\lcm(q,r)$ and greatest common divisor $\gcd(q,r)$ are well defined.
	In particular $R$ is right Ore. Note that if $\vfi$ is onto then $R$ is also left euclidean and 
	if it is not (which is the case that we consider here), $R$ is not even left Ore. The reader can see \cite{cohn} or \cite{valmod1} Section 2, for more details.

	For our quantifier elimination Theorem \ref{boule_up} we will import some definitions and results from \cite{ddp1} on the ring $R$ and 
	on $R$-modules. We give the proofs of  some results not to be self-contained but rather, to initiate the reader with the style of the computations that will permit us to refer   easily  to \cite{ddp1}.
	\begin{fait}
		The field $K$ is recursively enumerable since $\mathbb{F}_d(X)$ is, and $K$ is formed by adding to $\mathbb{F}_d(X)$ the unique root  of the  each polynomial satifying the hypothesis of the Hensel's Lemma . Consequently the ring $R$ and its matrix rings are recursively enumerable.
	\end{fait}
	
	\begin{lemma}\label{lem:decompq}
		Let $q\in R$ and $n$ be a positive integer. Let $\alpha$ be a basis of the
		$K^\vfi$-vector space $K$. Then 
		\begin{enumerate}[1.]
			\item  \label{qi} $q$ can be uniquely written as:
			\begin{equation}
				\sum_{i \in d^n} q_i\alpha_i \quad (q_i \in R),
			\end{equation}
			
			\item for all $n>0$ there is an endomorphism $\nqrt{n}{\cdot}:(R,+) \to (R,+)$ such that
			if  $q=\sum_{i \in d^n} q_i\alpha_i$  then 
			\begin{equation}
				t^nq = \sum_{i \in d^n} \nqrt{n}{q_i}t^n\alpha_i.
			\end{equation}
			In addition, if $q$ is separable there exists $i \in d^n$ such that
			$\nqrt{n}{q_i}$ is separable.
		\end{enumerate}\medskip
	\end{lemma} 
	\begin{proof}
		\item[1.]  Let $0 \leq k\leq s=\deg{q}$. For each monomial $t^ka_k$, by expressing $a_k$ with respect to the basis $\alpha(n)$,  $q$ can be written as
		$$t^k(\sum_{i\in d^n} a_{k,i}^{\varphi^n}\alpha_i) .$$ Set 
		$q_{k,i}=t^ka_{k,i}^{\varphi^n}$ and  $q_i = \sum_{0\leq k \leq s} q_{k,i}.$ 
		So we get 
		$$q=\sum_{i \in d^n}q_i\alpha_i.$$
		
		\item[2.] Let $a \in K$. Write 
		$$a=\sum_{i \in d^n} a_i^{\vfi^n}\alpha_i .$$ We define
		$$\nqrt{n}{a}:=\sum_{i \in d^n} a_i\alpha_i$$ and extend it to $r \in R$ by applying it to its coefficients. By definition $\nqrt{n}{\cdot}$ preserves the addition.
		
		To see that \begin{equation}\label{nqrt} 
			t^nq = \sum_{i \in d^n} \nqrt{n}{q_i}t^n\alpha_i
		\end{equation} it is enough to show that
		$$t^nq_i=\nqrt{n}{q_i}t^n .$$ 
		Note that
		$\nqrt{n}{q_{k,i}}=t^ka_{k,i}$ and hence  
		$$\nqrt{n}{q_{k,i}}t^n=t^{n+k}a_{k,i}^{\vfi^{n}}=t^{n}q_{k,i} .$$ 
		Using additivity,
		we get $\nqrt{n}{q_i}t^n=t^nq_i$.
		
		Suppose now $q$ is separable. If none of the  $\nqrt{n}{q_i}$ is separable, then  for all  $i$, we have  $\nqrt{n}{q_i}=tq_i'$. But then 
		$t^nq \in t^{n+1}R$ by the equality (\ref{nqrt}). Hence
		$q \in tR$,  which is a contradiction.
	\end{proof}
	

	\begin{defn}\label{matriss} An $m\times n$  matrix  $A=(q_{i,j})$ over $R$ is said to be 
		\begin{enumerate}[1.]
			\item  lower triangular if, 
			$j>i$ implies $q_{i,j}=0$,
			\item   lower triangular diagonally  separable 
			if it is lower triangular, $n\leq m$ and the  $q_{ii}$ ($i\leq n$) are separable,
			\item lower triangular separable  if $A=(A_1,0)$ where $A_1$ is an $m\times k$ lower triangular diagonally separable matrix, 
			and $0$ is the $m\times l$ null matrix with $k+l=n$.
			
		\end{enumerate}
	\end{defn}

	\begin{prop}\label{mattri}
		For any matrix $A$ there exists an invertible matrix  $P$  with coefficients 
		in  $\{0,1\}$ and an invertible matrix
		$Q$ such that $PAQ$ is lower triangular.
		\begin{proof}
			See the Proposition $6.1$ in  \cite{ddp1}.
		\end{proof}
	\end{prop}
	
	\subsection{$R$-modules}
	For the rest of the article we will always understand the expression $R$-module as right $R$-module. In an $R$-module $M$, scalar multiplication will be denoted as $x.r$, for $x\in M$ and $r\in R$.
	\begin{defn}
		Let $M$ be an  $R$-module and $\beta$ a basis of the $K^\vfi$-vector-space $K$. 
		\begin{enumerate}[1.]
			\item $M$ is said to be $t\beta$-decomposable if $x\mapsto x.t$ is injective and we have
			\begin{equation}
				M=\bigoplus_{i \in d}M.t\beta_i.
			\end{equation}
			(where $\oplus$ indicates the direct sum as abelian subgroups). 
			\item $M$ is said to be $t$-decomposable if it is $t\alpha$-decomposable with
			$\alpha=(1,X,\dots, X^{d-1})$.  Furthermore, for $i\in d$ we then define 
			$\lambda_i(x)=x_i$ where 
			$$x=x_0.t\alpha_0 + \dots + x_i.t\alpha_i + \dots x_{d-1}.t\alpha_{d-1}.$$
		\end{enumerate}\medskip
	\end{defn}
	
	\begin{remark}\label{lambdascompose}
		The direct sum 
		$$M=\bigoplus_{i \in d}M.t\beta_i$$ induces the direct sum below for every positive $s$
		$$M=\bigoplus_{i \in d^s} M.t^s\beta_i.$$
	\end{remark}
	\begin{nota}
		For the rest of the article we set  
		$\alpha:=\{1,X,\dots,X^{d-1}\}$.
	\end{nota}

	\begin{remark}\label{extunideflambdas}
		The functions $\lambda_i$ are both existentially and universally definable in the language of right $R$-modules:\begin{equation}
			y=\lambda_i(x) \iffl \forall x_1, \dots,\forall x_i, \dots, \forall x_d \quad x=\sum_{j=1}^d x_i.t.\alpha_j \longrightarrow x_i=y
		\end{equation}
		and
		\begin{multline}
			y=\lambda_i(x) \iffl \exists x_1, \dots, \exists x_{i-1}, \exists x_{i+1}, \dots, \exists x_d \quad \\ 
			x=\sum_{j=1}^{i-1} x_j.t\alpha_j + y.t\alpha_i + \sum_{j=i+1}^d x_j.t\alpha_j.
		\end{multline}
		
		For a positive $s$, by Remark \ref{lambdascompose} above, we get canonically the $\lambda$ functions of level $s$,
		defined for all $i\in d^s$, in an obvious way.
	\end{remark}

	{\bf The language $L(\lambda)$.} For the rest of the article we let $L(\lambda)$ be the language of $R$-modules together with  the functions (-symbols) $\lambda_i$. 	
	
	\begin{defn} We denote by $T_\lambda$ the $L(\lambda)$-theory of $t$-decomposable $R$-modules, that is, the theory of $R$-modules together with the axioms
		
		{\bf $\lambda$-decomposition :}\label{ax:lambdas} 
		$$\begin{matrix}  \forall x &  x = \sum \lambda_i (x) . t \alpha_i \\
		\forall  x \forall (x_i )_{(i\in d)} & (x = \sum_i x_i. t \alpha_i \to  \bigwedge_i x_i = \lambda_i (x)). \end{matrix}$$
	\end{defn}
	
	\begin{lemma}\label{Lambdaterms} In any $t$-decomposable $R$-module $M$, any $L(\lambda)$-term can be evaluated on the tuple $(x_i)_i$ from $M$,   as 
		$$\sum_i \sum_j \lambda_j(x_i).r_{ij}$$ where 
		$r_{ij} \in R$.
	\end{lemma}
	\begin{proof}
		This is  Corollary 3.3 in \cite{ddp1}.
	\end{proof}
	
	\begin{lemma}\label{lem:passagelambdasterms}
		Let  $m > 0$, and $q_ j , q_ j' \in R$ such that  $q_j = t^m .q_j'$ . Then the equation
		$\sum y_j.q_j = u$ is equivalent to
		$$ \bigwedge_{i \in d^{m}}\sum_j y_i.\nqrt{n}{q'_{j_i}}=\lambda_i(u).$$ 
		in any $t$-decomposable $R$-module.
	\end{lemma}
	\begin{proof}
		This is  Lemma 3.4 in \cite{ddp1}
	\end{proof}
	
	\begin{obs}
		Let $q=(q_0, \dots, q_{n-1})$ be a non zero tuple from $R$. We set
		$$e:=\min\{k \tq q_i \in t^kR, \mesp \text{for all}\, i\}.$$
		Notice that  $e=0$ means that at least one of the $q_i$ is separable. 
		
		Suppose $e>0$, one can write 
		$$q_i=t^eq'_i=\sum_{k \in d^e} \nqrt{e}{q'_i}_{k}t^e\alpha_k .$$
		Since at least one of the $q'_i$ is separable, for some $(i,k)$, 
		the polynomial  $\nqrt{e}{{q'_{i}}_{k}}$ is separable by Lemma \ref{lem:decompq} (2.).
		Set $q_{i,k}= \nqrt{e}{q'_i}_k $  and  let $q^e_\lambda$ be the  $n\times d^e$-matrix $(q_{i,k})$ whose $i$-th line consists of 
		the sequence $(q_{i,k})_{k \in d^e}$.
		By  this process we have replaced the tuple $q$ by a matrix which  has at least one separable coefficient. Iterating this process and using  Lemma \ref{lem:passagelambdasterms} above we get the following result. 
	\end{obs}
	
	\begin{lemma}\label{ddp6.4}
		Let  $A$ be a non zero $n\times k$ matrix over $R$. Then, 
		the system  $y.A = u$ is equivalent to 
		$$y.PQ = w(u)$$ modulo $T_{\lambda}$,
		where $P$ is a permutation matrix (i.e. invertible with coefficient in $\{0,1\}$),  $Q$ is lower triangular separable  and $w$ is a tuple 
		consisting of $L({\lambda})$-terms. 
	\end{lemma}
	\begin{proof}
		This is a reformulation of  Lemma  6.4.4 in \cite{ddp1}.
	\end{proof}
	
	\subsection{Baur-Monk Elimination}\label{sec:baur-monk}
	
	The following is a reminder of Theorem A.1.1, Corollary A.1.2 and the discussion which follows in \cite{hodges}, p. 653-656.
	
	Let $\mathcal{L}$ be any language which contains the language $\{+,-,0\}$ of abelian groups.  A {\it positive primitive formula} ($p.p.$) $\phi$ of 
	$\mathcal{L}$, is the one of the form
	$$\exists \bar{y} \mesp \left(\bigwedge_i \psi_i(\bar{x},\bar{y})\right)$$
	where the $\psi_i$ are atomic.

	A {\it group-like $\mathcal{L}$-structure} $A$, is an $\mathcal{L}$-structure whose base set is a group
	with respect to $\{+,-,0\}$. A {\it basic formula} for an $\mathcal{L}$-theory $T$, is a $p.p.$ formula which defines  a subgroup of the corresponding cartesian power of any $\mathcal{L}$-structure $A\models T$.

	Note that if $S$ is any ring and $\mathcal{L}$ is the language of $S$-modules then in any $S$-module $N$ any $p.p.$ formula defines a subgroup
	of the corresponding cartesian power of $N$, hence any $p.p.$ formula is basic for the theory of $S$-modules.
	
	Let $T$ be an $\mathcal{L}$-theory such that every model of $T$ is group-like, 
	and every $p.p.$ formula is a basic formula for $T$. An invariant 
	$\mathcal{L}$-sentence  of $N\models T$ is  an $\mathcal{L}$-sentence 
	$\Theta$ satisfied by $N$, 
	such that for some   $p.p.$ formulas $G(x)$ and $H(x)$ of one variable $x$,   for some $m\in \mathbb{N}$,
	$$T\models \Theta \Longleftrightarrow  \vert G/G\wedge H \vert =m$$
	where the right-hand side of the equivalence is an abbreviation of the formula:
	\begin{multline}
		\exists  (x_i)_{(i=1..m)}  (\bigwedge_i G(x_i) \wedge \\
		\bigwedge_{i\neq j} \neg (G\wedge H)(x_i-x_j)\wedge \forall z \mesp (G(z) \rightarrow \bigvee_i (G\wedge H)(z-x_i))).
	\end{multline}
	An invariant sentence is an invariant sentence of some $N$.
	\begin{theorem}[Baur-Monk]\label{thm:baur-monk}
		Every $\mathcal{L}$-formula is equivalent modulo $T$ to a boolean combination of
		$p.p.$ formulas and invariant sentences. Hence for all models $N,M \models T$, 
		$N \equiv M$ if and only if $N$ and $M$ have same invariant sentences. 
	\end{theorem}
	
	\begin{corr}\label{corr:baur-monk}
		A completion of $T$ is model-complete if and only if every $p.p.$ formula is equivalent to a universal formula modulo $T$.
	\end{corr}

	\section{The tropical action of $R$ on $\mathbb{Z}$ and Filtration}

	We recall some elementary facts about henselian valued fields. Let $(F,v)$ be a valued field with value group $\mathfrak{G}$ and valuation ring
	$\mathcal{O}$. We set $\Gamma=\mathfrak{G}\cup\{\infty\}$,    extend the usual addition of $\mathfrak{G}$ to $\Gamma$ by letting 
	$$
	\infty+\infty=\infty + a=a+\infty=\infty
	$$ for all $a \in \mathfrak{G}$.
	We recall 
	the {\it tropicalisation} of a one variable polynomial $Q(T)$ over $F$:
	Write $$Q(T)=\sum a_iT^i.$$ Then tropicalisation of $Q$ is the map
	$$Q_v: \Gamma \to \Gamma$$
	$$\gamma \mapsto  \min_i\{i\gamma + v(a_i)\}.$$

	A {\it jump value} (or a tropical zero) of $Q$ is some $\gamma \in \Gamma$ such that
	$$ \vert \{i \tq  i\gamma + v(a_i)= Q_v(\gamma)\} \vert \geq 2 .$$
	We denote by $\jump{Q}$ the set of jump values of $Q$. Note that this set is finite and has at most $n-1$ element if $Q$ is of degree $n$.

	\begin{fait}[Newton's Lemma] Let $(F,v)$ be a valued field and $f$ be a polynomial with coefficients from the valuation ring  
		$\mathcal{O}$. Consider the following property $h(f)$ of $f$,
		\begin{multline}
			v(f(0))>2v(f'(0)) \Rightarrow (\exists b \mesp 
			f(b)=0 \quad \text{and} \quad v(b)=vf(0)-vf'(0)) \tag{{\it h(f)}}
		\end{multline}
		Then $(F,v)$ is henselian if and only if $h(f)$ holds for every $f$ over the valuation ring $\mathcal{O}$ of $F$.
	\end{fait}
	\begin{proof}
		See \cite{prestel} Theorem 4.1.3.
	\end{proof}
	
	Consider a polynomial $G$,  such that 
	$G(0)=0$ and $G'(0)\neq 0$. Then $G$ is of the form 
	$$G(T)=aT + \text{sum of monomials of higher degree}.$$
	Consider the set $A_1:=\{ \gamma \tq  G_v(\gamma) <   (G-aT)_v(\gamma) \}$. This is a non-empty final segment of $\Gamma$. Also let 
	$$A_2:= \{\gamma + v(a) \tq \gamma \in A_1\}.$$ 
	We set 
	$B_1(G):=v^{-1}(A_1)$ and $B_2(G):=v^{-1}(A_2)$. 
	
	Note that $B_i$  $(i=1,2)$ are convex subsets of $F$, that is, inverse images by $v$ of convex subsets of $vF$.
	
	\begin{theorem}\label{hensel-onay} $(F,v)$ is henselian if and only if  each polynomial $G$ with coefficients in the valuation ring  $\mathcal{O}$, such that
		$G(0)=0$ and $G'(0)\neq 0$,  induces a bijection $B_1(G) \to B_2(G)$.
	\end{theorem}
	\begin{proof}
		Suppose $(F,v)$ is henselian. By the definitions of $B_1$ and $B_2$ we have 
		$G(B_1)\subseteq B_2$. We will show the converse inclusion. Let $z \in B_2$ and set $f=G-z$. We have 
		$vf(0)=v(z)=\gamma + v(a)$ where $\gamma \in A_1$. Hence 
		$ G_v(\gamma) =\gamma + v(a) <  (G-aT)_v(\gamma)$. Since a monomial of $G-aT$ is at least of degree $2$ and $G$ is over $\mathcal{O}$, we have 
		$$\gamma + v(a)<2\gamma.$$
		It follows that $\gamma > v(a)$ and $v(z)>2v(a)$. By the above fact, there is root $b$ of 
		$f$, of valuation $\gamma$. In other words, $G(b)=z$ and hence $G(B_1) \supseteq B_2$. 
		
		Now we show that $G \restriction B_1$ is 1-1. Let $z$ and $f$ be as above and $G(x_1)=G(x_2)=z$ with $x_1,x_2 \in B_1$. Then   $x_1,x_2$ have the same  valuation $\gamma$, but it is well-known that $f$ has at most one root of valuation $\gamma$ (or it can be proven by using Taylor expansions that $v(G(x_1-x_2))>G_v(v(x_1-x_2)$ but this is not possible since $v(x_1-x_2)$ is bigger than then all the jump values of $G$).
		
		The converse is similar and easy (and we won't need).   
	\end{proof}

	\begin{corr}
		Let $G$ be any polynomial over $(F,v)$ such that $G(0)=0$ and $G'(0)\neq 0$ then there exist convex subsets $B_1$, $B_2$ such that
		$G\restriction B_1:B_1 \to B_2$ is a bijection.
	\end{corr}
	\begin{proof}
		Divide $G$ by the  coefficient which has the minimal valuation among the coefficients of $G$. Let $H$ be the obtained polynomial (which has coefficients over $\mathcal{O}$). Then $B_1:=B_1(H)$ and $B_2:=aB_2(H)$ suits. 
	\end{proof}
	
	\begin{nota}
		Let $G$ be as above. If the value group of $F$ is discrete of rank 1, then $A_1,A_2$ are closed intervals; in this case, we set $$\h(G):=\min A_1 \quad \text{and} \quad \hens(G):=\min{A_2}.$$
	\end{nota}
	{\bf The tropical action of $R$.} From now on we set 
	$$\Gamma:=\mathbb{Z}\cup \{\infty\}$$ and we equip $\Gamma$ with a right action of $R$ using tropicalisations:
	
	Let 
	$q=\sum_i t^ia_i \in R$.  We define the tropicalisation 
	of $q$, as the tropicalisation of the $Q(T)=\sum_i a_iT^{d^{i}}$, that  is,  as  the map 
	\begin{equation}
		\cdot q:\Gamma \to \Gamma; \quad \gamma \mapsto \gamma\cdot q= \min_i\{d^i\gamma + v(a_i)\}.
	\end{equation}
	In particular 
	$\cdot q$ is strictly increasing  if  $q\neq 0$.
	Note that $\gamma \cdot a =  (aT)_v(\gamma) =\gamma + v(a)$ for $a \in K$.

	We set  $\jump{q}:=\jump{Q}$, $\h(q):=\h(Q)$ and $\hens(q):=\hens(Q)$. For instance,
	if $q=t-1$ then $A_1=A_2:=\{\gamma \in \Gamma \tq \gamma>0\}$. Hence
	$\h(q)=\hens(q)=1$.

	For our interests, we also introduce the tropicalisations of the $\lambda$-functions. Remark that we have for all $x\in K$,
	$x=\sum_{i \in d} \lambda_i(x)^dX^i$, 
	$$v_K(x)=\min_i\{v_K(\lambda_i(x)^d)+i\},$$
	and the minimum is attained for a unique $i\in d$.
	We define 
	$$ \lambda_i(\gamma):=\begin{cases} \frac{\gamma - i}{d}\quad \text{if} \mesp \gamma \in d\mathbb{Z} + i \\
	0 \quad \text{else}
	\end{cases}$$
	and we set 
	$$\lambda(\gamma):=\sum_i \lambda_i(\gamma).$$

	\begin{remark}\label{rem:lambdaregular} For $x\in K$ 
		$$v_K(x)\geq \gamma \Leftrightarrow \quad \text{for all}\mesp i\in d \quad v_K(\lambda_i(x)) \geq \lambda(\gamma).$$
	\end{remark}
	\begin{proof} Note that $\Leftarrow$ is clear and if $\gamma=v_K(x)$ then the assertion is trivial. Suppose $v_K(x)>\gamma$. Let $j$ be such that 
		$v_K(\lambda_j(x).t^{d}\alpha_j)=dv_K(\lambda_j(x))+j=v_K(x)$. Let $i$ be such that 
		$\gamma=d\lambda(\gamma)+i$.  It follows that
		$$d(v_K(\lambda_j(x))-\lambda(\gamma))>j-i.$$
		Hence $v_K(\lambda_j(x))-\lambda(\gamma)>-1$ and
		$v_K(\lambda_j(x))-\lambda(\gamma)\geq 0$.
	\end{proof}
	\begin{corr}\label{corr:lambas(s)}For all non zero $x\in K$ and $\gamma \in \Gamma$, 
		$$v_K(x)\geq \gamma \Leftrightarrow \quad \text{for all}\mesp i\in d^s \quad v_K(\lambda_i(x)) \geq \lambda^s(\gamma)$$
		where $\lambda^s=\underbrace{\lambda\circ\ldots\circ\lambda}_{s-times}.$
	\end{corr}

	\begin{remark}Let $\gamma\in \Gamma$ and $r,q \in R$. Then
		\begin{enumerate}[1.]
			\item $\gamma\cdot rq = (\gamma \cdot r)\cdot q$,
			\item $\gamma\cdot (r + q)\geq\min\{\gamma\cdot r, \gamma \cdot q\}$,
			\item $\cdot r$ is strictly increasing for all non zero $r$, 
			\item $\infty \cdot r =  \gamma \cdot 0 = \infty$ for all $r$ and $\gamma$.
		\end{enumerate}
	\end{remark}
	\begin{proof}
		This follows by direct computations from the definition. Note that (1) follows more generally from  \cite{gonenc}, Corollary 4.1.14.
	\end{proof}
	
	\begin{fait}
		The theory of $\Gamma$ together with the tropcial action of $R$ and
		$\lambda$ is decidable since this structure  is definable in the ordered abelian group structure of $\mathbb{Z}$, together with $\infty$ and constants for the elements of $\mathbb{Z}$.
	\end{fait}
	
	\begin{defn}	
		We call this structure the tropical structure of $\Gamma$.	
	\end{defn}
	\subsection{Quantifier Elimination Near 0}

	\subsection*{\bf The languages $L_{\mathcal{O}}$ and $L_{\mathcal{O}}(\lambda)$} The language  
	$L_{\mathcal{O}}$ is obtained by adding to $L$ a unary predicate $\mathcal{O}$, and the language $L_\mathcal{O}(\lambda)$ is the language $L_\mathcal{O}$ together with the functions (-symbols) $\lambda_i$ added to $L_\mathcal{O}$. We want to study the divisibility conditions for $t$-decomposable $R$-modules which can be seen
	as  properties reflecting a kind of henselianity, analog to Theorem \ref{hensel-onay} just above.
	
	\begin{nota}
		For the rest of this article, for $\gamma \in \Gamma$,
		we write $P_\gamma$ 
		for the predicate  $\mathcal{O}.X^{\gamma}$, i.e. in any $L_\mathcal{O}$-structure
		$M$, $x\in P_\gamma$ if and only if $x.X^{-\gamma} \in \mathcal{O}$.
	\end{nota}
	\begin{defn}
		An henselian filtered module, is a $t$-decomposable $R$-module which is an $L_\mathcal{O}(\lambda)$-structure satisfying the following axioms:
		
		\item[0. {\bf Balls:}] $P_{\infty}=\{0\}$  and the $P_{\gamma}$ form a chain of subgroups decreasing with $ \gamma$  such that  the  inclusions are proper.
		\item[1. {\bf Ultrametric:}]  $\forall x  \forall y  \mesp x \in P_{\gamma}\wedge y \in P_{\delta}  \rightarrow (x.r + y) \in P_{\min\{\gamma\cdot r,\delta\}}$   for all $\gamma,\delta \in \Gamma$ and $r\in R$.
		\item[2. {\bf Regularity:}]   $\forall x \mesp  x \in P_{\gamma} \leftrightarrow x.r \in P_{\gamma \cdot r}$, for all $\gamma$ and $r\neq 0$ such that
		$\gamma \notin \jump{r}$.
		\item [3. $\lambda$-{\bf regularity:}]  $\forall x \, x\in P_{\gamma} \leftrightarrow \bigwedge_i \lambda_i(x) \in P_{\lambda(\gamma)}$, for all $\gamma$. 
		\item [4. {\bf Henselianity:}] $\forall x  \in P_{\hens(s)}\setminus \{0\} \exists ! \, y \in P_{\h(s)} \mesp y.s=x$ for all separable $s$.
		
	\end{defn}
	Note that axiom 3 implies that
	$$\forall x \, x\in P_{\gamma} \leftrightarrow \bigwedge_{i\in d^s} \lambda_i(x) \in P_{\lambda^s(\gamma)},$$ a consequence analog to the one expressed in Corollary
	\ref{corr:lambas(s)}.   
	\begin{nota}
		We denote by $T_{hens}$ the theory of henselian filtered modules.
		
		We also isolate some theories of the $R$-modules that are already considered in \cite{ddp1}:
		\begin{itemize}
			\item we denote by $T_{free}$ the theory $T_\lambda$ together with the following (scheme of-) axioms:
			\label{ax:separablyclosed} $$\forall x \exists y \mesp x=y.s,$$
			for all separable $s$,
			\item  we denote by $T^{0}_{free}$, the $L_{\lambda}$-theory of torsion-free non zero models of $T_{free}$.
		\end{itemize}
	\end{nota}
	\begin{remark}\label{P_M} In a henselian filtered module $M$, we denote by 
		$P_M$ the subgroup defined by  the intersection 
		$$P_M:=\bigcap_{\gamma\neq \infty} P_\gamma(M).$$
		Note that by  ultrametric and regularity axioms, $P_M$ is a 
		$L(\lambda)$-substructure of $M$ and it is torsion-free as an $R$-module. Moreover, it is straightforward to check that by Hensel's axioms 
		$$P_M \models T^0_{free}.$$
	\end{remark}
	\begin{proposition}\label{prop:t0freedecidable}  Given $L$-$p.p.$ formulas	$\mathfrak{a}(x)$ and $\mathfrak{b}(x)$  of one variable $x$, 
		the quotient $\mathfrak{a}/(\mathfrak{a}\wedge \mathfrak{b})$ is either trivial or infinite in every model of $T^0_{free}$. In particular $T^0_{free}$ is complete and decidable. 
		\begin{proof}
			This follows from Lemma 6.8 in \cite{ddp1}.
		\end{proof}
		
	\end{proposition}
	
	\begin{defn}
		A ball of $L_\mathcal{O}$, is an atomic formula $W(x_1,\ldots,x_k)$ of the form
		$$W(x_1,\ldots, x_k): \bigwedge_{i=1}^k x_i \in P_{\gamma_i}.$$ We will write it rather as a product 
		$$W=P_{\gamma_1}\times \dots \times P_{\gamma_k}$$ of predicates. 
		$W$ is said to be proper if none of the ${\gamma_i}$ is equal to $\infty$. 
	\end{defn}
	
	\begin{remark}
		A  positive primitive  formula $\phi(\bar{x})$ of $L_{\mathcal{O}}$ is equivalent to one in the form
		$$\exists \bar{y} \mesp\bar{x}.B-\bar{y}.A\in \mathcal{O}^k\times\{0\}^n$$ which it self is equivalent to a formula 
		$$\exists \bar{y} \mesp\bar{x}.B'-\bar{y}.A'\in W$$ 
		where $A,B,A',B'$ are matrices over $R$ and $W$ is a ball. Hence any $p.p.$ formula is a basic formula. For our purposes will rather use the latter
		equivalence.
	\end{remark}
	
	Note that modulo $T_{hens}$, the set of  definable sets  by an $L_{\mathcal{O}}$-$p.p.$ formula contains the set
	of definable sets by an $L$-$p.p.$ formula 
	since $W$ can be
	chosen equal to $P_{\infty}^k$ for some $k$.  
	
	In this section, we will prove the following:
	\begin{theorem}\label{boule_up}
		Let $\phi(x_1,\dots,x_m)$ be a $p.p.$ formula of
		$L_\mathcal{O}$. Then there is some computable $\delta\neq \infty$ and some positive quantifier free
		$L(\lambda)$-formula $\psi$ which depends only to the theory $T_\lambda$, such that, with
		$V=P_\delta^m$,  we have
		$$T_{hens} \models \phi \wedge V \leftrightarrow \psi \wedge V .$$ 
		Moreover,
		If $\phi_1$ is another $p.p.$ formula with the same arity, and
		if $\phi$ and $\phi_1$ are equivalent modulo $T^0_{free}$, then for some computable $\gamma$ and with $W:=P_\gamma^m$ 
		$$ T_{hens} \models \phi\wedge W \leftrightarrow \phi_1 \wedge W .$$
	\end{theorem}

	\begin{lemma}\label{croissant} Given $\delta\neq \infty$ and an $L(\lambda)$-term $u(\bar{x})$, there exists $\gamma\neq \infty$ such that 
		\begin{equation}
			T_{hens} \models \forall \bar{x} \mesp \bar{x} \in P_\gamma^{\vert \bar{x} \vert} \rightarrow u(\bar{x}) \in P_{\delta}.
		\end{equation}
	\end{lemma}
	\begin{proof} Put $u(\bar{x})$ in the form $\sum_{i,j} \lambda_j(x_i).r_{ij}$ using Lemma \ref{Lambdaterms}. 
		By regularity, $\lambda$-regularity and ultrametric axioms, for any $\rho$,
		$$\bar{x} \in P_\rho^{|\bar{x}|} \rightarrow u(\bar{x}) \in P_{\min\{\lambda(\rho)\cdot r_{ij}\}}$$ holds. Since all the tropical functions $\lambda$ and $\cdot r_{ij}$ are unbounded and increasing,
		one can choose $\gamma$ such that $\gamma \cdot r_{ij} \geq \delta$ for all $i$.
	\end{proof}
	\begin{remark}\label{rem:gamma_computable}
		The value $\gamma$  is computable from $\delta$ and $u(\bar{x})$ in the tropical structure of $\Gamma$ .
	\end{remark}

	As a consequence of the henselianity axioms we  observe the following.
	\begin{lemma}
		Given  $\gamma\neq \infty$ and a separable $s$ there exists  $\delta\neq \infty$ such that
		
		$$ T_{hens} \models \forall y \mesp  [(y \in P_\delta) \rightarrow (\exists x  \in P_\gamma \wedge x.s=y)].$$
		
	\end{lemma}
	\begin{proof}
		If $ \gamma \leq \h(s)$ then set $\delta=\hens(s)$. Otherwise $\gamma \notin \jump{s}$ and by the 
		regularity axioms $x \in P_{\gamma} $ if and only if $x.s \in P_{\gamma \cdot s}$. 
		Since $P_\gamma \subseteq P_{\h(s)}$, $P_{\gamma \cdot s} \subseteq  P_{\hens(s)}$. 
		Hence if $y \in P_{\gamma \cdot s}$ then the unique solution $x$ such that $x.s=y$ lies in 
		$P_\gamma$.  So $\delta=\gamma \cdot s$ fits for our requirements.
	\end{proof}
	
	\begin{remark}\label{rem:delta_computable}
		As above, $\delta$ is computable in the tropical structure $\Gamma$. 
	\end{remark}

	\begin{corr}\label{bouledivisible}
		Let $A=(a_{ij})$  be an $m\times k$ lower triangular diagonally separable matrix (in particular  $k \leq m$). Then,  
		for every proper ball $W$, there exists a proper computable ball $W_1 $ such that  
		$$T_{hens} \models \forall \bar{x} \mesp (\bar{x}\in W_1 \rightarrow \bar{x}\in  W.A) .$$ 
	\end{corr} 
	\begin{proof}
		Write $W = \prod_{i=1}^{m} P_{\delta_i}$. By Remark \ref{croissant},  for fixed $i$, chose  $\gamma_i  \in \Gamma\setminus\{\infty\}$ such that,
		$$\sum_{j \neq i} P_{\gamma_i}.a_{ji} \subset P_{\delta_i}.$$
		Since the $a_{ii}$ are separable, by the above lemma there exist proper balls $U_i$ such that whenever   $z_i \in U_i$,  
		there exists $y_i \in P_{\gamma_i}$, such that $z_i = y_i.a_{ii}$; hence 
		$$z_i - \sum_{j=1}^{n} y_i.a_{ji}= \sum_{j\neq i} y_i.a_{ji} \in P_{\delta_i}$$ for all $1\leq j \leq k$. Take $W_1=
		\prod_{i} U_i.$ 
	\end{proof}
	
	{\bf Proof Theorem \ref{boule_up}.}
	Let $\bar{y}:=(y_1,\dots,y_k)$, $\bar{x}:=(x_1,\dots,x_m)$,  $\phi(\bar{x})$,
	$W:=\prod P_{\gamma_i}$ and
	$$\phi(\bar{x}): \exists \bar{y} \mesp \bar{x}.B -\bar{y}.A \in W.$$
	
	Let $I:=\{i \tq \gamma_i=\infty \}$ and $J:=\{j \tq \gamma_j\neq \infty \}$, 
	and $A_I$  be the matrix formed by the columns $C_{i \in I}$ (resp. 
	$A_J$ be the matrix formed by the columns $C_{j\in J}$ of $A$). We set
	$u(\bar{x}):=\bar{x}.B$ and  by $u_I(\bar{x})$ (resp.  $u_J(\bar{x})$) we denote the
	tuple formed by $I$-coordinates (resp. $J$-coordinates) of $u(\bar{x})$. We may assume that $\phi$ is of the form
	\begin{equation}\label{firstred}
		\exists \bar{y} \mesp \left( {u_I(\bar{x})} = \bar{y}.A_I \wedge u_J(x) -
		\bar{y}.A_J \in W_J \right)
	\end{equation}
	where $W_J$ is the obvious projection of $W$ to its non zero coordinates.
	By Lemma \ref{ddp6.4}, there exists a lower triangular separable matrix  $\widetilde{A}_I=(S,0)$ 
	such that the formula 
	$$u(\bar{x})=\bar{y}.A_I$$
	is equivalent, modulo $T_\lambda$, to 
	$$
	(t_1(u_I(\bar{x})),\dots,t_{n-l}(u_I(\bar{x})))=\bar{y}.PS \wedge (t_{n-l+1}(u_I(\bar{x})), \dots ,t_{n}(u_I(\bar{x})))=0,
	$$ 
	where $P$ is a permutation matrix and the $t_1(u_I(x)),\ldots t_n(u_I(\bar{x}))$ are some $L(\lambda)$-terms.
	
	By remark \ref{croissant}, chose a proper ball $U'$ such that $U'.A_J \subseteq W_J$
	and by Corollary \ref{bouledivisible} chose a proper ball $V'$ such that $V' \subseteq U'.PS$.
	
	By remark \ref{croissant} again, we choose $V$ such that for all 
	$\bar{a} \in V$,
	$$(t_1(u_I(\bar{a})),\dots,t_{n-l}(u_I(\bar{a}))) \in V'$$ and  $u_J(\bar{a}) \in W_J$.
	
	Hence:
	\begin{equation} T_{hens} \models \forall \bar{x} \mesp  \left(\phi(\bar{x})\wedge \bar{x}\in V\right) \longleftrightarrow \left((t_{n-l+1}(u_I(\bar{x})), \dots ,t_{n}(u_I(\bar{x})))=0\wedge \bar{x}\in V\right).
	\end{equation}
	
	Denote by $\psi$ the formula 
	$$(t_{n-l+1}(u_I(\bar{x})), \dots ,t_{n}(u_I(\bar{x})))=0.$$
	
	The first statement is now proved.

	Now if $\phi_1$ is another $p.p.$ formula, modulo $T^0_{free}$, $\phi_1$ is equivalent  to some positive quantifier free formula $\psi_1$. Since $T^{0}_{free}$ is decidable, we have an algorithm which checks if $T^0_{free} \models \psi \leftrightarrow \psi_1$. 
	Modifying this algorithm we can remember the finitely many
	non zero $r \in R$, such that algorithm uses the axiom
	$$\forall x \mesp  x\neq 0 \to  x.r\neq 0,$$ and the finitely many separable $s$, such that, the algorithm uses the axiom 
	$$\forall x \exists y \mesp y.s=x.$$ Choose $\gamma$ bigger than
	all the $\max \jump{r}$ and $\max\{\h(s),\hens(s)\}$ for all the $r$ and $s$ as above. Set $W:=P_\gamma^m$. Hence for any non zero $x \in P_\gamma$,
	$x.r\neq 0$ and there exists   $y \in P_\gamma$, with $y.s=x$. It follows that the same algorithm computes a proof of
	$ W\wedge \psi(\bar{x}) \leftrightarrow W\wedge \psi_1(\bar{x})$ from $T_{hens}$.  In particular
	we have $\phi_1 \wedge W \leftrightarrow W\wedge \phi(W)$ modulo $T_{hens}$.
	\begin{remark}
		By passing to an $\omega_1$-saturated model
		$M$, since $P_M\models T^0_{free} $, it is trivial that if $T^0_{free} \models \phi \leftrightarrow \phi_1$, then for some proper $W$, $$\phi(W) = \phi_1(M).$$ What we show above is that the decidability of $T^0_{free}$ yields the computability of $W$.
	\end{remark}
	
	\begin{corr}\label{corr:boule_up}
		Let $\phi$ and $\phi_1$ be $p.p.$ formulas of $L_{\mathcal{O}}$ and $\psi$ and $\psi_1$ are given as in the proof of the above theorem. Then, there is a computable ball $W$ such that 
		$$\vert (\phi\wedge W)/(\phi_1 \wedge W)\vert =1,$$ or 
		for all proper ball $V \subseteq W$, 
		$$\vert (\phi\wedge V)/(\phi_1 \wedge V)\vert = \infty.$$
	\end{corr}
	\begin{proof}
		Follows by Proposition \ref{prop:t0freedecidable} and by the above theorem.
	\end{proof}
	
	
	\section{pseudo-complements}
	We will introduce the notion of {\it valued module}  to study henselian filtered modules. These are  $R$-modules  $M$, equipped with a function $v:M\to \Gamma$ inducing the ultrametric topology.  After investigating elementary properties of valued modules, we will get the consequences that can be expressible in the language $L_{\mathcal{O}}$.
	
	\begin{defn}\label{def:valmod} A valued module  is a $t$-decomposable $R$-module together with a surjective map $v:M \to \Gamma$ such that for all $x,y\in M$, 
		
		\begin{enumerate}[1.]
			\item  $v(x\pm y)\geq \min\{v(x),v(y)\}$
			\item  $v(x)=\infty \leftrightarrow x=0$
			\item  $v(x)\notin \jump{r} \to v(x.r)=v(x)\cdot r$,  for all $r\in R$.
		\end{enumerate}\medskip
	\end{defn}
	
	\begin{remark}\label{rvR} Let $r=t^na + \dots + t^ka_k \in R$ where monomials are written following decreasing degrees,  then
		\begin{enumerate}[1.]
			\item $v(x.t^ia_i)=v(x)\cdot t^ia_i$, for all $x \in M$,
			\item  $v(x.r)=v(x)\cdot t^ka_k <v(x.(r-t^ka_k))$ whenever $v(x)>\max \jump{r}$,
			\item  $v(x.r)=v(x)\cdot t^na_n<v(x.(r-t^na_k))$ whenever $v(x)<\min \jump{r}$.
		\end{enumerate}
	\end{remark}
	\begin{proof}
		1. Follows from  Definition \ref{def:valmod} (3.) since a monomial has no jump value.
		
		2. By Definition \ref{def:valmod} (3.), if $v(x)>\max \jump{r}$ then $v(x.r)=v(x)\cdot r$. 
		Let $\gamma > \max \jump{r}$. Then for some $i$, $\gamma\cdot r=\gamma \cdot t^ia_i<\gamma \cdot t^ja_j$ for all $j\neq i$.  In other words the line  $\{(\delta,d^i\delta +v(a_i)\}_\delta$ does not intersect any other line $\{(\delta,d^j\delta +v(a_j)\}_\delta$ in the area $(\max \jump{r}, \infty]\times \Gamma$. This can happen only if $i<j$ for all $j\neq i$. Hence $i=k$.
		
		3. The proof is very similar to (2.)
	\end{proof}

	Let $(M,v)$ be a  valued module.
	We define the equivalence relation
	$\mathbf{RV}$ on $M$ by 
	$$x \mathbf{RV} y \iffl v(x)=v(y)<v(x-y).$$ We denote the $\mathbf{RV}$-class of an element $x$ by $\rv{x}$ whereas $\rv{A}$ stands for the set $\{\rv{x} \tq x\in A\}$ for $A\subseteq M$.
	We also set the notation $P_\gamma$  for the closed ball of radius $\gamma$ centered at $0$.
	
	\begin{defn} For subgroups $A$ and $B$ of $M$, we say that $A$ and $B$
		\begin{itemize}
			\item   are $m$-immediate (m stands for \emph{mutually}) if $\rv{A\setminus P_\gamma}=\rv{B\setminus P_{\gamma}}$ for some $\gamma \in \Gamma$, and we write
			$A \approx B$,
			\item   are pseudo-orthogonal if $vA\cap vB \subseteq [\gamma, \infty]$ for some $\gamma \in \Gamma$, and we write $A | | B$.
		\end{itemize}\medskip
	\end{defn}
	\begin{remark}
		$A\cap C \subseteq P_{\gamma}$ for some $\gamma$, whenever $A || C$.
	\end{remark}
	\begin{remark}\label{rem:modulogamma} If $A\approx B$ and $B\subseteq A$ then for some $\gamma$,
		$A+P_\gamma=B+P_\gamma$.
	\end{remark}
	\begin{proof}
		Let $\gamma$ be such that  $\rv{A\setminus P_\gamma}=\rv{B\setminus P_{\gamma}}$.
		It is enough to show $\{x\in A \tq v(x)\leq \gamma\} \subseteq B+P_\gamma$.
		We proceed by induction w.r.t. the dual order on the initial segment 
		$(-\infty ,\gamma]$:
		Let $a\in A$. If $v(a)=\gamma$ there is nothing to do. 
		Suppose $v(a)<\gamma$ and for all $a' \in A$ of valuation $>v(a')$ there is some
		$b' \in B$ such that $a-b'\in P_\gamma$. Since $A\approx B$, there is $b\in B$ such that
		$v(a-b)>v(a)$. Since $b\in B$ and $B\subseteq A$, $a-b \in A$. Now by applying 
		the induction hypothesis to $a-b$ we have $a-b-b'\in P_\gamma$ for some $b'\in B$. That is, $a \in B+P_\gamma$.
	\end{proof}
	\begin{defn} A pseudo-complement of a subgroup $A$, is a subgroup $C$ such that,
		for some $\gamma$ 
		\begin{equation}
			M=A + C + P_\gamma \quad \text{and} \quad A|| C.
		\end{equation}
	\end{defn}

	\begin{remark}\label{rem:invimage}
		If $A,C$, and $P_\gamma$ are as above and $f:M\to M$ is an additive map such that $f^{-1}(P_\gamma) \subseteq P_\delta$ for some $\delta$, then $f^{-1}(C)$ is a pseudo-complement to $f^{-1}(A)$. In particular this is the case when 
		$f$ is given by a scalar multiplication.
	\end{remark}
	
	\begin{remark}\label{rem:+ballpseudocomp}
		If $C$ is a pseudo-complement of $A$ then $C$ is a pseudo-complement
		of $A+P_{\gamma}$ for any $\gamma$.
	\end{remark}
	
	\begin{remark}\label{rem:samecomplement} It is straightforward to see that if $A$ and $B$ have the same pseudo-complement then $A \approx B$. The following lemma establishes the converse using that the value set is $\mathbb{Z}\cup \{\infty\}$.
	\end{remark}
	\begin{lemma}\label{changecomp} If $A\approx B$ then $C$ is a pseudo-complement of $A$ if and only if it is a pseudo-complement of 
		$B$.
	\end{lemma}
	\begin{proof}
		Suppose $C$ is a pseudo-complement of $A$  satisfying
		$$   M=A + C + P_\gamma$$
		and $A\approx B$. Let $\delta$ be such that  $\rv{A\setminus P_\delta}=\rv{B\setminus P_\delta}$. We may assume that  $\delta \leq \gamma$. We claim that 
		$$M=B+C+P_\delta .$$
		
		We proceed by induction on $\{\beta \in \Gamma \tq \beta\leq \delta\}$ as in the proof of Remark \ref{rem:modulogamma}. If $v(z)=\delta$ then trivially $z \in B + C + P_\delta$. Let $z \in M$ and suppose that for all $x$ with $v(x)>v(z)$ there exist 
		$(b,c,x_\delta) \in B\times C\times P_{\delta} $ such that 
		$$x=b+c+x_{\delta}.$$
		Write $z=a+c+z_\delta$ for some $(a,c,z_\delta) \in A\times C \times P_\delta$ (in fact $z_\delta$ can be chosen in $P_\gamma \subseteq P_\delta$). If $v(a)\geq \delta$ there is nothing to do. If $\delta>v(a)>v(z)$ then by induction hypothesis
		$a=b+c'+z'_\delta$ with $(b,c',z_\delta) \in B\times C \times P_\delta$, and hence 
		$$z=b+c+c'+z_\delta + z'_\delta. $$ The only possibility which remains to be considered is $v(z)=v(a)<\delta$ since $A || C$ implies that
		$v(a+c+z_\delta)=\min\{v(a),v(c)\}$ whenever $v(a)<\delta$. Pick  $b\in B$ such that
		$\rv{b}=\rv{a}$. Then, since $v(a-b)>v(a)$, by induction hypothesis the equality
		$$a-b=b'+c'+z'_\delta$$ holds for some $ (b,c',z_\delta') \in B\times C \times P_\delta$. Hence $z \in B+C+P_\delta$. 
		
		Now to see that $B || C $ we claim that $vB \cap vC \in [\delta, \infty]$. Suppose for a contradiction that
		for some $b \in B $ and $c\in C$, $v(b)=v(c)<\delta$. Then we can choose $a \in A$ such that
		$v(a)=v(b)$ hence $v(a)=v(c) \in vA \cap vC$. But then $v(a) \in [\gamma, \infty]$. 
		This is a contradiction since $\gamma \geq \delta$.  
	\end{proof}
	\begin{lemma} Let $A,A',B,B'$ be  such that 
		$A \approx A'$,  $B\approx B'$ and $A || B$. Then the following holds: 
		\begin{enumerate}[1.]
			\item  $A' || B'$, 
			
			\item  every pseudo-complement of $A+B$ is a  pseudo-complement of $A'+B'$.
		\end{enumerate}\medskip
	\end{lemma}
	\begin{proof}
		1. Let $\gamma=v(a')=v(b')$ with $a'\in A$ and $b'\in B$. Let $\delta$ be such that $vA\cap vB \subseteq [\delta,\infty]$ and $\gamma_1$, $\gamma_2$ are respectively the values yielding $A\approx A'$ and $B\approx B'$. We claim that  
		$\gamma\geq \min\{\delta,\gamma_1,\gamma_2\}$. Let $a \in A$ and $b\in B$ be such that $\rv{a}=\rv{a'}$ and
		$\rv{b}=\rv{b'}$. Then $\gamma=v(a)=v(b)$. Hence $\gamma$ must be $\geq \delta$.

		2. Let $C$ be a pseudo-complement of $A+B$. Let $\delta,\gamma_1,\gamma_2$ be as above. Let
		$z=a+b$ with $a\in A$ and $b \in B$ such that $v(z)<\min\{\delta,\gamma_1,\gamma_2\}$. Then $v(a)\neq v(b)$ and $v(z)=\min\{v(a),v(b)\}$ necessarily. Suppose $v(a)=v(z)$.
		Then $\rv{a}=\rv{z}$ hence for some $a'\in A'$, $\rv{z}=\rv{a'}$.  If $v(z)=v(b)$ we can choose  in the same way $b'$ such that $\rv{b'}=\rv{z}$. Hence $A+B \approx A' + B'$.  We apply now Lemma \ref{changecomp}.
	\end{proof}
	
	\begin{defn}
		A valuation independent basis $\beta$ of the $K^{\vfi^{n}}$-vector space $K$,  is a basis such that the members of $\beta$ have all different valuations in the finite set $\{0,\ldots,d^{n}-1\}$. 
	\end{defn}
	For example, the basis  $\alpha$ of the $K^{\vfi}$-vector space $K$ is valuation independent.
	\begin{remark}
		If $\beta$ is a valuation independent basis of the $K^\vfi$-vector space $K$ then for all $s>0$, $\beta(n)$ is a valuation independent basis of $K^{\vfi^n}$-vector space $K$. Moreover,  any valued module is 
		$t\beta$-decomposable for all valuation independent $\beta$.
	\end{remark}
	
	
	\begin{lemma}
		Let $g_1,\dots,g_m$ be all of the same degree $s$  such that the leading coefficients $b_1, \dots, b_m$ have distinct valuations  in
		$\{0,\dots,d^s-1\}$. Then $M.g_i || M.g_j$ whenever $i\neq j$. 
	\end{lemma}
	\begin{proof}
		By \ref{rvR}, we have $v(x.g_i)=v(x)\cdot g_i=v(x)\cdot t^sb_i=p^{d^s}v(x)+v_K(b_i)$ whenever $v(x)<\jump{g_i}$ for $i=1\dots m$. In particular 
		if $v(x)<\min\{\jump{g_i},\jump{g_j}\}$ for a $j\neq i$ then $v(x.g_i)\neq v(x.g_j)$.
	\end{proof}
	\begin{lemma}\label{lem:addtoMta} Let $\beta$ be a valuation independent basis of the $K^{\vfi}$-vector space $K$. Let $q=t^sa+\dots \in R$ be of degree $s$. 
		Then there is a unique $j=j(q)$ such that $\beta_j \in \beta(s)$ and $M.q \approx M.t^s\beta_j$. As a consequence $C:=\sum_{j' \neq j(q)} M.t^s\beta_{j'}$ is a pseudo-complement for $M.q$.
	\end{lemma}
	\begin{proof}
		Write
		$$a=\sum_{j \in d^s} a_j^{\varphi^s}\beta_j. $$
		
		Let $j$ be such that $v_K(a)=v_K(a_j^{\vfi^s}\beta_j)<v_K(a_{j'}^{\vfi^s}\beta_{j'})$ for all $j' \neq j$ (such a $j$ exists since $\beta$ is valuation-independent). 
		Then $$v(x.t^sa)=v(x)\cdot t^s(a_j^{\vfi^s}\beta_j)<v(x)\cdot t^s(a_{j'}^{\vfi^s}\beta_{j'})$$ for all $x\in M$ and $j'\neq j$. 
		
		Let $\gamma < \min
		\jump{q}$. Then by Remark \ref{rvR}, for all $x$ of such that $v(x)\leq \gamma$, we have $\rv{x.q}=\rv{x.t^sa}$, hence 
		$$\rv{x.q}=\rv{x.t^sa_j^{\vfi^s}\beta_j}=\rv{(x.a_j).t^s\beta_j}.$$ Since $x\mapsto x.a_j$ is a bijection it follows that $M.q \approx M.t^s\beta_j$. The consequence follows by Lemma \ref{changecomp}.
	\end{proof}

	\begin{lemma}\label{vddku}
		Given a non zero matrix $Q$ with coefficient over $R$, say with $k$ rows, there exists a matrix $Q'$ such that
		the first column of $Q'$ consists of polynomials which have all the same degree $s$,
		such that the leading coefficients in this column have distinct valuations  in $\{0,\ldots, d^s-1\}$  and
		$$M^k.Q =  M^k.Q' \mesp .$$
	\end{lemma}
	\begin{proof} This proof is essentially a slight generalization  of the proofs of 
		Lemma 3 and Lemma 4 in \cite{vdd-kuhlmann}.
		
		Let $Q=(q_{ij})$ be a matrix over $R$ with $k$-many rows.
		We will proceed by induction on $f=\sum_{\{(i,j)\tq q_{ij}\neq 0\}}(\deg(q_{ij})+1)$.
		
		Since $Q$ is non zero, $f>0$. Suppose $f=1$. We may  assume 
		$q_{11}=c \in K^{\times}$ and all other entries of $Q$ are zero. Then 
		$\bar{x}=\bar{y}.Q$ for some $\bar{y}$ if and only if, all the coordinates of $\bar{x}$ except possibly the first one, are $0$. Hence we can take $Q'$ the matrix which has $1$ at the position
		$(1,1)$ and has all other entries equal to zero.
		
		Now we suppose $f>1$. Let $e:=\max\{\deg{q_{ij}}\}\geq 0$. We may suppose that $q_{11}$ has degree $e$.
		Set $$e_{ij}:=\deg{q_{ij}},  \mesp e_i:=e_{i1}$$ and 
		$$c_{ij}:=\text{the leading coefficient of}\mesp q_{ij},  \mesp c_i=c_{i1}.$$
		
		\noindent
		{\bf Claim 1}: We may assume that for all $a_1,\ldots a_k \in K$, not all are $0$,
		$\sum_i a_i^{\vfi^{e_i}}c_i \neq 0$.
		
		Suppose $\sum_i a^{\vfi^{e_i}}c_i = 0$. We may also suppose that $a_1=1$.
		We define for all $j$,  
		$$\widetilde{q_{1j}}=\sum_{i=1}^{k} t^{e-e_{i}}a_iq_{ij}.$$
		
		Since $a_1=1$ and $e=e_1$ 
		\begin{equation}\label{eq:matricepassage}
			\widetilde{q_{1j}}=q_{1j} + t^{e-e_2}a_2q_{2j} + \dots + t^{e-e_k}a_kq_{kj}.
		\end{equation}  
		We also define $q'_{ij}$ by the equality
		$$q_{ij}=t^{e_i}c_{ij} + q'_{ij}.$$
		
		We claim that $\widetilde{q_{11}}$ has degree  $<e$: We have
		\begin{equation}
			\widetilde{q_{11}}=\sum_{i=1}^k t^{e-e_i}t^{e_i}a_i^{\vfi^{e_i}}c_i + t^{e-e_i}a_iq'_{i1} = \sum_{i=1}^k t^ea_i^{\vfi^{e_i}}c_i + t^{e-e_i}a_iq'_{i1}.
		\end{equation}
		Since each $q'_{i1}$ has degree $<e_i$, each  $t^{e-e_i}a_iq'_{i1}$ has degree $<e$. Hence the coefficient $t^e$ in $\widetilde{q_{11}}$ is 
		$\sum_i a_i^{\vfi^{e_i}}c_i=0$ and $\deg (\widetilde{q_{11}})<e$.

		Let $\tilde{Q}$ be the matrix where we have replaced $q_{1j}$ by  $\widetilde{q_{1j}}$. Now the sum of degrees of the non zero entries of $\tilde{Q}$ is less than the sum of degrees of the non zero entries of $Q$. Hence in order apply induction it is enough to prove that the solvability of the system
		$\bar{x}=\bar{y}.Q $ is equivalent to the solvability of $ \bar{x}=\bar{z}.\tilde{Q}$. 
		But this follows by expressing the equations (\ref{eq:matricepassage}) by the equality
		$$\tilde{Q}=PQ$$
		where $$P=\begin{pmatrix} 
		1 & t^{e-e_2}a_2 & t^{e-e_3}a_3 & \cdots & t^{e-e_k}a_k \\
		0 & 1            & 0            & \cdots &       0       \\
		\vdots & 0 & \ddots \\
		0 & \cdots & \cdots & \cdots & 1
		\end{pmatrix}					  $$ 
		is clearly invertible in $R$.
		
		{\it Claim 2}: Assume Claim 1. We may assume that the polynomials 
		$(q_{i1})_i$ of the first  column have same degree and leading coefficients
		of the $(q_{i1})_i$ are $K^{\vfi^e}$ linearly independent.
		
		We will show that we can change $Q$ to some $S$, possibly having more rows, such that 
		the system $\bar{x}=\bar{y}.Q$ is equivalent to $\bar{x}=\bar{z}.S$ with the first column of $S$ has required properties.

		Recall that $q_{i1}$ has degree $e_i$ with $e_1=e$. For all  $1\leq i \leq k$,  using the basis $\alpha(e-e_i)$, we write the equality of terms 
		$$y_i=\sum_{u \in d^{e-e_i}}\lambda_u(y_i).t^{e-e_i}\alpha_u .$$
		Now $$y_i.q_{i1}=y_i.t^{e_i}c_i + y_i.q'_{i1} = \sum_u \lambda_u(y_i).(t^e\alpha_u^{\vfi^{e_i}}c_i + r_{i1}(u))$$
		where $r_{i1}(u)= t^{e-e_i}\alpha_uq'_{i1}$ which has degree $<e$. Set $z_i(u)=\lambda_u(y_i)$ and 
		$$s_{i1}(u):=t^e\alpha_u^{\vfi^{e_i}}c_i + r_{i1}(u).$$ 
		Since the leading coefficient of $s_{i1}(u)$ is $\alpha_u^{\vfi^{e_i}}c_i$, for any $i$ and $u$, 
		$s_{i1}(u)$ has the degree  $e$.
		
		For $j>1$ and $u\in d^{e-e_i}$, set $s_{ij}(u):=t^{e-e_i}\alpha_u.q_{ij}$.   
		Note that we keep $e_i=e_{i1}$ but $j$ varies. Let $S$ be the matrix obtained from $Q$, by replacing 
		$i$-th row by the matrix $(s_{ij}(u))_{u,j}$ where $u$ is the row-index and $j$ is the column index.
		Then the system $(x_j=\sum_i y_i.q_{ij})_j$ is equivalent to the system
		$(x_j=\sum_i \sum_u \lambda_u(y_i).s_{ij}(u))_j$, which can be written as
		$$\bar{x} = \left((\lambda_u(y_1))_u,\ldots,(\lambda_u(y_k))_u \right).S .$$
		
		Now we will show that the leading  coefficients $\alpha_u^{\vfi^{e_i}}c_i$'s are 
		$K^{\vfi^{e}}$-linearly independent.
		
		Suppose $$\sum_i \sum_u a_{u_i}^{\vfi^e}\alpha_u^{\vfi^{e_i}}c_i= 0$$ for some
		tuple $(a_{u_i})_{u,i}$ from $K$. It follows   by Claim 1 that, 
		$\sum_u (a_{u_i}^{\vfi^{e-e_i}}\alpha_u)=0$ for each $i$. Since the $\alpha_u$ are $K^{\vfi^e}$-linearly independent and $\vfi$ is injective
		$a_{u_i}=0$ for all $i,u$. The Claim 2 is  proved.
		
		We assume now that the first column of $Q$ consists of 
		polynomials of degree $e$ with leading coefficients being $K^{\vfi^e}$-linearly independent.
		By section 3  and by the last paragraph of Lemma 4 of \cite{vdd-kuhlmann}, there exists an invertible matrix over $P$ over $K$
		such that $P.Q_1$, where $Q_1$ the first column of $Q$, consists of polynomials with leading coefficients has all different valuations in 
		$\{0,\ldots,d^{e}-1\}$. Hence considering $P.Q$ finishes the proof.
	\end{proof}

	For the following lemma, we will use the above result with $Q$ a column matrix and then in the following corollary we will use it in whole generality.
	
	\begin{lemma}\label{lem:pseudocompbeta}
		Let $A \subseteq M$, of the form $A=\sum_i M.q_i$ then for some integer $s$, $A$ has a pseudo-complement of the form 
		$\oplus_{i\in I} M.t^sb_i$ where $I\subseteq d^s$, and the $b_i$ are valuation independent.  
	\end{lemma}
	\begin{proof}
		Use the above lemma to chose $g_j$ such that $\sum_j M.g_j =A$,  all of degree $s$ with leading coefficients $b_j$'s have different valuations in $\{0,\dots,d^s-1\}$. Then
		by the lemma \ref{lem:addtoMta}, $M.g_j \approx M.t^sb_j$ for all $j$. Complete the $b_j$ to a valuation independent basis of $K^{\vfi^s}$-vector space $K$. We write the new elements of this basis as  the $b_i$.    Let $C:=\sum_{i} M.t^sb_i$. 
		Since $C\oplus \sum M.t^sb_j=M$, $C$ is in particular a pseudo-complement for $\sum M.t^jb_j$. Hence by Lemma
		\ref{changecomp}, $C$ is a pseudo-complement for $A$. 
	\end{proof}
	\begin{remark}\label{rem:pseudocomplements} 
		The pseudo-complement $C$ is $p.p.$ definable  subgroup by an $L$-formula $C(x)$, which does not depend on $M$. In other words, if $\phi(x)$ is the formula 
		$$\exists y_1, \dots , y_m \mesp x=\sum_i y_i.q_j$$ then in any valued module $(M,v)$,
		$C(M)$ is a pseudo-complement to $\phi(M)$.
		
	\end{remark}

	\begin{theorem}\label{pseudocomplements}
		Let $\phi(\bar{x})$ be a $p.p.$ formula of $L_\mathcal{O}$ of the form 
		$$\bar{x} - \bar{y}.Q \in W,$$ where $Q$ a matrix with coefficients from $R$ and $W$ a ball. Set $k:=\vert \bar{x} \vert$. Then there exists computable $\gamma$ and an existential $L$-$p.p.$ formula $D(\bar{x})$ such that 
		$$M^k=\phi(M^k) + D(M^k) + P_\gamma^k$$ $\&$
		$$D(M^k)\cap \phi(M^k) \subseteq P_\gamma^k$$
		for all valued module $(M,v)$.
	\end{theorem}
	\begin{proof}
		
		By Lemma \ref{vddku}
		we may assume that the first column $Q_1$ of $Q$ consists of polynomials having same degree $s$ with leading coefficients having distinct valuations in $\{1,\dots, d^s-1\}$. In addition, by Remark \ref{rem:+ballpseudocomp},
		we may suppose that $W=\{0\}^k$.
		
		Take $C$ a pseudo-complement to $M.Q_1$ as in the above lemma. 
		Let $x=(x_1, \dots, x_k)$ such that $x_1 \in C \cap M.Q_1$. Then $x_1 \in 
		P_\delta$ for some computable $\delta$. Write $x_1$ also as
		$$x_1=y_1.q_{11} + y_2.q_{21} + \dots y_mq_{m1}.$$
		Since the leading coefficients of the $q_{1j}$ have different valuations
		in  $\{0,\dots, d^s-1\}$,  if 
		$$v(y_j)\leq \rho:=\min\{\min_{ij} \{\jump{q_{ij}}\}\}-1;$$ we have $v(y_j.q_{1j})\neq v(y_{j'}.q_{1j'})$ for $j\neq j'$. Hence the $y_j$ can not have indefinitely small valuations since $v(x_1)\geq \delta$. Now  for $i>1$, since $x_i=\sum_{j}y_j.q_{ij}$ we have 
		$v(x_i)\geq \min_j \{\rho\cdot q_{ij}, \delta\}$. Setting $\gamma:=\min_{ij}\{\rho\cdot q_{ij}\}$ and $D:=C\times M^{k-1}$ yields our claim since $D$ and $\gamma$ depends only to $Q$.
	\end{proof}
	
	\begin{remark}\label{boundedsolutions}
		The above proof shows that whenever a matrix $Q$ with the first column  consists of polynomials whose dominant coefficients are valuation independents and $\gamma \in \Gamma$, as above are given, there is a computable $\delta \in \Gamma$, such that   
		$$\bar{y}.Q \in P_\gamma^m \Rightarrow \bar{y} \in P_\delta^k  .$$  
	\end{remark}
	\subsection*{\bf  The theory $T_{\Psi}$}

	Let $(M,v)$ be a valued module and set 
	$$A:=\sum M.q_i$$
	where the $q_i$ are all of degree $s$ such that for a  valuation independent basis $\beta$  of $K^{\vfi}$-vector space $K$, the leading coefficients of the $q_i$ are from the basis $\beta(s)$. Let $\mathbb{A}(\beta)$ be the tree consisting of subgroups $M.t^i\beta_i(n) \, (i,n \in \omega)$, ordered by inclusion. At the level $n$ we have the subgroups 
	$M.t^i\beta_i(n) \mesp (i \in d^n)$. 
	
	\begin{defn}
		We call the the pseudo-complement for $A:=\sum_j M.q_j$ 
		the unique pseudo-complement which can be written as the sum of some elements of level
		$\deg(q_i)$ of the tree $\mathbb{A}(\beta)$.
	\end{defn} 
	This definition only depends on  the decomposition 
	of the leading coefficients of the $(q_i)_i$ in the basis $\beta(n)$. Hence, given $A$ as above, by Lemmas \ref{vddku} and \ref{lem:pseudocompbeta}, there is a recursive function 
	\begin{equation}\label{metaf}
		\mathfrak{f}: \bigcup_n K^n \to \bigcup_n K^n
	\end{equation} 
	which computes the basis $\beta$ and the pseudo-complement $C$ of  in every valued 
	$t\beta$-decomposable $R$-module $M$. We  write 
	$C_{\mathfrak{f}(q_1,\ldots,q_n)}$ for the pseudo-complement computed by $\mathfrak{f}$.
	
	On the other hand, the ball $P_\gamma$ such that $M=A+C + P_\gamma$, can be chosen 
	by letting $\gamma:=\min\{\min_i\{\jump{q_i}\} \} - 1$. This yields another recursive function
	\begin{equation}\label{metaj}
		\mathfrak{j}:  \bigcup_n K^n \to \Gamma.  
	\end{equation} 
	
	Hence we can express the statement of Corollary \ref{pseudocomplements} by introducing an axiom $\psi(q_1,\ldots, q_n)$  for any matrix $Q$ which has the first column  $(q_,\ldots,q_n)$, in the language 
	$L_{\mathcal{O}}$, which says that 
	\begin{equation}
		M^k=M^k.Q + C_{\mathfrak{f}(q_1,\ldots,q_n)}\times M^{k-1} + P_{\mathfrak{j}(q_1,\ldots,q_n)}^k \mesp \quad \text{and}
	\end{equation}
	$$ M^k.Q \cap  C_{\mathfrak{f}(q_1,\ldots,q_n)}\times M^{k-1} \subseteq P_{\mathfrak{j}(q_1,\ldots,q_n.)}^k.$$

	Let $T_\Psi$ be the $L_\mathcal{O}$-theory of  $R$-modules together with the sentences  $\psi(q_1,\ldots, q_n)$. Hence
	$T_\Psi$ is recursively enumerable.

	\begin{theorem}\label{approx}
		Let $Q$ be a  $m\times k$  matrix over $R$, $W=\prod_{i=1}^k P_{\gamma_i}$  and $\phi$ 
		be the $L_{\mathcal{O}}$-formula
		$$\phi(x_1,\ldots,x_k): \exists y_1\ldots y_m \mesp (x_1,\ldots x_k) - (y_1,\ldots, y_m).Q \in W.$$  Then for some  computable $\gamma$, and some positive primitive  $L_\mathcal{O}$-formula $D(x)$ we have  
		$$M^k =  \phi(M^k) + D(M^k) + \prod_{i=1}^k P_\gamma   \quad \& \quad D(M) \cap \phi(M^k)\subseteq \prod_{i=1}^k P_\gamma.$$ 
		for all $M\models T_{\Psi}$.
	\end{theorem}

	\section{Decidability and model completeness of $\mathbb{F}_d((X))$}
	We will introduce a new theory $T_1$, containing $T_{\Psi}$, augmented by sentences counting the number of solutions of the $p.p.$ formula 
	$$\bar{x}-\bar{y}.Q \in B$$ in $B_0$ modulo $B_1$, for proper balls $B_1 \subseteq B_0$.
	
	For  $Q$, an $m\times n$ matrix over $R$, such that its first column consists of 
	polynomials whose leading coefficients are valuation independent and $B$ a ball, set $$A(\bar{x}):= \exists  \bar{y} \mesp \bar{x}-\bar{y}.Q \in B.$$ Let $\delta$ be the value computed by 
	Theorem \ref{boule_up}, such that $A\wedge P_\delta^m$ is quantifier-free definable
	in the language  $L(\lambda)$. 
	
	Let    $B_0:=P_\gamma^m$, where $\gamma$ is given by Theorem \ref{approx}.  $B_1:=P_\delta^m$,  and 
	$$k:=\vert (A(K)\cap B_0(K))/A(K)\cap B_1(K)\vert$$ if $\gamma\leq \delta.$ We set the sentence  $\theta(A,\gamma,\delta)$ expressing 
	\begin{equation}
		k=  \left| \frac{(A\wedge B_0)(K)}{(A \wedge B_1)(K)}\right| .
	\end{equation}
	Let  $\Theta$  be the set of sentences  $\theta(A,\gamma,\delta)$ and we set 
	$$T_1:=T'_{hens}\cup T_{\Psi} \cup \Theta$$
	where $T'_{hens}$ is the $L_\mathcal{O}$-theory composed by the axioms of $T_{hens}$, where we have replaced any $L(\lambda)$-term by its equivalent modulo the theory of $R$-modules in language $L$ (recall Remark \ref{extunideflambdas}). 
	
	Note that $T_1$ implies that $P_{\gamma +1}/P_{\gamma}$ has exactly $d$ elements for all $\gamma \neq \infty$.
	\begin{remark} $\Theta$ and hence $T_1$ is a recursively enumerable theory. In fact, by Remark \ref{boundedsolutions}, if $\bar{x},\bar{y},\gamma$ are such that
		$\bar{y}.Q=\bar{x}\in P_{\gamma}^n$,  then $\bar{y}$ is in some $P_{\theta}^k$ for a computable $\theta$, hence searching the solutions 
		$\bar{y}$, of $\bar{y}.Q=x$ with $\bar{x} \in P_\gamma^n$ can be bounded to searching $y$'s in some ball. Hence
		searching such solutions modulo another ball can be done in some finite $\mathbb{F}_d$-vector space by an algorithm.
	\end{remark}
	\begin{prop}\label{laurentsatisfyT1}
		Let $M:=\mathbb{F}_d((X))$ then as an $L_\mathcal{O}$-structure $M\models \Theta$  and hence $M \models T_1$.
	\end{prop}
	\begin{proof}
		By \cite{kuhlmanntame}, Theorem 5.14, $K$ is existentially closed as a ring  in $M$. In addition, by \cite{sylvyarno17} Corollary 6.18, there exists an  existential ring-formula
		without parameters which defines the maximal ideal both in $K$ and in $M$. Since the valuation ring is the complement of the set of inverses of the elements in the maximal ideal, we have a universal ring-formula  which defines uniformly the valuation ring in $K$ and in $M$. Hence the balls 
		centered at $0$ are definable  universally with the parameter  $X$ both in $K$ and $M$.

		Suppose $\vert (A(K)\cap B_0(K))/(A\cap B_1(K))\vert = k$. Consider the sentence $$\sigma: \vert (A\wedge B_0)/(A\wedge B_1)\vert \geq k.$$
		Since by theorem \ref{boule_up}$, (A(K)\cap B_1(K))$ is definable both universally
		and existentially in $L_{\mathcal{O}}$,  and and the valuation ring is universally definable both in $K$ and $M$,  $\sigma$ is equivalent
		to an existential ring-formula with parameters in $K$. Since this quotient is finite, we must have $\vert (A(M)\cap B_0(M))/A(M)\cap B_1(M)\vert = k$.
	\end{proof}
	
	To prove Theorem \ref{th:pp2universal}, we will use  a lemma from Rohwer's thesis 
	(\cite{rohwer}, Lemma 8.2). This lemma is a generalization of the following fact:
	
	\begin{obs} In an abelian group $G$ with existentially definable subgroups $A,B$ such that $A+B=G$, if $A\cap B$ is definable by a universal formula then $A$  is definable by the following universal formula $\psi(x)$:
		
		\begin{equation}
			\psi(x): \forall y \mesp (x-y \in A \wedge y \in B) \rightarrow y \in A\wedge B.
		\end{equation}
		By iterating this observation,  we have the following lemma.
	\end{obs}
	\begin{lemma}[Rohwer]
		Let $\mathcal{T}$ be a theory expanding  the theory of abelian groups. For each
		$M \models \mathcal{T}$ and for  $A, A_c, A_m, A_s, B_0, B_1$  definable subgroups of $M$ satisfying  the following configuration,
	\end{lemma}
	\begin{enumerate}[1.]
		\item $A + A_c=M$,
		\item $A\cap A_c \subseteq B_0$,
		\item $A\cap B_1 = A_s\cap B_1$,
		\item $A\cap B_0 \subseteq A_m \subseteq A + B_1$.
	\end{enumerate}\medskip
	where $A,A_c,B_1$ are definable by existential formulas, and $A_m, A_s$ 
	by universal formulas (where all formulas in question do not depend on $M$), $A$ is definable by a universal formula (which does not depend on $M$).   
	\begin{proof} See \cite{rohwer} Lemma 8.2. 
	\end{proof}
	
	\begin{theorem}\label{th:pp2universal}
		Any $p.p.$ formula of $L_{\mathcal{O}}$ is equivalent modulo $T_1$ to a universal $L_\mathcal{O}$-formula.   
	\end{theorem}\label{modcomplete}
	\begin{proof}
		Let $$A(\bar{x}): \exists \bar{y} \mesp \bar{x}.S -\bar{y}.Q \in W$$ be a $p.p.$ formula. For our purposes, we may assume that 
		$S$ is the identity matrix since we way replace $A(x)$ by
		$$\forall \bar{z} \mesp (\bar{z}=\bar{x}.S \to \exists \bar{y} \mesp \bar{z}-\bar{y}.Q \in W).$$  
		By Lemma \ref{vddku}, we may suppose that the first column of $Q$ consists of polynomials which have coefficients
		in a valuation independent basis. 
		By Theorem \ref{boule_up}, there is some proper ball $B_1$ determined by $T_1$ such that the formula $A_s:=A\wedge B_1$ is equivalent 
		to  a quantifier free $L(\lambda)$-formula. Hence by Remark \ref{extunideflambdas}, $A_s$ is equivalent to a universal
		$L$-formula. Now, let  $B_0$ be given by Theorem \ref{approx} such that
		$$A \wedge A_c \to B_0$$ where $A_c$ is of the form $A_c=D+P_\gamma^l$ with $D$    as in Theorem \ref{approx}.
		
		\emph{Claim} :  Set $A_m:=(B_0 \wedge  A) + A_s$. Then $A_m$ is equivalent to a  universal 
		$L_{\mathcal{O}}$-formula modulo $T_1$.
		
		\underline{proof of the claim}:
		Let $M \models T_1$. Let $k$ be the cardinality of $$(B_0(M) \cap A(M))/
		A_s(M)$$ and $y_1, \dots, y_k$ ($y_i$ are tuples of variables)  be  representatives of the classes.
		Note that $k$ is determined by a sentence 
		in $\Theta$, hence depends only on the theory $T_1$.
		Then $$\{y_1, \dots y_k \}  + A_s(M) = A(M)\cap B_0(M)+A_s(M).$$ Moreover,
		for all $z_1,\dots,z_k \in B_0$ satisfying $z_i-z_j \notin A_s(M)$ for $i\neq j$,
		$$\{ z_1, \dots z_k \}  + A_s(M) = A(M)\cap B_0(M) + A_s(M)$$
		if and only if,  
		$$\{z_1, \dots z_k \}  + A_s(M) \subseteq A(M) \cap B_0(M) + A_s(M).$$ Hence 
		the formula 
		\begin{multline}
			x\in B_0 \wedge \forall y_1, \dots, \forall y_k \\ (
			(\bigwedge_{i=1}^k y_i \in B_0 ) \wedge ( \bigwedge_{i\neq j} y_i - y_j \notin A_s) \wedge  \\
			[\exists_{i=1}^{k}z_i (\mesp \bigwedge_{i=1}^k z_i\in A \wedge \bigwedge_{i=1}^k y_i - z_i \in A_s)]) \longrightarrow \bigvee_{i=1}^k y_i -x \in A_s.
		\end{multline}
		is equivalent to $A_m$, which is equivalent to a universal $L_\mathcal{O}$-formula, thus the claim is proved.
		
		Now $A,A_c,A_m,A_s,B_0,B_1$ are in the configuration of Rohwer's Lemma above.
	\end{proof}
	
	\begin{corr}\label{decidable}
		Every completion of $T_{1}$ is model-complete in language $L_{\mathcal{O}}$.  In particular the complete $L_\mathcal{O}$-theory of $\mathbb{F}_d((X))$ is model-complete.
	\end{corr}
	\begin{proof}
		It follows by Corollary \ref{corr:baur-monk}.
	\end{proof}
	Note that 
	$K$ embeds  (via an $L_{\mathcal{O}}$-embedding) to any model $N\models T_{1}$. In fact, choose any $s \in N\setminus\{0\}$ such that $s.t=s$. Consider 
	the $L$-embedding $k\mapsto s.k$. We show that it is an $L_{\mathcal{O}}$-embedding:  It is easy to see that the assertion
	$$ \mathcal{O} = A \oplus \mathcal{O}.X$$ where  $A$ 
	is the formula $x.t=x$, is a consequence of $T_{1}$ and it is clear that 
	if $a\in \mathcal{O}_K$ then $s.a \in \mathcal{O}_N$. Suppose now that 
	$s.a \in \mathcal{O}_N$ for some $a \in K$. Then
	$s.a=y  +  m$ where $y.t=y$ and $m\in \mathcal{O}_N.X$. Then
	$(s.a).(t-1)=m.(t-1) \in \mathcal{O}_N.X$. Note that 
	$(s.a).(t-1)=s.ta^{\vfi} -sk=s.(a^{\vfi} -a) \in \mathcal{O}.X$ since $s.t=s$.
	This can only happen if $(a^{\vfi} -a) \in \mathcal{O}_K$, only if 
	$a\in \mathcal{O}_K$.

	By Corollary \ref{corr:baur-monk}, we have:
	\begin{corr}\label{K}
		The models of $T_1$ in which $K$ is existentially closed as an $L_\mathcal{O}$-structure are elementary equivalent to $K$. 
	\end{corr}
	
	\begin{proof}
		
		Let $N\models T_{1}$ in which $K$ is existentially closed. 
		Let $A, B$ be $p.p.$ formulas with one free variable  such that  
		$$T_1\models \forall x \mesp A_1(x) \rightarrow A(x).$$
		
		Set $k:=\vert A(K)/A_1(K) \vert$. It is enough to show that 
		we have 
		$$\vert A(N)/A_1(N) \vert =k.$$ 
		
		Let $l>k$ and suppose
		$$M \models \exists x_1,\ldots \exists x_l \mesp (\bigwedge_i x\in A)\wedge
		(\bigwedge_{i\neq j} x_i-x_j \notin A_1).$$ 
		Since modulo $T_{1}$ the formula $x \in A_1$ is equivalent to a universal formula, its negation is equivalent to an  existential $L_{\mathcal{O}}$-formula. Since
		$K$ is existentially closed in $N$, there exists at least $l$ element in $A(K)/A_1(K)$. Contradiction.

	\end{proof}

	In particular we have:
	\begin{corr}
		$K$ is the prime model of the complete theory of $\mathbb{F}_d((X))$.
	\end{corr}
	
	\begin{proof}
		$K$ is existentially closed in $\mathbb{F}_d((X))$ as an $L_{\mathcal{O}}$-structure since it is existentially closed in $\mathbb{F}_d((X))$ as a ring. Hence $$K\equiv \mathbb{F}_d((X)).$$ Since any completion of $T_1$ is model-complete $K$ is an $L_{\mathcal{O}}$ elementary substructure of $\mathbb{F}_d((X))$.
	\end{proof}
	
	\subsection*{Decidability of $\mathbb{F}_d((X))$}
	
	The following fact is an easy exercise using Hensel's lemma:
	\begin{fait}
		$\mathbb{F}_d[[X]]$ is definable by the $L$-$p.p.$ formula 
		$$ \exists y \mesp x.tX=y.(t-1)$$ inside $\mathbb{F}_d((X))$. 
	\end{fait}
	
	Hence, the decidability of $\mathbb{F}_d((X))$ as an 
	$L_\mathcal{O}$-structure and as an $L$-structure are equivalent. Recall that the $L$-theory of $K$ (hence of $\mathbb{F}_d((X)))$, are given by the sentences
	stating that 
	$$\vert A(K)/(A(K)\wedge B(K)) \vert =k \quad (k\in \mathbb{N}\cup \{\infty\})$$
	where $A, B$ are $L$-$p.p.$ formulas of  
	with one free variable.
	
	We will show that  a recursively enumerable subset of these sentences forms 
	a complete axiom system which implies all of them. Hence the $L$-theory of $K$ is decidable.

	For $A$ and $B$ as above, we set $D:=A\wedge B$. Let 
	$A$ be given by $\exists \bar{y} \mesp x.p=\bar{y}\bar{q}$ and 
	$B$ by $\exists \bar{y} \mesp x.r=\bar{y}\bar{s}$.  We may suppose that 
	both $p$ and $r$ are unitary. Hence, by Lemma \ref{lem:addtoMta}, for  $k=\deg(p)$ and $s=\deg(r)$  we have
	$$K.p \approx K.t^k$$ and $$K.r \approx K.t^s .$$

	
	Now by Remark \ref{rem:samecomplement}, \ref{rem:invimage} and Lemma \ref{changecomp}, $A\approx D$ if and only if  the preimage by $.t^k$ of $C_{\mathfrak{j}(\bar{q})}$(which is a pseudo-complement of $A$) is equal to
	the preimage  of $C_{\mathfrak{j}(\bar{s})}$ by $.t^s$ (which is a pseudo-complement of $D$). Hence it is decidable if $A\approx D$.
	
	Notice that  If  $A \not \approx D$ then
	$A/D$ is infinite:  Let $a\in A \setminus D$.  Then for  $a' \in A$ with $v(a')<v(a)$  either $a'$ or $a-a'$ is not in $D$. 
	
	Now suppose $A\approx D$. It follows by Remark \ref{rem:modulogamma} that
	$$ A+P_\alpha=D+P_\alpha $$  for some computable $\alpha$: $\alpha$ can be chosen less then every jump values of all polynomials appearing in the definition of $A$ and $B$. 
	Hence $\vert A/D \vert = \vert A\cap P_\alpha / D\cap P_\alpha\vert $.

	Now by Corollary \ref{corr:boule_up}, there is a computable $\gamma$ such that $A/D$ is finite if and only if
	$$A\cap P_\gamma = D \cap P_\gamma$$

	Consider the following algorithm: given $A$ and $B$  setting $D=A\wedge B$, the algorithm check if $A\approx B$, if not it  sets
	$\vert A/D \vert =\infty$. Otherwise computes $\alpha $ such that $A+P_\alpha=B+P_\alpha$.  Then it computes the value $\gamma$,
	so that it can check whether $A\cap P_{\gamma}/D\cap P_{\gamma}$  is trivial or infinite. If it is infinite it sets 
	$\vert A/D \vert = \infty$. If not, by Remark \ref{boundedsolutions} the algorithm can compute the number of the elements of 
	\begin{equation}\label{t0}
		\frac{(A\cap P_\alpha)/P_{\gamma}}{(D\cap P_\alpha)/P_{\gamma}}.
	\end{equation}
	
	Hence the $L$-theory consisting of sentences of type $(\ref{t0})$ can be
	recursively enumerable and implies the $L$-theory of $K$, hence of $\mathbb{F}_d((X))$. We have proved:
	\begin{theorem}\label{thm:finaldecidable}
		Both the $L$-  and the $L_{\mathcal{O}}$-theories of $\mathbb{F}_d((X))$ are decidable.
	\end{theorem}
	
	\bibliographystyle{amsalpha}
	\bibliography{onayg}
	
\end{document}